# Time-independent equivalent matrix representation for time-domain boundary element analysis of large-scale 3D acoustic problems

Weiyu Zhou[a], Yongqiang Chen[b,*]

[a] Institute of Nuclear and New Energy Technology, Tsinghua University, Beijing, 100084, China

[b] School of Mechanics and Engineering Science, Peking University, Beijing, 100871, China

**Abstract:** This paper proposes a time-independent equivalent matrix representation to enable large-scale three-dimensional (3D) acoustic analysis using the time-domain boundary element method (TDBEM). In conventional TDBEM, pair-dependent retarded times couple spatial interactions with the retrieval of historical boundary data, resulting in substantial storage and history-convolution costs. The proposed representation decouples the compressible spatial interaction structure from propagation history, allowing the same compressed spatial representation to be reused across propagation delays and throughout time marching, without additional temporal interpolation of the kernel. For conventional space-time discretizations, blockwise low-rank compression of the equivalent matrices leads to the compressed equivalent matrix method in the $\mathcal{H}$-matrix format (CEM-$\mathcal{H}$), with dynamic cluster-pair tracking used to avoid inactive interactions during wavefront propagation. For low-frequency, long-wavelength problems with only a bounded number of nonzero time-lag matrices, the framework is extended to the $\mathcal{H}^2$-matrix format (CEM-$\mathcal{H}^2$) using delay-consistent partitioning, shared coupling matrices, and nested cluster bases. Under the sparsity and bounded-parameter assumptions, the storage requirement and computational cost per time step of CEM-$\mathcal{H}^2$ both scale as $O(N_s)$, excluding preprocessing. Numerical results agree well with conventional TDBEM results and semi-analytical reference solutions. For CEM-$\mathcal{H}$, the storage requirement for a sphere with $N_s = 20{,}480$ is reduced from approximately $349.0\ \mathrm{GB}$ to $1.6\ \mathrm{GB}$, while a computational speedup of approximately $45.4\times$ is achieved. For CEM-$\mathcal{H}^2$, simulations with up to $N_s = 1{,}310{,}720$ are performed in the low-frequency, long-wavelength regime. A million-degree-of-freedom underwater-vehicle model is also successfully analyzed. These results demonstrate the capability of the proposed framework for large-scale 3D transient acoustic analysis using TDBEM.

# 1 Introduction

Transient wave propagation is encountered in a wide range of engineering applications. In acoustics, representative examples include pulse scattering, sound radiation, and acoustic-structure interaction. Unlike frequency-domain approaches, in which transient responses are reconstructed from solutions at multiple frequencies through an inverse transform, time-domain methods directly describe the evolution of wave fields in physical time, including their propagation, reflection, and decay, while naturally retaining excitation timing and propagation delays. For such problems, the time-domain boundary element method (TDBEM) requires only boundary discretization and inherently satisfies the radiation condition for exterior domains through the retarded fundamental solution. It therefore avoids truncation of the unbounded domain and the introduction of artificial absorbing boundaries, making it particularly attractive for transient radiation and scattering analyses [1-5].

However, the large number of boundary interactions and their history dependence impose substantial computational

* Corresponding author.

*E-mail address*: zwy97@tsinghua.edu.cn (Weiyu Zhou), chenyq@pku.edu.cn (Yongqiang Chen).

and storage costs on conventional TDBEM, limiting its application to large-scale transient acoustic problems [2, 6]. To address the computational bottleneck associated with spatial interactions, frequency-domain acoustic BEM has incorporated fast multipole methods and associated preconditioning techniques [7-11], as well as hierarchical low-rank compression methods [12, 13]. Model order reduction has also been used to reduce the cost of repeated solutions at multiple frequencies [14]. In contrast, accelerating direct time-domain computations additionally requires handling the coupling between propagation delays and access to historical boundary data. Various fast algorithms for transient wave fields have been developed to address these features [1, 2, 15-18].

The plane-wave time-domain (PWTD) algorithm accelerates long-range interactions through plane-wave representations of propagating fields, temporal grouping, and hierarchical translations; both two-level and multilevel variants have been developed [1, 16, 18]. The non-uniform grid time-domain (NGTD) algorithm reconstructs far fields from sparse-grid samples by exploiting their spatial smoothness after delay and amplitude compensation. Subsequent studies have developed multilevel variants and implementations for acoustics and acoustic–structure coupling [2, 15, 17]. For interactions at sub-wavelength scales, the accelerated Cartesian expansions (ACE) method combines retarded potentials with temporal basis functions and reformulates their evaluation as the fast summation of spatial potentials [19]. In addition, spherical Kirchhoff representations, spatial spectral reconstruction, and time-shifted interpolation have been used to accelerate the evaluation of transient acoustic fields generated by known sources [20]. These methods exploit different aspects of transient wave propagation but share a common feature: they construct space-time approximations of interaction kernels or propagation operators containing propagation-time information, thereby providing effective approaches to large-scale time-domain computations.

Importantly, in low-frequency or sub-wavelength problems, the size of the spatial discretization may still be dictated by fine geometric features, whereas temporal sampling is governed primarily by the frequency band of the response. In this regime, the interaction partitioning in conventional PWTD depends on the time-step size, and a large proportion of interactions may be assigned to direct evaluation, reducing its acceleration efficiency. Previous work has shown that, under the sub-wavelength conditions considered, the evaluation of such retarded potentials can be reformulated as the fast summation of spatial potentials, enabling potential evaluation with linear complexity per time step [19]. Nevertheless, these methods must still explicitly handle the coupling between temporal basis functions and retarded potentials within their propagation expansions.

Kernel-independent fast algorithms further reduce the reliance on analytical expansions of specific integral kernels. Approaches based on equivalent densities determined by field matching on check surfaces, or on kernel interpolation [21, 22], provide a relatively unified numerical framework for representing far-field interactions. In the time domain, the interpolation-based fast multipole method (IFMM) proposed by Takahashi interpolates the temporally integrated kernel in both space and time and employs a multilevel space-time hierarchy to accelerate retarded interactions [6]. Subsequent studies have extended the integral equations, temporal basis functions, and fast translation procedures [23]. Although auxiliary kernel interpolation reduces computational costs, it also introduces additional approximation errors. Takahashi [6] showed that, with fixed numbers of spatial and temporal interpolation points, the error may approach a plateau as the underlying discretization is refined; further improvements in accuracy require a higher interpolation resolution. In addition, temporal grouping is constrained by the spatial box size, and an excessively large time step can invalidate the grouping for a fixed spatial hierarchy.

In addition to propagation expansions and kernel interpolation, hierarchical matrices offer a compression approach based on spatial low-rank structure. Combining $\mathcal{H}$-matrices with adaptive cross approximation (ACA) enables low-rank representations of far-field blocks to be constructed from a small number of matrix entries. Under appropriate conditions, the storage and computational costs are reduced to $O(rN_s \log N_s)$ and $O(r^2 N_s \log N_s)$, respectively [24-27], where $r$ denotes the maximum rank of the admissible low-rank blocks. $\mathcal{H}^2$-matrices further share spatial representations across blocks and levels through nested cluster bases, reducing the corresponding costs to $O(N_s)$ [28-32]. These techniques have been applied to both iterative and direct solvers for acoustic BEM [12, 13] and extended to time-dependent boundary integral equations. For example, Seibel combined spatial $\mathcal{H}^2$-matrix representations with ACA in the frequency direction within a convolution quadrature framework [33], whereas Aimi et al. applied partially pivoted ACA to the corresponding temporal blocks in an energetic time-domain BEM formulation [34]. In these studies, compression is performed along the frequency direction or on temporal blocks, and the matrices to be compressed remain organized as frequency slices or temporal blocks. In contrast to these joint space-time compression strategies, the present work targets a time-marching formulation and separates propagation delays from the spatial coefficient matrices, yielding a set of equivalent matrices independent of the delay index. Combined with retarded boundary values, this representation describes the discrete action of the original family of time-lag matrices during time marching. Delay information is carried by the boundary values and index tables rather than being included in the spatial compression.

In conventional TDBEM, each coefficient matrix depends only on the difference between time-step indices and can therefore be reused as the current time step advances. This structure is rooted in the finite propagation speed of waves: boundary values are evaluated at retarded times, and the distance between a collocation point and a source element varies with their indices, so each spatial interaction has its own propagation delay. Consequently, coefficient matrices associated with different time lags must be stored separately even though they repeatedly contain similar spatial interaction information. Although the spatial interactions possess an inherent low-rank structure, the delay indices partition this compressible structure, hindering the sharing of its compressed representation across propagation delays and time steps. Existing fast methods generally reduce computational costs using propagation expansions or space-time kernel approximations, but their fast representations must still account for the propagation-time structure. To address this difficulty, we decouple propagation delays from compressible spatial interactions and construct time-independent spatial equivalent matrices that can be reused for different propagation delays, without introducing additional temporal interpolation of the kernel. The same compressed spatial representation can thus be reused throughout time marching.

Firstly, for conventional space-time discretizations, we propose a time-independent equivalent matrix compression method based on $\mathcal{H}$-matrices, termed CEM-$\mathcal{H}$ (compressed equivalent matrix in the $\mathcal{H}$-matrix format). The method retains the conventional treatment of propagation delays at quadrature points in the near field and introduces an element-centroid-based delay approximation only for temporal sampling in the far field. Spatial integration is still carried out at the actual quadrature points. The far-field contribution is thereby expressed in terms of three time-independent equivalent matrices and the corresponding retarded boundary values, and partially pivoted adaptive cross approximation (PACA) is used to compress the spatial matrices. Under standard hierarchical-matrix assumptions, the storage requirement for the far-field coefficients is $O(rN_s \log N_s)$. When the errors introduced by the centroid-based delay approximation and low-rank compression are sufficiently small relative to the original discretization error, the method can retain essentially the same accuracy as conventional TDBEM. In addition, an optional local dynamic cluster-pair strategy updates the initial wavefront arrival states to reduce unnecessary accesses to regions not yet reached by the wavefront.

For low-frequency, long-wavelength problems in which the wavelength is large relative to the computational domain, we propose an equivalent matrix compression method based on $\mathcal{H}^2$-matrices, termed CEM-$\mathcal{H}^2$ (compressed equivalent matrix in the $\mathcal{H}^2$-matrix format). Provided that the temporal sampling and discretization accuracy requirements are satisfied, a relatively large time step can make the distance traveled within one time step comparable to the maximum propagation distance across the computational domain, so that only a small number of time-lag matrices are nonzero [19, 35]. Through delay-consistent partitioning, the method retains quadrature-point propagation delays within cluster pairs satisfying the required conditions, combines the original piecewise-linear temporal weights with the spatial kernel, and reorganizes the treatment of propagation delays into blockwise matrix operations. On this basis, simplification of the double-layer matrix terms, sharing of coupling matrices between single- and double-layer contributions, and the use of nested cluster bases reduce the number of independent spatial operators and repeated transformations. The method accelerates both the evaluation of history terms and the application of the current-step operator within a unified framework. Under the sparsity and bounded-parameter assumptions stated in Section 5.4, it requires $O(N_\mathrm{s})$ storage and $O(N_\mathrm{s})$ computational cost per time step, excluding preprocessing.

The main contribution of this work is a time-independent equivalent matrix representation that reformulates the far-field retarded potentials in direct time-marching TDBEM using three spatial equivalent matrices, together with the corresponding retarded boundary quantities. Unlike conventional time-lag-dependent coefficient matrices, these equivalent matrices can be compressed independently of propagation delays, allowing the same spatial representation to be reused throughout time marching. However, the pair-dependent retarded boundary quantities prevent the resulting history interactions from being directly evaluated as standard matrix-vector products. To address this challenge, two complementary hierarchical implementations are developed. CEM-$\mathcal{H}$ reduces coefficient storage through PACA compression while retaining pairwise access to historical data, whereas CEM-$\mathcal{H}^2$ further enables hierarchical matrix-vector evaluation through delay-consistent partitioning and nested cluster bases in the low-frequency, long-wavelength regime. The proposed framework thus provides a way to integrate reusable spatial hierarchical compression with retarded boundary interactions without additional temporal interpolation of the kernel.

The remainder of this paper is organized as follows. The basic formulation of the three-dimensional (3D) acoustic TDBEM is presented in Section 2. The far-field time-independent equivalent matrices are constructed in Section 3, together with an analysis of the error introduced by the delay approximation. The CEM-$\mathcal{H}$ method and the local dynamic cluster-pair tree strategy are described in Section 4, followed by the CEM-$\mathcal{H}^2$ method in Section 5. Accuracy, storage requirements, and computational performance are evaluated through numerical examples in Section 6, and the main conclusions are summarized in Section 7.

# 2 Time-domain boundary element method for 3D acoustics

This section presents the basic formulation and discretization of conventional 3D acoustic TDBEM. The resulting convolution involving time-lag matrices provides the theoretical basis for the equivalent matrix representation and compression algorithms developed in the following sections.

## 2.1 Boundary integral equation

Consider a 3D acoustic domain $\Omega \subset \mathbb{R}^3$ with boundary $\Gamma = \partial\Omega$. The medium is assumed to be homogeneous and isotropic, with density $\rho$ and sound speed $c$. The acoustic pressure $p(\boldsymbol{x}, t)$ satisfies

$$\nabla^2 p(\boldsymbol{x}, t) - \frac{1}{c^2}\frac{\partial^2 p(x,t)}{\partial t^2} = 0,\ \boldsymbol{x} \in \Omega,\ t > 0 \tag{1}$$

where $t$ denotes time. Zero initial conditions are assumed, and the acoustic field is generated by boundary excitation applied for $t > 0$. On $\Gamma$, the outward normal derivative of the acoustic pressure is defined as $q(\boldsymbol{x}, t) = \partial p(\boldsymbol{x}, t)/\partial \boldsymbol{n}$, where $\boldsymbol{n}$ is the unit normal to $\Gamma$ pointing outward from $\Omega$.

Applying Green's identity in the time domain yields the following time-domain boundary integral equation:

$$\frac{1}{2}p(\boldsymbol{x}, t) + \int_\Gamma \int_0^t \frac{\partial G}{\partial \boldsymbol{n}_y}(\boldsymbol{x}, t; \boldsymbol{y}, \tau) p(\boldsymbol{y}, \tau)\, \mathrm{d}\tau \mathrm{d}\Gamma_y = \int_\Gamma \int_0^t G(\boldsymbol{x}, t; \boldsymbol{y}, \tau) q(\boldsymbol{y}, \tau)\, \mathrm{d}\tau \mathrm{d}\Gamma_y \tag{2}$$

where $G(\boldsymbol{x}, t; \boldsymbol{y}, \tau)$ is the fundamental solution of the 3D wave equation,

$$G(\boldsymbol{x}, t; \boldsymbol{y}, \tau) = \frac{\delta(t-\tau-r/c)}{4\pi} \tag{3}$$

where $\delta(\cdot)$ is the Dirac delta function, and $r = \| \boldsymbol{x} - \boldsymbol{y} \|$ is the distance between the source point $\boldsymbol{y}$ and the field point $\boldsymbol{x}$. The variables $t$ and $\tau$ denote the observation time at $\boldsymbol{x}$ and the source time at which the excitation is generated at $\boldsymbol{y}$, respectively, while $r/c$ is the propagation time from $\boldsymbol{y}$ to $\boldsymbol{x}$. The time delay in the fundamental solution ensures that the boundary integral representation inherently accounts for finite propagation speed and causality. For exterior scattering problems, the prescribed incident acoustic field is added to the right-hand side of Eq. (2) as a known excitation. This term does not affect the subsequent derivation of the boundary integral operators or the history contributions.

The time-domain boundary integral equation involves only the boundary acoustic pressure $p$ and its normal derivative $q$, so spatial discretization is required only on the boundary. Since the boundary quantities are evaluated at the retarded time $t - r/c$, different propagation distances between collocation points and source points require boundary data from different past times. Thus, while reducing the spatial dimension of the problem, 3D acoustic TDBEM also introduces coupling between spatial integration and the time history.

## 2.2 Temporal and spatial discretization

The collocation method is adopted to solve the time-domain boundary integral equation in Eq. (2). The boundary is discretized into $N_s$ three-node planar triangular elements $\Gamma_j$, $1 \le j \le N_s$, with the element centroids $\boldsymbol{x}_j$ serving as collocation points. The acoustic pressure $p$ and its normal derivative $q$ are approximated as constant over each element and interpolated using piecewise-linear functions between uniformly spaced time nodes $t_n = n\Delta t$. Denoting either boundary quantity by $\phi$, we write

$$\phi(\boldsymbol{y}, \tau) \approx \sum_{n=0}^{\infty} \phi_j^n \psi_n(\tau), \ \ \boldsymbol{y} \in \Gamma_j, \ \ \phi \in \{p, q\} \tag{4}$$

where $\phi_j^n$ is the discrete value of $\phi$ on element $\Gamma_j$ at time $t_n$, and $\psi_n(\tau) = \max\left(1 - \frac{|\tau - t_n|}{\Delta t}, 0\right)$ is the piecewise-linear temporal basis function centered at $t_n$. Accordingly, the boundary quantity at any time is obtained by linear interpolation between its discrete values at the two endpoints of the interval containing that time.

Although the boundary quantities are approximated as constant within each element, the propagation delay is still determined by the distance from each quadrature point to the collocation point. The spatial approximation of the boundary quantities introduces no further approximation to the propagation delay. Consequently, different quadrature points within the same element still correspond to different retarded times. For a collocation point $\boldsymbol{x}_i$ and a quadrature point $\boldsymbol{y}_{jg}$ on element $\Gamma_j$, let $r_{ijg} = \|\boldsymbol{x}_i - \boldsymbol{y}_{jg}\|$. At the current time $t_n$, the retarded time is $t_n - r_{ijg}/c$, and the corresponding boundary quantity is obtained by evaluating Eq. (4) at this time. Owing to the local support of the piecewise-linear basis functions, this evaluation involves only the historical boundary values at the two endpoints of the interval containing the retarded time.

Substituting Eq. (4) into Eq. (2) and evaluating the contributions for all collocation points and quadrature points yields the time-domain single-layer matrix $\boldsymbol{S}^\ell$ and double-layer matrix $\boldsymbol{D}^\ell$, whose entries are given by

$$S_{ij}^\ell = \int_{\Gamma_j} \int_0^{t_n} G(\boldsymbol{x}_i, t_n; \boldsymbol{y}, \tau) \psi_{n-\ell}(\tau) \, \mathrm{d}\tau \, \mathrm{d}\Gamma_y \tag{5}$$

$$D_{ij}^\ell = \int_{\Gamma_j} \int_0^{t_n} \frac{\partial G(\boldsymbol{x}_i, t_n; \boldsymbol{y}, \tau)}{\partial n_y} \psi_{n-\ell}(\tau) \, \mathrm{d}\tau \, \mathrm{d}\Gamma_y \tag{6}$$

The Dirac delta function in the fundamental solution localizes the temporal integration to the retarded time associated with wave propagation. For uniform time steps, the matrix coefficients depend only on the difference $\ell$ between time-step indices.

Nonsingular spatial integrals are evaluated by uniformly subdividing each triangular boundary element into 64 subtriangles and applying a symmetric three-point Gaussian quadrature rule [36] to each subtriangle, giving 192 quadrature points per element. For weakly singular integrals, each triangular element is divided into three subtriangles sharing the collocation point, and a $12 \times 12$ tensor-product Gauss-Legendre rule is applied to each subtriangle after a Duffy transformation. The singular integrals associated with the double-layer operator are treated analytically, with the Cauchy principal-value contribution vanishing for planar triangular elements and the jump term accounted for separately.

After spatial and temporal discretization, the 3D acoustic boundary element equation can be written as

$$\frac{1}{2}\boldsymbol{I}\boldsymbol{p}^{n} + \sum_{\ell=0}^{\min(n,\gamma_{\max})}\left(\boldsymbol{D}^{\ell}\boldsymbol{p}^{n-\ell} - \boldsymbol{S}^{\ell}\boldsymbol{q}^{n-\ell}\right) = \boldsymbol{0} \tag{7}$$

where $\boldsymbol{p}^{n-\ell}$ and $\boldsymbol{q}^{n-\ell}$ contain the acoustic pressure and its normal derivative, respectively, on all elements at time step $n-\ell$, and are referred to as historical boundary data. The maximum time-lag index, $\gamma_{\max}$, is determined by the maximum propagation distance and the time-step size:

$$\gamma_{\max} = \left\lceil \frac{\max_{\boldsymbol{x},\boldsymbol{y}\in\Gamma}\|\boldsymbol{x}-\boldsymbol{y}\|}{c\Delta t} \right\rceil \tag{8}$$

where $\max_{\boldsymbol{x},\boldsymbol{y}\in\Gamma}\|\boldsymbol{x}-\boldsymbol{y}\|$ is the maximum propagation distance, and $\lceil\cdot\rceil$ denotes the ceiling function.

The temporal discretization scale in TDBEM is generally related to the spatial discretization scale. To characterize the distance traveled by the wave within one time step relative to the characteristic boundary element size, we introduce the time-step parameter [37, 38]

$$\beta = \frac{c\Delta t}{h} = \frac{N_{\lambda}}{N_{T}} \tag{9}$$

where $h$ is the characteristic boundary element size, $N_{\lambda} = \lambda/h$ is the average number of boundary elements per wavelength, and $N_{T} = T_{0}/\Delta t$ is the average number of time steps per period. The quantities $\lambda$ and $T_0$ denote the wavelength and period, respectively. Similar parameters have been widely used in time-domain BEM to examine the effects of the time-step size on stability and computational accuracy. Previous studies have shown that the relation between the time-step size and the spatial mesh size significantly affects the stability and accuracy of acoustic TDBEM. For some classical collocation-based marching-on-in-time (MOT) schemes, an excessively small $\beta$ may induce late-time instability, whereas an excessively large time step reduces temporal resolution [2, 39].

Eq. (9) shows that the time-step parameter $\beta$ essentially relates the spatial and temporal discretization resolutions. In addition to measuring the relation between the time-step size and the spatial mesh size, $\beta$ can characterize the ability of the discrete system to resolve wave motion at different frequencies in time. This relation will be used in the following sections to discuss temporal interpolation accuracy and the applicability of hierarchical fast methods under different frequency conditions.

### 2.3 Discrete time-domain equations and computational characteristics of history terms

According to the boundary condition type, the boundary variables at the current time step are rearranged into an unknown vector $\boldsymbol{x}^{n}$ and a known vector $\boldsymbol{y}^{n}$, yielding the discrete equation

$$\boldsymbol{A}_{0}\boldsymbol{x}^{n} = \boldsymbol{b}^{n} = \boldsymbol{B}_{0}\boldsymbol{y}^{n} + \sum_{\ell=1}^{\min(n-1,\gamma_{\max})}\left(\boldsymbol{S}^{\ell}\boldsymbol{q}^{n-\ell} - \boldsymbol{D}^{\ell}\boldsymbol{p}^{n-\ell}\right) \tag{10}$$

where $\boldsymbol{A}_0$ and $\boldsymbol{B}_0$ are the coefficient matrices associated with the unknown and known boundary quantities at the current time step, respectively. These matrices can be assembled during preprocessing and reused throughout subsequent time marching. At each time step, $\boldsymbol{b}^{n}$ is formed from the current known boundary data and the historical boundary data. After Eq. (10) is solved, $\boldsymbol{p}^{n}$ and $\boldsymbol{q}^{n}$ are updated.

The history term on the right-hand side of Eq. (10) is a discrete convolution involving a sequence of spatial interaction matrices indexed by time lag and the corresponding historical boundary data. Conventional TDBEM typically precomputes and stores the coefficient matrices according to propagation delay and repeatedly accesses historical boundary data during time marching. Consequently, the costs of matrix storage and history-convolution evaluation increase rapidly with the number of spatial degrees of freedom and the number of propagation delays.

For well-separated far-field cluster pairs, the purely spatial kernels generally exhibit favorable compressibility and can be represented in low-rank form using methods such as hierarchical matrices. However, the time-lag matrices obtained after temporal discretization also contain piecewise structure arising from propagation delays and the local support of the temporal basis functions. Different collocation point-element pairs therefore often require boundary data from different past times. Compressing each time-lag matrix independently makes it difficult to fully exploit the compressibility of the spatial interactions and limits the reuse of compressed information across propagation delays and time steps.

To address these issues, the present work separates spatial interactions from propagation delays. Conventional integration and pointwise propagation delays are retained in the near field, whereas time-independent spatial interactions shared across different delays are extracted in the far field, with the corresponding propagation delays organized at the block or cluster level. In Section 4, $\mathcal{H}$-matrices are combined with PACA to compress far-field interactions and accelerate history-convolution evaluation. In Section 5, the nested-basis structure of $\mathcal{H}^2$-matrices is further exploited to reuse spatial information across levels.

# 3 Far-field time-independent equivalent matrices and hierarchical structure

The conventional space-time discretization for 3D acoustics presented in Section 2 incorporates spatial integration, propagation delays, and temporal interpolation into the coefficient matrices. To extract reusable spatial interactions, this section first establishes a geometric hierarchy for the matrix indices and then introduces two representations: directly integrated coefficients and far-field time-independent equivalent matrices. The separation of the resulting contributions into current-step terms and history terms is also described. Building on this framework, coefficient organization strategies and fast evaluation methods based on $\mathcal{H}$- and $\mathcal{H}^2$-matrices are developed in Sections 4 and 5, respectively.

## 3.1 Spatial cluster tree and geometric near-field and far-field partitioning

Let $I = \{1, \ldots, N_s\}$ and $J = \{1, \ldots, N_s\}$ denote the row index set associated with the collocation equations and the column index set associated with the boundary elements, respectively. Here, $i \in I$ corresponds to the $i$-th collocation point $\boldsymbol{x}_i$, and $j \in J$ corresponds to the $j$-th boundary element $\Gamma_j$. A matrix index pair $(i, j) \in I \times J$ represents the contribution of boundary element $\Gamma_j$ to the boundary integral equation at collocation point $\boldsymbol{x}_i$. Throughout this paper, a given combination of $\boldsymbol{x}_i$ and $\Gamma_j$ is referred to as a collocation point-element pair. The complete matrix index set is partitioned into the near-field index set $\mathcal{P}_\mathrm{n}$ and the far-field index set $\mathcal{P}_\mathrm{f}$, satisfying $\mathcal{P}_\mathrm{n} \cap \mathcal{P}_\mathrm{f} = \emptyset$, $\mathcal{P}_\mathrm{n} \cup \mathcal{P}_\mathrm{f} = I \times J$. These two sets form a disjoint partition of the complete matrix index space. The geometric near field is treated by direct integration, whereas the geometric far field serves as the candidate region for hierarchical representation and compression. The specific computational and storage strategies are presented in Sections 4 and 5.

**Cluster tree and block cluster tree.** Since the $i$-th collocation point is located at the centroid of the $i$-th triangular element, the row and column indices can be organized using the same cluster tree $\mathcal{T}$, whose root cluster is denoted by $C_0$. For any cluster $C \in \mathcal{T}$, let $\mathcal{K}_C$ denote its element index set. If $C$ is a non-leaf cluster, the index sets of its child clusters are pairwise disjoint and together form $\mathcal{K}_C$. The elements are recursively partitioned using the geometric bounding box $B_C$ of each cluster until every leaf cluster contains at most $n_\mathrm{leaf}$ elements. A global permutation $\boldsymbol{\pi}$ is then constructed according to a fixed traversal order, such that each cluster $C$ corresponds to a contiguous index range in the permuted ordering:

$$\mathcal{K}_C = \{\boldsymbol{\pi}(b_C), \ldots, \boldsymbol{\pi}(b_C + n_C - 1)\} \tag{11}$$

where $b_C$ and $n_C$ denote the starting position of the index range and the number of elements, respectively. This permutation preserves the ranks of matrix blocks and facilitates subsequent hierarchical partitioning and contiguous storage of low-rank data. Based on the cluster tree $\mathcal{T}$, each block $b$ is specified by an ordered pair of clusters $(C_c^b, C_s^b)$, where $C_c^b$ and $C_s^b$ denote the collocation-point cluster and source-element cluster, hereafter referred to as the row and column clusters, respectively. These cluster pairs form a block cluster tree according to the parent-child relationships in the cluster tree. The root block $b_0$, specified by $(C_0, C_0)$, corresponds to the complete matrix.

**Geometric admissibility.** Let $B_{C_c^b}$ and $B_{C_s^b}$ denote the geometric bounding boxes of the row and column clusters, respectively. The following geometric admissibility condition is adopted [24-26]:

$$\mathcal{M}_\sigma \left\{\operatorname{diam}\left(B_{C_c^b}\right), \operatorname{diam}\left(B_{C_s^b}\right)\right\} \leq \eta_{\text{adm}} \operatorname{dist}\left(B_{C_c^b}, B_{C_s^b}\right) \tag{12}$$

where $\operatorname{diam}(B)$ is the Euclidean diameter of bounding box $B$, $\operatorname{dist}(B_1, B_2)$ is the minimum Euclidean distance between two closed bounding boxes, and $\eta_{\text{adm}} > 0$ is the admissibility constant. The operator $\mathcal{M}_\sigma\{a_1, a_2\}$ selects a cluster size, with $\mathcal{M}_-\{a_1, a_2\} = \min\{a_1, a_2\}$ and $\mathcal{M}_+\{a_1, a_2\} = \max\{a_1, a_2\}$. A matrix block satisfying this condition is termed a geometrically admissible block. Sections 4 and 5 select the cluster partitioning strategy and the type of admissibility condition according to the spatial hierarchy employed. Thus, the two methods follow the same geometric partitioning principles but need not produce identical block structures or near-field and far-field index sets. The geometric admissibility condition partitions the matrix index space into geometric near-field and far-field regions. Blocks satisfying the condition are treated as far-field candidates. Blocks that do not satisfy it are subdivided further until both clusters are spatial leaves; terminal blocks that remain inadmissible are assigned to the geometric near field.

The block cluster tree is recursively constructed from the root block $b_0$, corresponding to the cluster pair $(C_0, C_0)$. If the current block satisfies the geometric admissibility condition, it is classified as a geometric far-field block, and subdivision stops. If the condition is not satisfied and both the row and column clusters are spatial leaves, the block is classified as a geometric near-field block. Otherwise, the clusters that can still be subdivided are partitioned according to the chosen strategy, and the same tests are applied recursively to the child blocks.

In the near field, the spatial kernels vary rapidly, and variations in propagation delay within an element may affect the evaluation of transient interactions. The direct integration approach described in Section 2.2 is therefore retained, preserving the propagation delay at each quadrature point. For regions classified geometrically as far field, the choice between a compressed representation and direct evaluation additionally depends on the relevant temporal conditions. The corresponding treatments are described in Sections 4 and 5.

### 3.2 Directly integrated near-field coefficients

For near-field matrix indices $(i, j) \in \mathcal{P}_\text{n}$, the conventional space-time discretization described in Section 2 is used. The propagation delay is determined from the actual distance between each quadrature point and the corresponding collocation point, hereafter referred to as the quadrature-point delay. The combined contributions of spatial integration and the temporal basis functions are accumulated into the corresponding matrix coefficients. The near-field part $\boldsymbol{S}_\text{n}^\ell$ of the conventional time-domain single-layer matrix $\boldsymbol{S}^\ell$ can be expressed as

$$\left(S_\text{n}^\ell\right)_{ij} = \begin{cases} S_{ij}^\ell, & (i,j) \in \mathcal{P}_\text{n} \\ 0, & (i,j) \notin \mathcal{P}_\text{n} \end{cases} \tag{13}$$

The near-field double-layer matrix $\boldsymbol{D}_\text{n}^\ell$ is obtained analogously. During matrix assembly, the pairs $(i, j) \in \mathcal{P}_\text{n}$ can be traversed directly, and the contribution from each quadrature point is assigned to the appropriate matrix coefficient according to its time-lag index. It is therefore unnecessary to construct the complete matrices first. A block-level compressed sparse row (CSR) format is used for sparse storage [37, 40].

### 3.3 Far-field time-independent equivalent matrices

For far-field matrix entries, the collocation point and the integration element are sufficiently separated in space, and the integral kernel generally varies smoothly within the element, providing favorable conditions for subsequent matrix compression. However, in conventional coupled space-time integration, the propagation delay varies with the integration position within the element. An appropriate approximation of the propagation distance in the time-dependent terms is therefore introduced so that the delay associated with a given collocation point-element pair remains constant during integration over the element. This allows the time-dependent terms to be separated from spatial integration.

For any far-field matrix index pair $(i,j) \in \mathcal{P}_{\mathrm{f}}$, let $\boldsymbol{y}_j^{\mathrm{cen}}$ denote the geometric centroid of boundary element $\Gamma_j$. The distances from the collocation point $\boldsymbol{x}_i$ to a quadrature point $\boldsymbol{y}_{jg} \in \Gamma_j$ and to the centroid $\boldsymbol{y}_j^{\mathrm{cen}}$ are $r_{ijg} = \|\boldsymbol{x}_i - \boldsymbol{y}_{jg}\|$ and $R_{ij} = \|\boldsymbol{x}_i - \boldsymbol{y}_j^{\mathrm{cen}}\|$, respectively. The far-field equivalent representation retains the dependence of the spatial kernel on the actual quadrature-point distance $r_{ijg}$. Only the propagation delay in the temporal argument is approximated uniformly using the centroid distance $R_{ij}$:

$$\phi_j\left(t - \frac{r_{ijg}}{c}\right) \approx \phi_j\left(t - \frac{R_{ij}}{c}\right), \ \phi \in \{p, q\} \tag{14}$$

The centroid distance in Eq. (14) is used only to determine the temporal sampling location; it does not replace the actual distance in the spatial kernel. Far-field spatial integration is therefore still performed at the quadrature points of the boundary element, rather than replacing the entire element with a lumped contribution located at its centroid. Once the propagation delay is determined uniformly using $R_{ij}$, hereafter referred to as the centroid-based delay, all quadrature points in element $\Gamma_j$ access boundary data at the same time $t - R_{ij}/c$. The corresponding boundary quantities are then independent of the quadrature point and can be taken outside the spatial integral.

Fig. 1 schematically illustrates the quadrature-point delay treatment and the centroid-based delay approximation on a triangulated boundary surface. The circles represent quadrature points distributed within each element, while the diamonds represent element centroids. Filled and open symbols indicate locations inside and outside the wavefront, respectively. In the quadrature-point delay treatment, propagation delays are evaluated at individual quadrature points, so different points within the same element may lie on opposite sides of the wavefront. In contrast, the centroid-based approximation assigns a common propagation delay to each element based on its centroid, while retaining the original quadrature points for spatial integration.

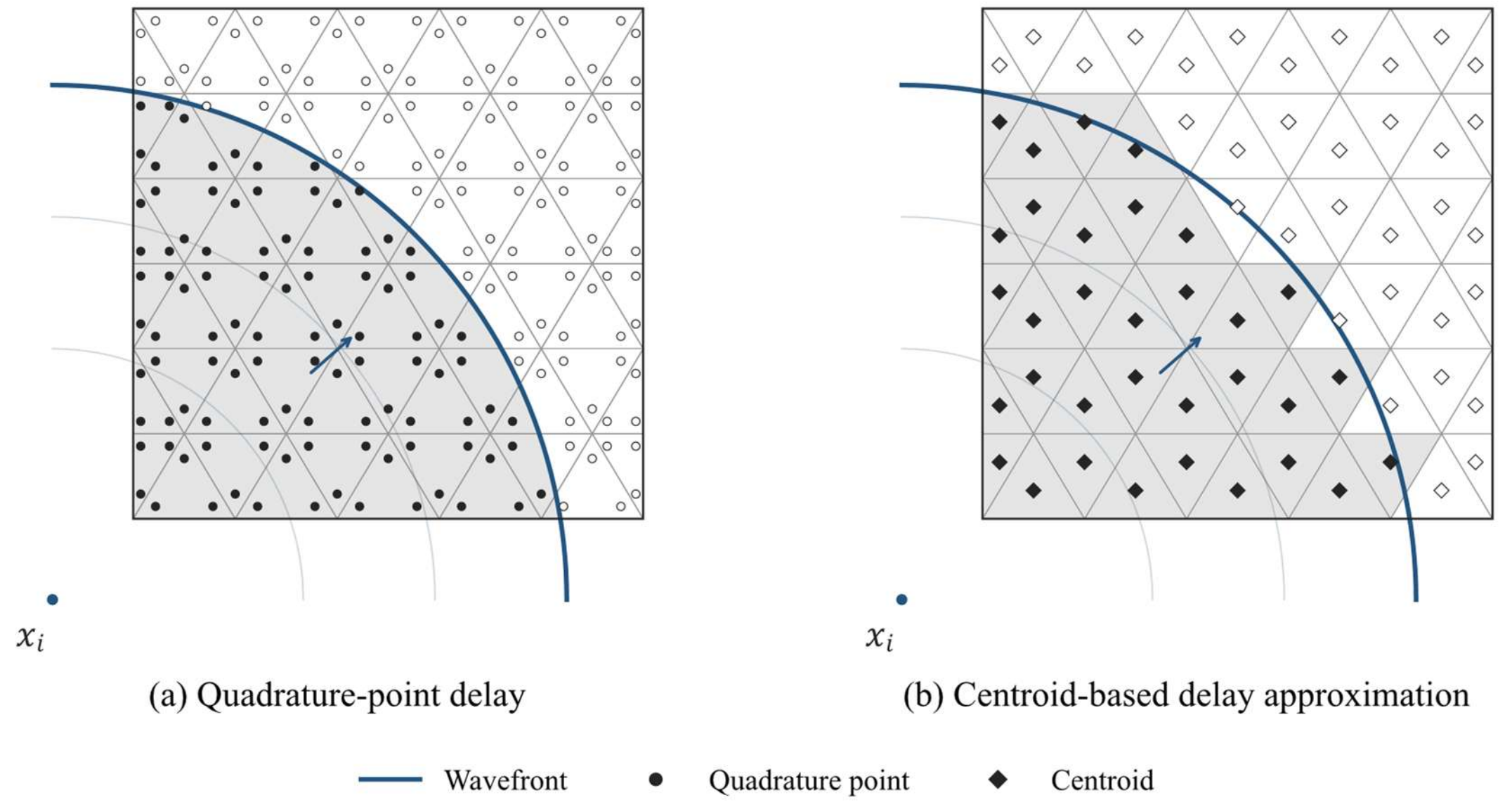

Fig. 1 Illustration of quadrature-point delays and the centroid-based delay approximation

To establish the spatial decomposition of the far-field integral terms, consider the causal fundamental solution $G(\boldsymbol{x}, t; \boldsymbol{y}, \tau)$ given in Eq. (3). Its spatial factor, $1/(4\pi r)$, corresponds to the single-layer integral term. The derivative of $G(\boldsymbol{x},t;\boldsymbol{y},\tau)$ with respect to the outward boundary normal at the source point $\boldsymbol{y}$ is

$$q^*(\boldsymbol{x},t;\boldsymbol{y},\tau) = \frac{\partial G}{\partial \boldsymbol{n}_y}(\boldsymbol{x},t;\boldsymbol{y},\tau) = \frac{(x-y)\cdot\boldsymbol{n}_y}{4\pi}\left[\frac{\delta(t-\tau-r/c)}{r^3} + \frac{\delta'(t-\tau-r/c)}{cr^2}\right] \tag{15}$$

where $\delta'(\cdot)$ denotes the derivative of the Dirac delta function with respect to its argument. After temporal convolution, the two terms in $q^*$ act on the retarded acoustic pressure $p(t-r/c)$ and its time derivative $\dot{p}(t-r/c)$, respectively.

We first consider the single-layer integral term. Substituting the fundamental solution $G(\boldsymbol{x},t;\boldsymbol{y},\tau)$ from Eq. (3) into the corresponding temporal convolution and setting $t_r = t - r/c$ gives

$$\int_0^t G\,(\boldsymbol{x},t;\boldsymbol{y},\tau) q(\boldsymbol{y},\tau)\, d\tau = \frac{1}{4\pi r}\int_0^t \delta\,(t_r-\tau) q(\boldsymbol{y},\tau)\, \mathrm{d}\tau \tag{16}$$

Using the sifting property of the Dirac delta function $\delta(\cdot)$, we obtain

$$\int_0^t \delta\,(t_r-\tau) q(\boldsymbol{y},\tau)\, d\tau = q(\boldsymbol{y},t_r) \tag{17}$$

Thus, the single-layer contribution of the $j$-th boundary element $\Gamma_j$ to the $i$-th collocation equation is

$$S_{ij}q_j = \int_{\Gamma_j} \frac{q_j(t-r_{ij}(y)/c)}{4\pi r_{ij}(y)} d\Gamma_y \tag{18}$$

For a far-field matrix index pair $(i,j) \in \mathcal{P}_\mathrm{f}$, the propagation delay is determined uniformly using the centroid distance $R_{ij}$. According to Eq. (14),

$$S_{ij}q_j \approx \int_{\Gamma_j} \frac{1}{4\pi r_{ij}(y)} q_j\left(t-\frac{R_{ij}}{c}\right) d\Gamma_y = \left[\int_{\Gamma_j}\frac{1}{4\pi r_{ij}(y)} d\Gamma_y\right] q_j\left(t-\frac{R_{ij}}{c}\right) \tag{19}$$

The retarded boundary quantity is thus separated from the spatial integral, whose coefficient depends only on the collocation-point position and boundary-element geometry. This leads to the time-independent single-layer equivalent matrix $\boldsymbol{S}$, with entries

$$S_{ij} = \int_{\Gamma_j}\frac{1}{4\pi r_{ij}(y)} d\Gamma_y \tag{20}$$

Similarly, for the double-layer integral term, temporal convolution of the two terms in Eq. (15) yields contributions involving the retarded acoustic pressure and its time derivative. Applying the same centroid-based delay approximation separates these boundary quantities from the corresponding spatial integrals, leading to two time-independent double-layer equivalent matrices, $\boldsymbol{D}_0$ and $\boldsymbol{D}_1$, whose entries are

$$D_{\mu,ij} = \int_{\Gamma_j}\frac{(x_i-y)\cdot\boldsymbol{n}_j}{4\pi r_{ij}^{3-\mu}(y)} d\Gamma_y \tag{21}$$

where $D_{\mu,ij}$, with $\mu = 0{,}1$, represents the acoustic pressure and pressure time-derivative terms in the double-layer potential, respectively. The detailed derivation is provided in Appendix A. Since $D_{\mu,ij}$ for $\mu = 0$ and $\mu = 1$ involves different powers of distance and acts on different boundary quantities, the two matrices must be constructed separately. The three equivalent matrices $\boldsymbol{S}$, $\boldsymbol{D}_0$, and $\boldsymbol{D}_1$ are independent of the temporal terms and describe only the boundary interactions determined by the geometry and spatial integral kernels. They contain no propagation-delay information and can be evaluated on demand during preprocessing and stored in the corresponding hierarchical formats. Unlike conventional time-lag matrix representations, an equivalent matrix is not associated with a particular time lag. Its compressed representation can be reused for different propagation delays and throughout all time steps. In particular, the definition of $D_{1,ij}$ does not include factors involving the sound speed $c$ or the time-step size $\Delta t$; these factors are incorporated when evaluating the history terms at each time step.

For far-field matrix indices, Eqs. (20) and (21) involve no singular integrals. Since all three matrices use the same quadrature points and geometric information, their entries can be evaluated together in a single numerical integration pass. For far-field indices using the centroid-based delay representation, the propagation time associated with $R_{ij}$ is expressed in units of time steps and decomposed into an integer delay $k_{ij}$ and a fractional delay $\alpha_{ij}$:

$$\frac{R_{ij}}{c\Delta t} = k_{ij} + \alpha_{ij} \tag{22}$$

where $k_{ij} = \left\lfloor \frac{R_{ij}}{c\Delta t} \right\rfloor$, and $0 \leq \alpha_{ij} < 1$, where $\lfloor \cdot \rfloor$ denotes the floor function.

Using the piecewise-linear temporal interpolation described in Section 2.2, the corresponding far-field boundary quantities are expressed as

$$\begin{aligned} \bar{\phi}_{ij}^n &= (1 - \alpha_{ij})\phi_j^{n-k_{ij}} + \alpha_{ij}\phi_j^{n-k_{ij}-1}, \phi \in \{p, q\} \\ \Delta p_{ij}^n &= p_j^{n-k_{ij}} - p_j^{n-k_{ij}-1} \end{aligned} \tag{23}$$

where $\bar{p}_{ij}^n$ and $\bar{q}_{ij}^n$ denote the acoustic pressure and its normal derivative obtained by linear interpolation of historical boundary data according to the propagation delay, while $\Delta p_{ij}^n$ is the pressure difference over the corresponding time interval. In particular, $\phi_j^{n-k_{ij}} = 0$ when $n - k_{ij} < 0$, ensuring consistency with the zero initial conditions and causality. For convenience, these three types of delay-dependent boundary quantities are collectively referred to as retarded boundary quantities. Unlike the quadrature-point delays used in Section 2.2, each far-field matrix index pair $(i, j) \in \mathcal{P}_\mathrm{f}$ requires only one integer delay $k_{ij}$ and one fractional delay $\alpha_{ij}$, both determined from the centroid distance $R_{ij}$.

Combining the three equivalent matrices with the far-field retarded boundary quantities, the far-field contribution $\boldsymbol{f}_\mathrm{f}^n$ at time step $n$ is expressed as

$$[\boldsymbol{f}_\mathrm{f}^n]_i = \sum\nolimits_{(i,j)\in\mathcal{P}_\mathrm{f}} \left[ S_{ij}\bar{q}_{ij}^n - D_{0,ij}\bar{p}_{ij}^n - \frac{D_{1,ij}}{c\Delta t}\Delta p_{ij}^n \right] \tag{24}$$

Eq. (24) shows that far-field spatial interactions are no longer represented by a sequence of matrices stored separately for different time lags $\ell$. Instead, they are evaluated using the three time-independent equivalent matrices $\boldsymbol{S}$, $\boldsymbol{D}_0$, and $\boldsymbol{D}_1$, together with the corresponding retarded boundary quantities. This separation yields time-independent far-field matrices that can be reused across time steps and compressed independently using low-rank approximations. It should be noted that, for a fixed mesh, sound speed, and time-step size, both $k_{ij}$ and $\alpha_{ij}$ are time-independent, but they depend on both the row and column indices. Consequently, Eq. (24) cannot generally be recast directly as a standard matrix–vector product. To evaluate this nonstandard history convolution, in Section 4, the time-independent equivalent matrices are compressed using PACA, and their contributions are evaluated with collocation-point-dependent retarded boundary quantities. In Section 5, low-frequency, long-wavelength problems with a small number of nonzero time-lag matrices are addressed by constructing an $\mathcal{H}^2$-matrix representation for fast evaluation of both history terms and current-step matrix contributions.

As shown in Eq. (23), a far-field contribution represented using the centroid-based delay involves only known historical boundary quantities when $k_{ij} \geq 1$. When $k_{ij} = 0$, it involves boundary quantities at both the current and previous time steps. According to the boundary conditions, the coefficients associated with the current unknown boundary quantities are combined with the corresponding current-step coefficients from direct near-field integration to form the coefficient matrix $\boldsymbol{A}_0$ for the unknowns.

### 3.4 Error analysis

This section analyzes the local error introduced by the propagation-delay approximation in regions where the centroid-based equivalent representation is used. In this representation, the quadrature-point distance is replaced by the centroid distance only in the propagation-delay term, so that a common delay is used to approximate the propagation times associated with all integration positions within the source element.

For a far-field matrix index pair $(i,j) \in \mathcal{P}_{\mathrm{f}}$, define $r_{ij}(\boldsymbol{y}) = \|\boldsymbol{x}_i - \boldsymbol{y}\|$, $\boldsymbol{y} \in \Gamma_j$. The maximum element-wise propagation-time deviation introduced by replacing the actual propagation distance with the centroid distance is defined as

$$\delta\tau_{ij} = \max_{\boldsymbol{y}\in\Gamma_j}\left|\frac{r_{ij}(\boldsymbol{y})}{c} - \frac{R_{ij}}{c}\right| = \frac{1}{c}\max_{\boldsymbol{y}\in\Gamma_j}\left|r_{ij}(\boldsymbol{y}) - R_{ij}\right| \tag{25}$$

Thus, $\delta\tau_{ij}$ is a theoretical element-wise measure and does not depend on the particular quadrature rule used to evaluate the spatial integral. For every quadrature point $\boldsymbol{y}_{jg} \in \Gamma_j$, its propagation-time deviation from the centroid-based delay is bounded by $\delta\tau_{ij}$. For a fixed quadrature rule, the realized quadrature-point maximum therefore cannot exceed the element-wise value defined in Eq. (25).

A dimensionless measure of the delay deviation relative to one time step is defined by $\eta_{ij}^{\tau} = \delta\tau_{ij}/\Delta t$. Suppose that the historical boundary quantity $\phi_j(t)$ on the $j$-th element is absolutely continuous over the relevant local retarded-time interval and has an essentially bounded first time derivative. Then, for any quadrature point $\boldsymbol{y}_{jg}$, the fundamental theorem of calculus and an $L^{\infty}$-norm estimate give

$$\left|\phi_j\left(t - \frac{r_{ijg}}{c}\right) - \phi_j\left(t - \frac{R_{ij}}{c}\right)\right| \le \delta\tau_{ij}\|\dot{\phi}_j\|_{L^{\infty}\left(I_{ij}(t)\right)} \tag{26}$$

where $I_{ij}(t)$ denotes the smallest time interval containing the centroid-based retarded time $t - R_{ij}/c$ and the retarded times $t - r_{ij}(\boldsymbol{y})/c$, with $\boldsymbol{y} \in \Gamma_j$. Eq. (26) shows that the local error in the retarded boundary quantity introduced by the centroid-based delay approximation is governed jointly by the element-wise propagation-time deviation $\delta\tau_{ij}$ and the temporal variation of the boundary quantity.

The geometrical meaning of $\delta\tau_{ij}$ can be further clarified by the reverse triangle inequality. If $\rho_j = \max_{\boldsymbol{y}\in\Gamma_j}\|\boldsymbol{y} - \boldsymbol{y}_j^{\mathrm{cen}}\|$, then $\delta\tau_{ij} \le \rho_j/c$. Hence, the propagation-time deviation is bounded by a quantity proportional to the element size and inversely proportional to the sound speed. The centroid-based approximation is therefore expected to be accurate when $\eta_{ij}^{\tau} \ll 1$ and the boundary quantities vary slowly over the corresponding retarded-time interval.

Eqs. (25) and (26) characterize only the local error introduced by the centroid-based delay approximation. They do not include spatial quadrature error, piecewise-linear temporal interpolation error, or error accumulation in the fully discrete time-marching procedure. A detailed derivation of the local error bound for the far-field element contribution is given in Appendix B.

In summary, this section establishes a geometric hierarchical partition of the matrix index space, presents the directly integrated coefficients and the three time-independent equivalent matrices, and describes the separation of current-step and history contributions. Building on these results, the following two sections develop fast methods based on $\mathcal{H}$- and $\mathcal{H}^2$-matrices, respectively. For conventional space-time discretizations, an $\mathcal{H}$-matrix representation of the equivalent matrices is constructed in Section 4 using PACA, together with an optional local dynamic cluster-pair tree strategy to reduce unnecessary accesses during initial wavefront propagation. For low-frequency, long-wavelength problems with a small number of nonzero time-lag matrices, the formulation in Section 5 combines delay-consistent partitioning with nested cluster bases to obtain an $\mathcal{H}^2$-matrix representation, thereby accelerating both history-convolution evaluation and the solution of the system of equations.

# 4 PACA compression based on $\mathcal{H}$-matrices

For the far-field time-independent equivalent matrices $\boldsymbol{S}$, $\boldsymbol{D}_0$, and $\boldsymbol{D}_1$ established in Section 3, this section develops an $\mathcal{H}$-matrix-based PACA compression method, termed CEM-$\mathcal{H}$. The method is intended for conventional space-time discretizations. PACA is used to compress the far-field equivalent matrices, while propagation delays are determined from the distances between collocation points and source-element centroids. The original temporal discretization is retained, with the centroid-based delay approximation and low-rank compression introduced only in the far field. The resulting computational accuracy remains essentially consistent with that of the conventional method.

Since the retarded boundary quantities depend on both the collocation-point and source-element indices, history convolution using the compressed representation must still account for pointwise propagation delays. To reduce unnecessary accesses during initial wavefront propagation, a local dynamic cluster-pair tree strategy is further proposed. This strategy updates the wavefront arrival states within each $\mathcal{H}$-matrix far-field block to identify the regions involved in the computation. It changes neither the original low-rank factors nor the delayed sampling rules and introduces no additional approximation error.

### 4.1 PACA compression of far-field equivalent matrices

A balanced binary cluster tree is constructed by median splitting along the longest axis. Far-field candidate blocks are selected using the min-type geometric admissibility condition $\mathcal{M}_\{a_1, a_2\}$ in Section 3.1, Eq. (12), with the admissibility constant $\eta_{\mathrm{adm}} = 0.6$. An additional far-field causal admissibility condition is introduced to ensure that the far-field contributions represented by the equivalent matrices depend only on known historical boundary quantities.

**Far-field causal admissibility.** For any block $b$, let $R_{\min,b}$ denote a conservative lower bound on the distances between collocation points and element centroids within the block, satisfying $R_{\min,b} \leq R_{ij},\ (i,j) \in I_b \times J_b$. A block $b$ is termed a causally admissible far-field block if $c\Delta t \leq R_{\min,b}$. In this case, $k_{ij} \geq 1$ holds for every collocation point–element pair within the block.

In the regime considered for CEM-$\mathcal{H}$, $\beta$ is of order unity, so $c\Delta t = \beta h$ is comparable to the characteristic boundary element size $h$. Under typical mesh and parameter settings, geometrically admissible far-field blocks also satisfy the far-field causal admissibility condition. Accordingly, the equivalent contributions of all far-field blocks are included entirely in the history term and do not contribute to the construction of the current-step system matrix $\boldsymbol{A}_0$.

For a far-field block $b$ satisfying both geometric and causal admissibility, let $I_b = \mathcal{K}_{C_c^b} \subseteq I,\ J_b = \mathcal{K}_{C_s^b} \subseteq J$ denote its row and column index sets, with $m_b = |I_b|$ rows and $n_b = |J_b|$ columns. For an equivalent matrix $\boldsymbol{K} \in \{\boldsymbol{S}, \boldsymbol{D}_0, \boldsymbol{D}_1\}$, denote the corresponding far-field matrix block by $\boldsymbol{A}_{K,b} \in \mathbb{R}^{m_b \times n_b}$. PACA approximates this block as [26, 34]

$$\boldsymbol{A}_{K,b} \approx \widetilde{\boldsymbol{A}}_{K,b} = \boldsymbol{U}_{K,b}\left(\boldsymbol{V}_{K,b}\right)^{\mathrm{T}} = \sum_{l=1}^{r_{K,b}} \boldsymbol{u}_{K,b,l}\left(\boldsymbol{v}_{K,b,l}\right)^{T} \tag{27}$$

where $r_{K,b}$ is the adaptively determined rank of far-field block $b$ of matrix $\boldsymbol{K}$, and $\boldsymbol{u}_{K,b,l} \in \mathbb{R}^{m_b}$ and $\boldsymbol{v}_{K,b,l} \in \mathbb{R}^{n_b}$ are the low-rank factor vectors associated with the $l$-th rank-one approximation term. Since the integral kernels of the three equivalent matrices have different characteristics, the matrices share the binary cluster tree, block structure, and index ordering, but their low-rank factors are computed independently.

The PACA variant adopted here alternates between rowwise and columnwise pivot selection. For notational simplicity, consider a far-field block $b$ and one matrix $\boldsymbol{K}$, and write $\boldsymbol{A} = \boldsymbol{A}_{K,b},\ m = m_b,\ n = n_b,\ \boldsymbol{u}_l = \boldsymbol{u}_{K,b,l}$, and $\boldsymbol{v}_l = \boldsymbol{v}_{K,b,l}$. At the $l$-th iteration, the contributions of the existing rank-one terms are first subtracted from the currently selected row. The entry with the largest absolute value in the resulting residual row determines the pivot column. The corresponding residual column is then computed, and the row containing its largest-magnitude entry is selected for the next iteration. After normalization by the pivot, new vectors $\boldsymbol{u}_l$ and $\boldsymbol{v}_l$ are formed, and the low-rank approximation is updated as $\widetilde{\boldsymbol{A}}^l = \widetilde{\boldsymbol{A}}^{l-1} + \boldsymbol{u}_l\boldsymbol{v}_l^T$. The iteration is terminated when the latest rank-one update satisfies

$$\|\boldsymbol{u}_l\|_2\|\boldsymbol{v}_l\|_2 \leq \varepsilon_{\mathrm{rel}}\left\|\widetilde{\boldsymbol{A}}^l\right\|_F \tag{28}$$

indicating that the update is sufficiently small relative to the current approximation. Here, $\varepsilon_{\mathrm{rel}}$ is the relative tolerance controlling the low-rank approximation accuracy. The Frobenius norm of the current approximation, $\left\|\widetilde{\boldsymbol{A}}^l\right\|_F$, can be evaluated from the low-rank factors already constructed. Because pivot searches are restricted to the current residual row and column, PACA avoids the full-matrix evaluation required by complete pivoting.

For $\boldsymbol{A} \in \mathbb{R}^{m\times n}$, a rank-$r$ approximation requires the factor matrices $\boldsymbol{U} \in \mathbb{R}^{m\times r}$ and $\boldsymbol{V} \in \mathbb{R}^{n\times r}$. The low-rank and dense representations require $r(m+n)$ and $mn$ double-precision numbers, respectively, to store their coefficients. Therefore, if the current approximation rank satisfies $r \geq r_{\text{stop}} = \lceil mn/(m+n) \rceil$, the low-rank representation requires at least as much storage as the dense representation. PACA is then terminated, and the block is stored in dense form.

For the same collocation point-element index pair $(i,j)$, the entries of the three equivalent matrices share the same integration domain, quadrature points, and geometric information, differing only in their integral kernels. These entries are therefore evaluated together when required, avoiding repeated quadrature evaluations. PACA is performed separately for each matrix, with independently determined pivot sequences and approximation ranks.

After PACA terminates, $\boldsymbol{U}_{K,b}$ and $\boldsymbol{V}_{K,b}$ are retained if the low-rank representation provides an actual storage benefit; otherwise, the block is stored in dense form. For blocks requiring a dense fallback, the remaining unevaluated entries are computed in a single joint traversal. These procedures do not alter the standard PACA iteration. Algorithm 1 summarizes the complete construction procedure for a single admissible far-field block $b$.

**Algorithm 1 PACA compression of equivalent matrices**

**Input**: Row and column index sets $I_b$ and $J_b$ of an admissible far-field block $b$; relative tolerance $\varepsilon_{\text{rel}}$; matrix-entry evaluation function $E(i,j) = (S_{ij}, D_{0,ij}, D_{1,ij})$.

Initialize a temporary row and column cache for block $b$.

for $\boldsymbol{K} \in \{\boldsymbol{S}_b, \boldsymbol{D}_{0,b}, \boldsymbol{D}_{1,b}\}$ do

  Evaluate the rows and columns selected by PACA on demand using $E(i,j)$.

  Reuse the previously evaluated entries of all three matrices stored in the cache.

  Construct $\boldsymbol{U}_{K,b}$ and $\boldsymbol{V}_{K,b}$ incrementally following the PACA pivoting strategy.

  Obtain the approximation rank $r_{K,b}$ when the stopping criterion is satisfied.

  if the low-rank factors provide a storage benefit then

    Store the matrix block in low-rank form (*LowRank*).

  else

    Mark the matrix block as *DenseFallback*.

  end if

end for

Evaluate the remaining row and column entries for matrix types requiring dense storage, and release the temporary cache.

**Output**: Storage objects for the equivalent matrices $\boldsymbol{S}$, $\boldsymbol{D}_0$, and $\boldsymbol{D}_1$ on far-field block $b$ after PACA compression.

4.2 Dynamic cluster-pair state classification and updating within $\mathcal{H}$-matrix blocks

The far-field causal admissibility condition in Section 4.1 ensures that the equivalent far-field blocks contain no current-step unknowns. Nevertheless, their active computational regions at the current time step are still determined by propagation delays. Here, the initial wavefront refers to the causal arrival front associated with the beginning of the boundary history at $t = 0$, rather than a physical wavefront of nonzero amplitude at that instant. It indicates which source–collocation interactions may contribute to the history convolution as the elapsed propagation time increases. Before this wavefront has covered the entire far field, checking every collocation point–element pair at every time step leads to repeated distance calculations, arrival tests, and accesses to historical data. An optional local dynamic cluster-pair tree $\mathcal{T}_{\text{d}}$ strategy is therefore introduced.

For a far-field block $b$, the root cluster pair of its local dynamic tree is defined as $\nu_b = (C_c^b, C_s^b)$. Starting from $\nu_b$, the row or column cluster is subdivided using the spatial binary cluster tree constructed during preprocessing, yielding a local dynamic cluster-pair tree $\mathcal{T}_b$. An arbitrary cluster pair is denoted by $\nu = (C_c^\nu, C_s^\nu)$. Subdivision and merging within each local tree are restricted to its associated far-field block and do not change the original matrix partition or low-rank factors.

The geometric bounding boxes $B_{C_c^\nu}$ and $B_{C_s^\nu}$ of the row and column clusters provide conservative lower and upper bounds, $R_{\min,\nu}$ and $R_{\max,\nu}$, on the propagation distances associated with cluster pair $\nu$. The corresponding integer-delay interval is defined as $[k_{\min,\nu}, k_{\max,\nu}]$, where $k_{\min,\nu} = \lfloor R_{\min,\nu}/(c\Delta t) \rfloor$, $k_{\max,\nu} = \min\{\lfloor R_{\max,\nu}/(c\Delta t) \rfloor, \gamma_{\max} - 1\}$. The maximum time-lag index, $\gamma_{\max}$, is incorporated into $k_{\max,\nu}$ to prevent the conservative distance bound $R_{\max,\nu}$ from extending the delay interval beyond its valid range. For the root pair $\nu_b$, the corresponding integer-delay interval is $[k_{\min,b}, k_{\max,b}]$.

At time step $n$, each cluster pair is classified into one of three states according to its integer-delay interval: *Unreached*, *Partially Arrived*, or *Fully Arrived*, indicating whether none, some, or all of its interactions have become causally reachable. The classification criteria and associated treatments are summarized in Tab. 1.

Tab. 1 State classification of dynamic cluster pairs

| State | Criterion | Treatment |
|---|---|---|
| *Unreached* | $n < k_{\min,\nu}$ | Skip. |
| *Partially Arrived* | $k_{\min,\nu} \le n < k_{\max,\nu}$ | Continue tree subdivision until a leaf is reached. |
| *Fully Arrived* | $k_{\max,\nu} \le n$ | Process the full block and merge upward. |

At each time step, the far-field block $b$ is first classified as a whole using the integer-delay interval of its root pair $\nu_b$. If $n < k_{\min,b}$, the block is skipped. If $k_{\max,b} \le n$, the block is *Fully Arrived*; its history contribution is subsequently evaluated over the complete index range, and its local dynamic cluster-pair tree records are released. For a block satisfying $k_{\min,b} \le n < k_{\max,b}$, an existing local tree $\mathcal{T}_b$ is reused and updated. If no tree exists, one is constructed by subdividing the root pair $\nu_b$ on demand.

The local tree is updated recursively from the root pair. For a *Partially Arrived* cluster pair requiring subdivision, the subdividable side containing more elements is selected. Its two spatial child clusters are paired separately with the other cluster, producing two child cluster pairs. Each newly generated pair is immediately classified at the current time step. Downward traversal stops when an *Unreached* or *Fully Arrived* pair is encountered. If a leaf cluster pair remains *Partially Arrived*, subdivision stops, and the condition $k_{ij} \le n$ is checked entrywise during history-convolution evaluation.

To reduce repeated traversal between consecutive time steps, each dynamic cluster pair stores its next possible update time step, $n_{\text{next},\nu}$. The rules for setting this value are given in Tab. 2. For clarity, a partially arrived pair whose two clusters are both spatial leaves is referred to as *Partially Arrived-Leaf*, whereas a partially arrived pair that can be further subdivided and retains child cluster pairs is referred to as *Partially Arrived-Refined*.

Tab. 2 Rules for setting $n_{\text{next},\nu}$

| State of $\nu$ | $n_{\text{next},\nu}$ |
|---|---|
| *Unreached* | $k_{\min,\nu}$ |
| *Partially Arrived-Leaf* | $k_{\max,\nu}$ |
| *Partially Arrived-Refined* | $\min(n_{\text{next},\nu_1}, n_{\text{next},\nu_2})$ |
| *Fully Arrived* | $+\infty$ |

Note: $\nu_1$ and $\nu_2$ denote the two immediate child cluster pairs of $\nu$, and $+\infty$ indicates that no further updates to the arrival state are required.

This variable specifies only the earliest time step at which the current cluster pair or one of its descendants may undergo a state change; it is not necessarily the time step at which the pair becomes *Fully Arrived*. A newly created pair is initialized as *Unreached*, with $n_{\text{next},\nu} = k_{\text{min},\nu}$. Updating is skipped if the pair is *Fully Arrived* or if $n < n_{\text{next},\nu}$. If $k_{\text{max},\nu} \le n$, the pair is marked as *Fully Arrived* and all descendant records are released. Otherwise, the pair is treated as partially arrived. If subdivision is possible, it is classified as *Partially Arrived-Refined*, and its child pairs are updated recursively; otherwise, it is classified as *Partially Arrived-Leaf*, and entrywise arrival tests are deferred to history-convolution evaluation. After updating, the next possible update time step is set according to Tab. 2.

During the return phase of the recursion, if both child pairs are *Fully Arrived*, they are merged into their parent pair, which is marked as *Fully Arrived*, and the descendant records are released. Because the child clusters provide tighter distance bounds, merging can establish that a parent pair is *Fully Arrived* before its conservative upper-bound criterion is satisfied. Once the root pair $\nu_b$ becomes *Fully Arrived*, the local dynamic tree for the block is released, and subsequent history contributions are evaluated directly over the entire block.

The evolution of the dynamic cluster-pair tree during initial wavefront propagation is illustrated in Fig. 2. The blue curves represent the initial wavefront at the corresponding time steps. The gray-filled blocks indicate previously fully arrived cluster pairs, while the beige-filled blocks represent newly fully arrived cluster pairs. The solid black lines indicate the fixed matrix partition, and the dashed lines denote temporary dynamic subdivisions. As the wavefront advances, partially arrived cluster pairs are locally subdivided, and their child pairs are merged once the parent pair becomes fully arrived.

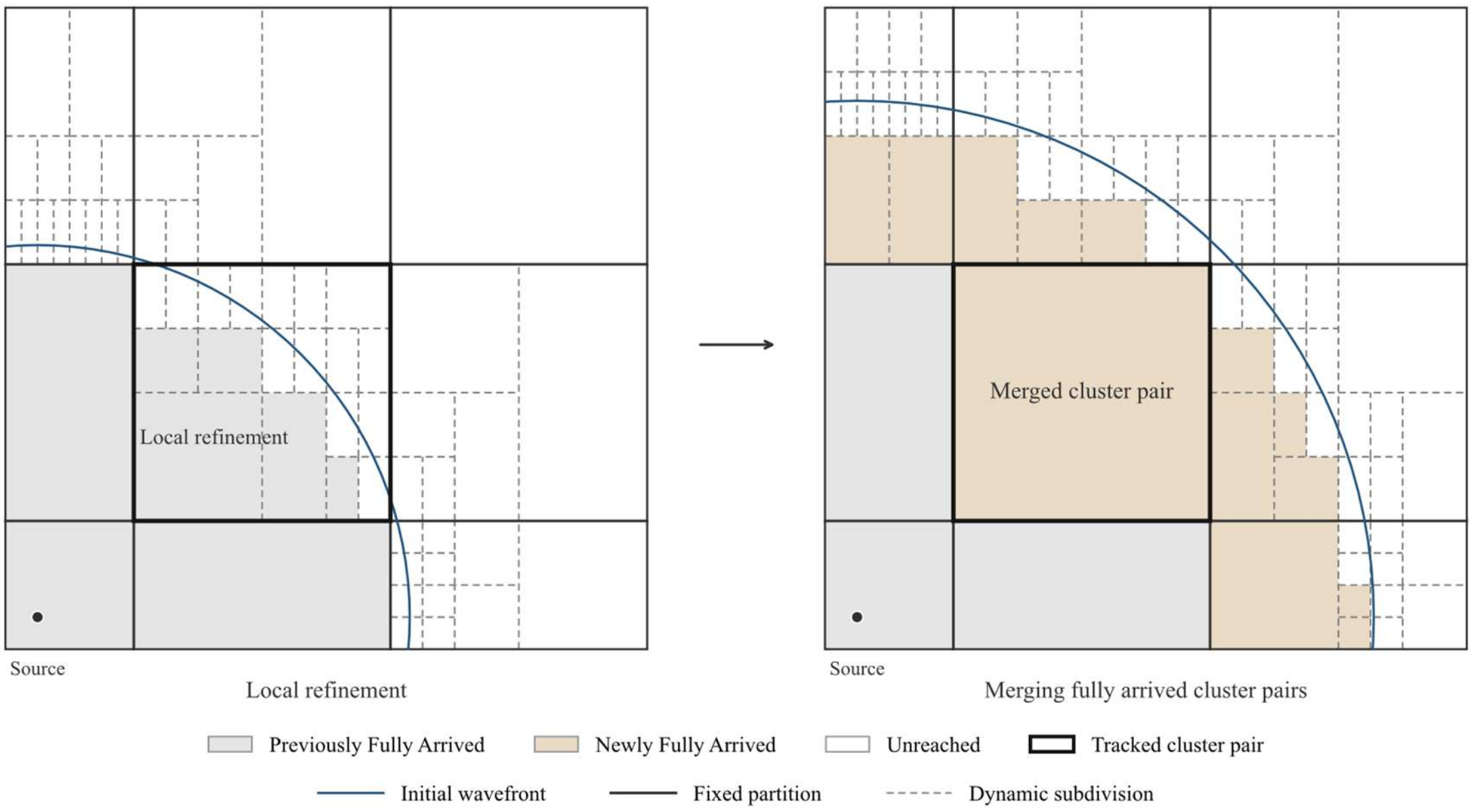

Fig. 2 Initial wavefront propagation and subdivision and merging of dynamic cluster pairs

Algorithm 2 details the update procedure.

**Algorithm 2 Updating the local dynamic cluster-pair tree of a far-field block**

**Input**: Far-field block $b$; local dynamic cluster-pair tree $\mathcal{T}_b$, which may be empty; current time step $n$; binary cluster tree; integer-delay interval $[k_{\min,b}, k_{\max,b}]$ of root pair $\nu_b$.

UpdatePair$(\nu_b, n)$

subroutine UpdatePair$(\nu, n)$
  if $\nu$ is *Fully Arrived* or $n < n_{\mathrm{next},\nu}$ then
    return
  end if
  if $k_{\mathrm{ma\ },\nu} \le n$ then
    Mark $\nu$ as *Fully Arrived* and release all descendants.
  else
    if $\nu$ can be further subdivided then
      Mark $\nu$ as *Partially Arrived-Refined*.
      if child pairs $\nu_1$ and $\nu_2$ do not exist then
        Create and initialize $\nu_1$ and $\nu_2$.
      end if
      UpdatePair$(\nu_1, n)$
      UpdatePair$(\nu_2, n)$
      if $\nu_1$ and $\nu_2$ are both *Fully Arrived* then
        Mark $\nu$ as *Fully Arrived* and release all descendants.
      end if
    else
      Mark $\nu$ as *Partially Arrived-Leaf*
    end if
  end if
  Update $n_{\mathrm{next},\nu}$ according to Tab. 2.

**Output**: Arrival state of $\nu_b$ and the updated tree $\mathcal{T}_b$.

After updating, history convolution is evaluated directly over the entire block for *Fully Arrived* blocks. For blocks that remain partially arrived, the local tree is traversed, *Unreached* pairs are skipped, and *Fully Arrived* regions and *Partially Arrived-Leaf* regions that cannot be subdivided further are treated separately, as described in Section 4.3. The local tree changes only the organization of the computational regions; it does not alter the spatial coefficients, propagation delays, or temporal interpolation rules.

## 4.3 History convolution using compressed matrices

This section presents the history-convolution procedure based on PACA-compressed equivalent matrices for *Fully Arrived* regions and *Partially Arrived-Leaf* regions within $\mathcal{H}$-matrix far-field blocks.

For a far-field block $b$, let $\Omega_b$ denote a computational region identified by its local dynamic cluster-pair tree, with row and column index sets $I_{\Omega_b} \subseteq I_b$ and $J_{\Omega_b} \subseteq J_b$, respectively. For any row index $i \in I_{\Omega_b}$ within block $b$, the column index set currently participating in the computation is defined as

$$J_{\Omega_b,i}^{n} = \begin{cases} J_{\Omega_b}, & \Omega_b \text{ is } Fully\ Arrived \\ \{ j \in J_{\Omega_b} : n \ge k_{ij} \}, & \Omega_b \text{ is } Partially\ Arrived - Leaf \end{cases} \tag{29}$$

For row index $i$, combining the three retarded boundary quantities $\bar{p}_{ij}^{n}$, $\bar{q}_{ij}^{n}$, and $\Delta p_{ij}^{n}$, defined in Section 3.3, with the PACA low-rank factors of the three equivalent matrices gives the far-field history contribution of the active region $\Omega_b$ at time step $n$:

$$
\begin{aligned}
f_{i,\Omega_b}^{n} \approx {} & \sum_{l=1}^{r_{S,b}} U_{S,b}(i,l) \sum_{j\in J_{\Omega_b,i}^{n}} V_{S,b}(j,l)\, \bar{q}_{ij}^{n} - \sum_{l=1}^{r_{D_0,b}} U_{D_0,b}(i,l) \sum_{j\in J_{\Omega_b,i}^{n}} V_{D_0,b}(j,l)\, \bar{p}_{ij}^{n} \\
& - \frac{1}{c\Delta t} \sum_{l=1}^{r_{D_1,b}} U_{D_1,b}(i,l) \sum_{j\in J_{\Omega_b,i}^{n}} V_{D_1,b}(j,l)\, \Delta p_{ij}^{n}
\end{aligned}
\tag{30}
$$

The delay parameters $k_{ij}$ and $\alpha_{ij}$ used to evaluate the retarded boundary quantities in Eq. (23) depend on both the collocation-point and source-element indices. Consequently, Eq. (30) cannot be reduced to a standard low-rank matrix–vector product in which all rows share the same input vector. The retarded boundary quantities must instead be evaluated entrywise for each pair $(i,j)$, row by row. Since the PACA low-rank factors are stored contiguously according to the binary cluster-tree ordering, the corresponding index ranges can be accessed directly when evaluating Eq. (30). Each computational region lies within its associated $\mathcal{H}$-matrix block, allowing direct access to the corresponding row and column ranges of the block's low-rank factors without repartitioning the matrix or reconstructing the factors.

For *Fully Arrived* regions, every collocation point–element pair satisfies $n \geq k_{ij}$. Entrywise arrival tests are therefore unnecessary, and all entries in the computational region $\Omega_b$ can be traversed sequentially. For *Partially Arrived-Leaf* regions, the condition $n \geq k_{ij}$ must be checked for each entry during evaluation. Contributions from arrived entries are evaluated using Eq. (30), whereas unreached entries are skipped. Eq. (30) is formulated specifically for PACA low-rank factors; blocks that cannot be compressed effectively are evaluated directly using the conventional approach.

Within each $\mathcal{H}$-matrix far-field block, the mutually disjoint computational regions identified by the local dynamic tree are traversed, and their history contributions are accumulated. *Fully Arrived* blocks are evaluated over their complete index ranges. The contributions from all far-field blocks are then summed to obtain the far-field history term $\boldsymbol{f}_{\mathrm{f}}^{n}$ at the current time step. This procedure skips unreached regions as a whole and restricts entrywise arrival tests to *Partially Arrived-Leaf* ranges containing the initial wavefront. It thereby reduces unnecessary distance calculations, wavefront localization, and accesses to historical boundary quantities, lowering the cost of far-field history convolution.

Once the initial wavefront has passed through the entire far field, all dynamic cluster pairs are merged into *Fully Arrived* regions. Subsequent time steps require neither wavefront classification nor subdivision; all far-field matrix blocks are traversed directly, and their history terms are evaluated using the *Fully Arrived* procedure. The same formulation also applies when the local dynamic cluster-pair tree is disabled. In that case, the entire far-field block is taken as the computational region, and entrywise arrival tests are performed for partially arrived blocks.

### 4.4 Computational complexity and storage requirements

The computational complexity is analyzed in terms of four components: near-field coefficient treatment, construction and storage of the far-field equivalent matrices, dynamic cluster-pair tree updates, and history convolution. For a far-field block $b$ with $m_b$ rows and $n_b$ columns, assume that the approximation rank of each time-independent equivalent matrix is bounded by $r$. The conventional PACA construction then requires approximately $O[r^2(m_b+n_b)]$ operations and $O[r(m_b+n_b)]$ storage. For an equivalent matrix with $N_s$ degrees of freedom, assuming a cluster-tree depth of $O(\log N_s)$ and a bounded block cluster-tree sparsity constant, the total costs of constructing and storing the compressed far-field equivalent matrices are approximately $O(r^2 N_s \log N_s)$ and $O(r N_s \log N_s)$, respectively. In addition, if the number of nonzero entries per row of the sparse current-step system matrix $\boldsymbol{A}_0$ is bounded, both its storage requirement and the cost of a single matrix-vector product are $O(N_s)$. Assuming that the number of iterations required by the iterative solver remains bounded with respect to $N_s$, the cost of solving the linear system at each time step is therefore also $O(N_s)$.

History convolution consists of far-field and near-field contributions. In the far field, the retarded boundary quantities still depend on both the row and column indices. With the rowwise evaluation procedure in Section 4.3, a fully arrived far-field block $b$ requires $O[r m_b n_b]$ operations per time step. Summation over all far-field blocks and time steps gives an upper bound of $O(rN_tN_s{}^2)$ for the far-field history-convolution cost. In the near field, quadrature-point delays are retained, and the corresponding coefficients are accessed in CSR format. If the number of near-field elements associated with each collocation point is bounded and the numbers of quadrature points and the temporal interpolation order are fixed, the near-field history-convolution cost over the full simulation is $O(N_tN_s)$, with a corresponding storage requirement of $O(N_s)$. Thus, the upper bound on the overall history-convolution cost is dominated by the far-field contribution and is $O(rN_tN_s{}^2)$.

These estimates provide upper bounds based on the complete far-field index range. During initial wavefront propagation, the local dynamic cluster-pair trees $\mathcal{T}_\mathrm{d}$ within the $\mathcal{H}$-matrix blocks can further reduce unnecessary accesses in history convolution by skipping unreached regions, at the expense of additional node storage and update operations. When $\beta$ is of order unity, the initial wavefront advances by $O(h)$ between consecutive time steps, as illustrated in Fig. 2. Local subdivision and state changes are therefore concentrated primarily near the wavefront. Unreached regions retain a coarse hierarchical representation, whereas *Fully Arrived* regions are merged into larger computational regions. Once an entire block becomes fully arrived, its history contribution is evaluated directly over the block, and its local dynamic records are released. The storage requirement of the local dynamic trees consequently depends on the number of nodes retained simultaneously.

Under uniform discretization, with no pronounced concentration in the propagation-delay distribution and a fixed leaf-cluster capacity, the total number of leaf cluster pairs in the local dynamic trees can be estimated as $O(N_s{}^2/\gamma_{\max})$. Here, the maximum time-lag index $\gamma_{\max}$, defined in Eq. (8), is determined jointly by the maximum propagation distance on the boundary and the time-step size. The storage requirement of the dynamic cluster-pair trees can therefore be estimated as $O(N_s{}^2/\gamma_{\max})$. Under this estimate, if $O(N_s{}^2/\gamma_{\max})$ cluster pairs within the narrow wavefront band actually require updating at each time step, and state classification, subdivision, and merging involve a fixed number of operations per pair, the cumulative dynamic-update cost during initial wavefront propagation is approximately $O\left(\gamma_{\max} \times \frac{N_s{}^2}{\gamma_{\max}}\right) = O(N_s{}^2)$. Compared with scanning all far-field indices at every time step, the local dynamic cluster-pair trees skip unreached regions. Under the stated assumption on the propagation-delay distribution, approximately $O(\gamma_{\max}N_s{}^2)$ pairwise checks can be avoided during initial wavefront propagation. The trees therefore primarily reduce unnecessary accesses during this initial stage. As the simulation duration increases, the effects of both the resulting savings and the update overhead on the average computation time per step gradually diminish, and full-block history convolution increasingly dominates the computation.

When the geometric dimensions and time-step size are fixed, $\gamma_{\max}$ is constant. When the time-step parameter $\beta$ is held fixed at a value of order unity, the time-step size decreases together with the element size, and $\gamma_{\max}$ increases with the number of elements [37, 40]. For a sphere or a general quasi-uniform surface mesh of fixed shape, the characteristic boundary element size satisfies $h \propto N_s{}^{-1/2}$. Thus, for fixed $\beta$, $\gamma_{\max} = O(N_s{}^{1/2})$. Compared with compact geometries such as spheres, slender geometries generally have larger values of $\gamma_{\max}$. In particular, if the element count increases primarily along one direction of the computational domain while remaining essentially constant in the other directions, $\gamma_{\max} = O(N_s)$ may occur. Consequently, for compact geometries such as spheres, the local dynamic cluster-pair trees $\mathcal{T}_\mathrm{d}$ associated with the $\mathcal{H}$-matrix blocks require $O(N_s{}^{3/2})$ storage. These trees store only a small amount of information per node, including cluster identifiers, child-node indices, integer-delay intervals, and states, so the storage constant is small.

The local dynamic cluster-pair tree is an optional acceleration strategy. Alternatively, arrival states can be classified directly using the delay interval of each block's root pair, with entrywise checks applied to partially arrived blocks. This avoids the storage and update costs of the within-block dynamic trees. Comparisons of the actual memory usage and computational costs are presented in Section 6.1.

In summary, the CEM-$\mathcal{H}$ method requires $O(rN_s \log N_s)$ matrix storage, and the combined cost of history convolution and system solution over $N_t$ time steps is bounded above by $O(N_t r N_s{}^2)$. In comparison, provided that the simulation duration covers the maximum propagation delay, a conventional method that stores the coefficient matrices densely for $\gamma_{\max}$ nonzero time lags and performs a dense matrix–vector product for each time lag requires $O(\gamma_{\max} N_s{}^2)$ storage and $O(N_t \gamma_{\max} N_s{}^2)$ operations for history convolution. The proposed local dynamic cluster-pair tree strategy further reduces unnecessary accesses during initial wavefront propagation at the expense of additional node storage and update operations. Its actual memory overhead and acceleration performance are evaluated in the numerical examples in Section 6.1.

This section primarily addresses conventional space-time discretizations, where the propagation distance within one time step is comparable to the characteristic element size ($\beta$ is of order unity), exploiting local access to PACA low-rank factors to handle the active regions as they evolve with the initial wavefront. The local dynamic cluster-pair tree changes only the organization of the computational regions, leaving the spatial coefficients and delayed sampling rules unchanged. It therefore introduces no approximation error beyond those associated with space-time discretization, the centroid-based equivalent representation, and PACA compression. For low-frequency problems in which the element size is much smaller than the wavelength, $\beta \gg 1$, or equivalently $c\Delta t \gg h$. When $L_{max}/(c\Delta t) = O(1)$, only a few integer-delay groups are required, and the wavefront may traverse multiple cluster pairs within a single time step. The computation therefore progressively shifts from local entries to complete cluster pairs. In this regime, $\mathcal{H}^2$-matrices exploit nested bases to share representation parameters among far-field blocks, making them well suited to matrix-vector operations involving complete cluster pairs.

# 5 Hierarchical far-field evaluation based on $\mathcal{H}^2$-matrices

This section develops an $\mathcal{H}^2$-matrix-based equivalent matrix compression method, termed CEM-$\mathcal{H}^2$, for low-frequency, long-wavelength problems in which the wavelength is large relative to the computational domain. Provided that the Nyquist sampling condition [15] and discretization accuracy requirements are satisfied, such problems permit relatively large time steps, making the distance traveled within one time step comparable to the maximum propagation distance across the domain. Consequently, only a small number of time-lag matrices are nonzero. In this regime, fewer integer-delay groups are required, and each group covers a larger spatial region, resulting in fewer but larger matrix blocks associated with far-field interactions. Compared with independent blockwise compression, $\mathcal{H}^2$-matrices employ nested cluster bases and interlevel transfer relations to propagate representations through the hierarchy. This reduces the repeated costs associated with independent representations and transformations for different cluster pairs, making the approach better suited to problems involving propagation over a large spatial region within a single time step. Unlike the $\mathcal{H}$-matrix method in Section 4, which focuses on local access to retarded data, the present method performs matrix-vector operations through a forward transformation on column clusters, intercluster interactions, and a backward transformation on row clusters.

Starting from the geometric partitioning in Section 3.1, geometrically admissible regions are further examined to determine whether the integer delay is consistent over the entire source-element integration domain. An $\mathcal{H}^2$-matrix representation is constructed for cluster pairs satisfying this condition. Based on the same geometric partitioning, the delay partitioning procedure in Section 5.2 identifies the far-field cluster pairs represented in the $\mathcal{H}^2$ format. Direct evaluation is retained for the geometric near field and for far-field regions that are ultimately excluded from the $\mathcal{H}^2$-matrix representation. Their propagation-delay treatment and coefficient construction are described in Section 5.2.

The method is applied both to history-term evaluation and to the action of the current-step coefficient matrix. At each time step, the right-hand side is assembled according to the boundary conditions, and the current unknown boundary quantities are obtained using an iterative solver. The spatial cluster bases, parent-child transfer matrices, and corresponding coupling matrices are constructed during preprocessing and reused throughout time marching.

5.1 Spatial hierarchical representation and transfer in the $\mathcal{H}^2$ format

In Section 4, a balanced binary cluster tree is employed for PACA low-rank compression. This section instead uses an adaptive octree-based cluster tree, retaining the definitions of the row and column index sets, cluster index sets, and contiguous ordering introduced in Section 3.1. During tree construction, the current spatial partitioning box is bisected at its midpoint along each of the three coordinate directions, producing eight subregions. Boundary elements are assigned to these subregions according to their centroid locations. Only nonempty child clusters are retained, and subdivision proceeds recursively until the stopping criterion is satisfied. The row and column clusters use the same spatial cluster tree. Geometric far-field candidate blocks are identified using the max-type geometric admissibility condition $\mathcal{M}_{+}\{a_1, a_2\}$ in Eq. (12), with $\eta_{\mathrm{adm}} = 2$. These candidate blocks are further partitioned according to the delay-consistency condition in Section 5.2 to obtain the cluster pairs actually used in the $\mathcal{H}^2$-matrix computations.

For any cluster $C$, its geometric bounding box $B_C$ is mapped affinely to the reference cube $[-1,1]^3$. Using $p_s$ Chebyshev nodes in each coordinate direction gives $n_s = {p_s}^3$ tensor-product interpolation nodes $\boldsymbol{\xi}_a^C$ and the corresponding Lagrange basis functions $L_a^C(\boldsymbol{x})$ [22, 32]. Unless otherwise stated, $i$ and $j$ denote collocation-point and source-element indices, respectively, while $a$ and $d$ denote interpolation-node indices for the row and column clusters, with $a, d = 1, \dots, n_s$. The two subscripts outside square brackets enclosing a matrix denote the row and column indices of the corresponding matrix entry. The construction of the two-sided Chebyshev-Lagrange interpolation, nested cluster bases, and parent-child transfer relations is detailed in Appendix C.

Based on this interpolation, the cluster bases for collocation points and source elements are constructed as $[\boldsymbol{P}_C]_{ia} = L_a^C(\boldsymbol{x}_i)$, $[\boldsymbol{Q}_C]_{ja} = \int_{\Gamma_j} L_a^C(\boldsymbol{y}) \mathrm{d}\Gamma_y$, $i, j \in K_C$. Here, $\boldsymbol{P}_C$ consists of basis-function values at the collocation points, whereas $\boldsymbol{Q}_C$ contains the integrals of the basis functions over the source elements. For the double-layer matrix, the outward normals of the source elements must also be incorporated. Since the normal to a planar triangular element is constant within the element, the normal-derivative column cluster basis is defined as $[\boldsymbol{Q}_C^D]_{ja} = \int_{\Gamma_j} \partial L_a^C(\boldsymbol{y}) / \partial n_y \mathrm{d}\Gamma_y = \int_{\Gamma_j} \boldsymbol{n}_j \cdot \nabla_y L_a^C(\boldsymbol{y}) \mathrm{d}\Gamma_y$, where the normal and the basis-function gradient are evaluated at the actual quadrature points of the source element.

For a far-field cluster pair $(C_c^b, C_s^b)$ represented in the $\mathcal{H}^2$ format, define the shared coupling matrix as

$$[\boldsymbol{G}_b]_{ad} = \frac{1}{4\pi \left\| \xi_a^{C_c^b} - \xi_d^{C_s^b} \right\|} \tag{31}$$

The single-layer equivalent matrix $\boldsymbol{S}_b$ and the double-layer equivalent matrix $\boldsymbol{D}_{0,b}$ are then approximated by

$$\begin{aligned} \boldsymbol{S}_b &\approx \boldsymbol{P}_{C_c^b} \boldsymbol{G}_b \boldsymbol{Q}_{C_s^b}^{\mathrm{T}} \\ \boldsymbol{D}_{0,b} &\approx \boldsymbol{P}_{C_c^b} \boldsymbol{G}_b \left( \boldsymbol{Q}_{C_s^b}^D \right)^{\mathrm{T}} \end{aligned} \tag{32}$$

The two representations share the row cluster basis $\boldsymbol{P}_{C_c^b}$ and coupling matrix $\boldsymbol{G}_b$, differing only in the source cluster bases $\boldsymbol{Q}_C$ and $\boldsymbol{Q}_C^D$. To further reduce the storage and multiplication costs of the coupling matrices, fully pivoted ACA is applied for recompression after $\boldsymbol{G}_b$ is constructed at the interpolation nodes [31]. In Section 5.3, the temporal interpolation relations are used to eliminate $\boldsymbol{D}_{1,b}$ for delay-consistent cluster pairs; hence, no separate $\mathcal{H}^2$-matrix representation is required for this matrix.

Interlevel transfer is achieved by exploiting the nested structure of the basis functions. Let $C_{\mathrm{p}}$ and $C_{\mathrm{ch}}$ denote a parent cluster and one of its child clusters, respectively. Define the transfer matrix as

$$\left[\boldsymbol{E}_{C_{\mathrm{ch}}\leftarrow C_{\mathrm{p}}}\right]_{da} = L_a^{C_{\mathrm{p}}}\left(\boldsymbol{\xi}_d^{C_{\mathrm{ch}}}\right) \tag{33}$$

Restricted to the child-cluster index set, the parent-cluster bases satisfy $\boldsymbol{P}_{C_{\mathrm{p}}}\big|_{K_{C_{\mathrm{ch}}}} = \boldsymbol{P}_{C_{\mathrm{ch}}}\boldsymbol{E}_{C_{\mathrm{ch}}\leftarrow C_{\mathrm{p}}}$, $\boldsymbol{Q}_{C_{\mathrm{p}}}\big|_{K_{C_{\mathrm{ch}}}} = \boldsymbol{Q}_{C_{\mathrm{ch}}}\boldsymbol{E}_{C_{\mathrm{ch}}\leftarrow C_{\mathrm{p}}}$. The normal-derivative column cluster basis $\boldsymbol{Q}_C^D$ satisfies the same nesting relation and can therefore share the parent-child transfer matrices with the scalar column cluster basis $\boldsymbol{Q}_C$. Accordingly, the column-cluster forward transformation uses the transpose of the transfer matrix, $\boldsymbol{E}_{C_{\mathrm{ch}}\leftarrow C_{\mathrm{p}}}^{\mathrm{T}}$, to aggregate source-side auxiliary vectors from child clusters to their parent. The row-cluster backward transformation uses $\boldsymbol{E}_{C_{\mathrm{ch}}\leftarrow C_{\mathrm{p}}}$ to transfer accumulated contributions from the parent to its children, and the row cluster bases recover the contributions at collocation points at the leaf level. Fig. 3 illustrates the two-sided spatial interpolation and nested cluster-basis transfers, including the projection of source-element contributions, intercluster coupling, and recovery of the resulting contributions at collocation points.

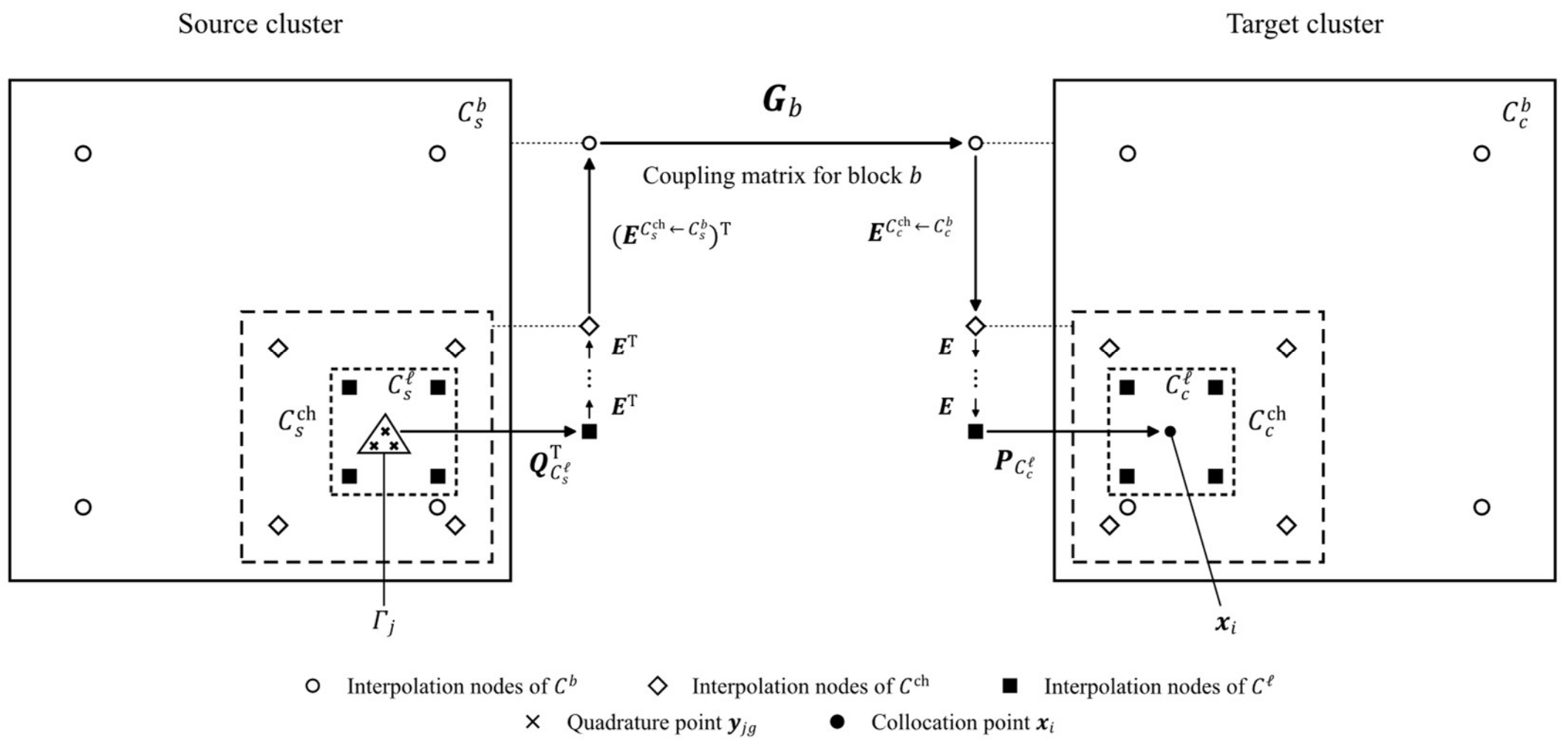

Fig. 3 Two-sided spatial interpolation and nested cluster-basis transfers in the $\mathcal{H}^2$-matrix representation

Thus, spatial interactions can be transferred between levels by explicitly storing only the leaf-level cluster bases and parent-child transfer matrices. These spatial structures are constructed during preprocessing and reused at every time step and during the iterative solution of the current-step system.

5.2 Delay-consistent partitioning and preprocessing of cluster pairs

In the regime considered here, $\beta \gg 1$ and $L_{max}/(c\Delta t) = O(1)$, so the propagation distance within one time step is comparable to the computational-domain size, and only a small number of integer-delay groups are required. Geometric far-field cluster pairs are therefore classified by integer delay during preprocessing, allowing accepted cluster pairs to share temporal sampling positions and employ an $\mathcal{H}^2$-matrix representation.

The delay-consistency procedure first examines geometric bounding boxes and then applies a centroid-based criterion to unresolved leaf cluster pairs, as detailed below.

**Delay consistency.** For any far-field block $b$, the integer-delay interval $[k_{\mathrm{min},b}, k_{\mathrm{max},b}]$ is determined from the geometric bounding boxes $B_{C_c^b}$ and $B_{C_s^b}$, as in Section 4.2. Since $\beta \gg 1$, we have $c\Delta t \gg h$. If $k_{\mathrm{min},b} = k_{\mathrm{max},b} = k_b$, the propagation times from all collocation points to all source-element integration positions within $b$ share the same integer delay $k_b$. The cluster pair is then said to satisfy the delay-consistency condition.

In this case, the entire cluster pair is added to the set $\mathcal{B}_{k_b}$ associated with integer delay $k_b$, and no further subdivision is performed. It should be emphasized that delay consistency requires only the integer delay to be identical. The fractional delay still varies with the collocation point and integration position: $\alpha_{ijg} = r_{ijg}/(c\Delta t) - k_b$. For these cluster pairs, propagation times are evaluated using quadrature-point delays. In Section 5.3, the fractional delay $\alpha_{ijg}$ is further combined with the spatial kernel to construct the corresponding $\mathcal{H}^2$ coupling matrices.

If $k_{\min,b} \neq k_{\max,b}$, the far-field cluster pair must be subdivided further. Among the clusters that can still be subdivided, priority is given to the one with the larger spatial extent. Each of its nonempty child clusters is paired with the other cluster, and the delay-consistency test is repeated.

When both the row and column clusters are leaves but $k_{\min,b} \neq k_{\max,b}$, further spatial subdivision is no longer available to improve the test. A centroid-based delay-consistency test is therefore applied. This test guarantees only that collocation point-element centroid pairs have the same integer delay; integer-delay consistency at all quadrature points within the elements is no longer strictly guaranteed. The centroid-based integer delay $k_{ij}$ is evaluated for every collocation point–element pair. If all values of $k_{ij}$ are identical, the entire cluster pair is added to $\mathcal{B}_{k_{ij}}$, and the common centroid-based integer delay is used as the time-interval index of the block.

If the centroid-based integer delays differ, the individual index pairs are assigned according to $k_{ij}$ to the corresponding geometric far-field direct-evaluation index set, $\mathcal{D}_k = \{(i,j) \in \mathcal{P}_{\mathrm{f,dir}}: k_{ij} = k\}$, where $\mathcal{P}_{\mathrm{f,dir}}$ denotes the direct-evaluation index region associated with far-field leaf clusters. The spatial coefficients $S_{ij}$, $D_{0,ij}$, and $D_{1,ij}$ are still obtained by spatial integration over the actual source elements. Only the temporal sampling positions use the centroid-based delay approximation.

To avoid unnecessary spatial subdivision, a bottom-up merging procedure is applied after the child cluster pairs have been processed. If all child blocks completely covering a parent block use an $\mathcal{H}^2$-matrix representation and share the same integer delay, they are replaced by the parent cluster pair, and the test continues upward. The merged parent pair retains the common integer delay of its children and does not alter the time-interval approximation described above.

After all geometric far-field candidate blocks have been traversed, a delay index table is constructed for each integer delay: $\mathcal{G}_k = (\mathcal{B}_k, \mathcal{D}_k)$, $k \in [0, \gamma_{\max}]$. For the particular class of problems considered in this section, $\gamma_{\max} \ll N_t$ generally holds. Here, $\mathcal{B}_k$ contains the cluster pairs evaluated in the $\mathcal{H}^2$ format using the common integer delay $k$, whereas $\mathcal{D}_k$ contains the direct-evaluation index pairs whose centroid-based integer delay is $k$. The two sets use different delay treatments, and their corresponding matrix indices do not overlap.

For leaf pairs accepted by the centroid-based criterion, assigning a common integer delay to the entire pair is an approximation, since individual quadrature points are not guaranteed to fall within the same delay interval. Quadrature points lying outside the selected interval are evaluated by extrapolating the linear temporal segment associated with that interval. This interval-assignment approximation is distinct from the centroid-based delay approximation introduced in Section 3.3.

The contributions associated with $\mathcal{G}_k$ involve the two time nodes $n-k$ and $n-k-1$. Coupling matrices for cluster pairs in $\mathcal{B}_k$ are constructed as described in Section 5.3. For direct-evaluation entries in $\mathcal{D}_k$, the interpolation weights are incorporated into the spatial coefficients in advance using Eq. (23). These coefficients are determined during preprocessing and reused at every time step. When $k = 0$, contributions associated with the current unknown boundary quantities enter the system solution, while known current-step contributions and history contributions are included in the right-hand side.

Fig. 4 illustrates the index partition after geometric block partitioning and delay treatment. The hatched regions represent leaf cluster pairs accepted after checking the centroid-based integer delays; their temporal weights are still constructed using quadrature-point distances. The gray-brown regions represent far-field index pairs retained for direct evaluation.

The entrywise delay calculations in Algorithm 3 are performed during preprocessing to construct the delay index tables and are not repeated during time marching.

**Algorithm 3 Cluster-pair traversal and construction of delay index tables**

**Input**: Octree-based cluster tree; set $Q$ of geometric far-field candidate cluster pairs.

while $Q \neq \emptyset$ do

  Remove the current cluster pair $(C_c, C_s)$ from $Q$, and compute $k_{\min,b}$ and $k_{\max,b}$.

  if $k_{\min,b} = k_{\max,b}$ then

    $\mathcal{B}_{k_{\min,b}} := \mathcal{B}_{k_{\min,b}} \cup \{b\}$

  else if $C_c$ or $C_s$ is not a leaf cluster then

    Subdivide the subdividable cluster corresponding to $\max\{\mathrm{diam}(B_{C_c}), \mathrm{diam}(B_{C_s})\}$, and add the resulting child cluster pairs to $Q$.

  else

    Compute $k_{ij}$ and $\alpha_{ij}$ for each index pair in the cluster pair using the centroid distance $R_{ij}$.

    if all $k_{ij}$ are identical then

      $\mathcal{B}_{k_{ij}} := \mathcal{B}_{k_{ij}} \cup \{b\}$

    else

      Assign each index pair to the corresponding $\mathcal{D}_{k_{ij}}$ according to $k_{ij}$.

      Use Eq. (23) to assemble the direct-evaluation coefficients for time lags $k$ and $k+1$.

    end if

  end if

end while

Proceed bottom-up, merging all child cluster pairs that completely cover a parent block and share the same integer delay into their parent cluster pair.

**Output**: Delay index tables $\mathcal{G}_k = (\mathcal{B}_k, \mathcal{D}_k)$, $k \in [0, \gamma_{\max}]$.

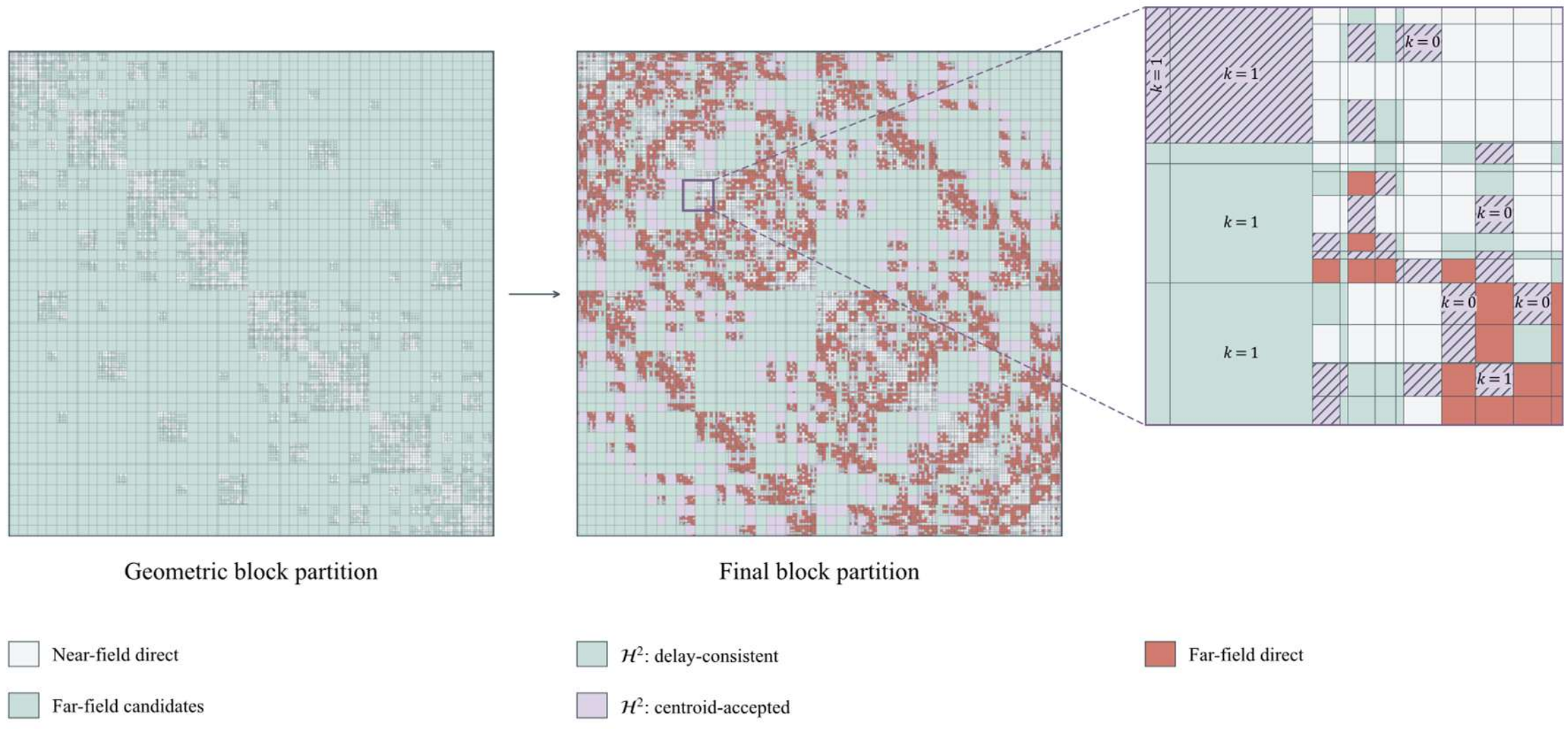

Fig. 4 Schematic of the matrix index partition after geometric block partitioning and delay treatment

## 5.3 Evaluation of retarded interactions and time marching using $\mathcal{H}^2$-matrices

For cluster pairs in $\mathcal{B}_k$, this section first combines the quadrature-point fractional delays with the spatial kernels to construct coupling matrices associated with adjacent time lags. The nested cluster bases introduced in Section 5.1 are then used to perform matrix-vector operations. Together with the directly evaluated contributions, these operations form the current-step equation and history terms.

For a cluster pair $(C_c^b, C_s^b) \in \mathcal{B}_k$ represented in the $\mathcal{H}^2$ format, the integer delay $k$ is determined as described in Section 5.2. The temporal weighting parameter is defined using the quadrature-point distance as $\alpha_{ijg} = r_{ijg}/(c\Delta t) - k$. The corresponding retarded boundary quantity is

$$\bar{\phi}_{ijg}^n = (1 - \alpha_{ijg})\phi_j^{n-k} + \alpha_{ijg}\phi_j^{n-k-1},\ \phi \in \{p, q\} \tag{34}$$

Unlike $\alpha_{ij}$ in Eq. (23), the fractional-delay parameter $\alpha_{ijg}$ is determined using quadrature-point distances and therefore depends additionally on the quadrature-point index $g$. However, for a cluster pair assigned to a common integer delay $k$, all contributions are evaluated using the same two time nodes, $n-k$ and $n-k-1$. The corresponding boundary values can therefore be taken outside the spatial integrals, while the spatially varying fractional-delay weights are incorporated into the equivalent matrix coefficients. This enables the history contributions to be evaluated through blockwise matrix-vector operations using nested cluster bases, rather than through rowwise access to individual retarded interactions.

To this end, define the single-layer matrix $\boldsymbol{S}_{r,b}$ for cluster pair $b$, weighted by the quadrature-point distance $r_{ijg}$, through the entries

$$S_{r,ij} = \int_{\Gamma_j} r_{ij} \frac{1}{4\pi r_{ij}} d\Gamma_y = \int_{\Gamma_j} \frac{1}{4\pi} d\Gamma_y \tag{35}$$

Using the fractional-delay weights defined in Eq. (34), the single-layer contribution can be expressed in terms of $\boldsymbol{S}_b$ and $\boldsymbol{S}_{r,b}$. For the double-layer contribution, the distance-weighted term introduced by temporal interpolation coincides with $\boldsymbol{D}_{1,b}$ and cancels the pressure time-derivative term when the same quadrature-point distances and common integer-delay interval are used. The detailed derivation is provided in Appendix C. Consequently, the double-layer contribution simplifies to

$$\boldsymbol{f}_{b,D}^n \approx \boldsymbol{D}_{0,b}\left[(k+1)\boldsymbol{p}_{C_s^b}^{n-k} - k\boldsymbol{p}_{C_s^b}^{n-k-1}\right] \tag{36}$$

For cluster pairs in $\mathcal{B}_k$, the resulting far-field interactions require only the $\mathcal{H}^2$-matrix representations of $\boldsymbol{S}_b$ and $\boldsymbol{D}_{0,b}$, which share the coupling matrix $\boldsymbol{G}_b$ as described in Section 5.1. Since $\boldsymbol{S}_{r,b}$ depends only on source-element areas, no additional full matrix is required. For the remaining far-field index pairs in $\mathcal{D}_k$, the three equivalent coefficients are retained using the centroid-based delay treatment of Section 3.3. The cluster bases, transfer matrices, and coupling matrices are constructed during preprocessing and reused throughout time marching.

In the implementation, a forward transformation is first performed on the column cluster tree for the boundary quantities at the required time nodes. For a leaf-level column cluster $C_s$, the auxiliary vectors at time node $m$ are

$$\begin{aligned} \hat{\boldsymbol{q}}_{C_s}^m &= \boldsymbol{Q}_{C_s}^{\mathrm{T}} \boldsymbol{q}_{C_s}^m \\ \hat{\boldsymbol{p}}_{C_s}^m &= \left(\boldsymbol{Q}_{C_s}^D\right)^{\mathrm{T}} \boldsymbol{p}_{C_s}^m \end{aligned} \tag{37}$$

where $\hat{\boldsymbol{q}}_{C_s}^m$ and $\hat{\boldsymbol{p}}_{C_s}^m$ are the auxiliary vectors obtained by applying the column-cluster forward transformation to the boundary quantities $\boldsymbol{q}_{C_s}^m$ and $\boldsymbol{p}_{C_s}^m$, respectively. They are used in subsequent intercluster coupling operations. The auxiliary vectors for the pressure normal derivative and pressure are formed using the scalar and normal-derivative column cluster bases, respectively, and then aggregated from child clusters to parent clusters level by level using the shared transfer matrices. The transformed vector for a given column cluster at a given time node can be shared by all associated cluster pairs.

For each cluster pair $(C_c^b, C_s^b) \in \mathcal{B}_k$, combining the single-layer contribution derived in Appendix C with the simplified double-layer contribution in Eq. (36), the intercluster contributions are accumulated into the row-cluster auxiliary vector as follows:

$$\begin{aligned} \boldsymbol{y}_{C_c^b}^n += \boldsymbol{G}_b \left[(k+1)\widehat{\boldsymbol{q}}_{C_s^b}^{n-k} - k\widehat{\boldsymbol{q}}_{C_s^b}^{n-k-1}\right] + \frac{1}{c\Delta t}\boldsymbol{G}_{S_r,b}\left[\widehat{\boldsymbol{q}}_{C_s^b}^{n-k-1} - \widehat{\boldsymbol{q}}_{C_s^b}^{n-k}\right] \\ -\boldsymbol{G}_b\left[(k+1)\widehat{\boldsymbol{p}}_{C_s^b}^{n-k} - k\widehat{\boldsymbol{p}}_{C_s^b}^{n-k-1}\right] \end{aligned} \tag{38}$$

The coupling matrix $\boldsymbol{G}_{S_r,b}$ need not be constructed separately. For a given cluster pair, the above calculation requires only the auxiliary vectors at time steps $n-k$ and $n-k-1$. When $k \geq 1$, Eq. (38) involves only known historical boundary quantities, and its contribution is included in the history term. When $k=0$, the coefficients associated with current unknown boundary quantities contribute to $\boldsymbol{A}_0$, while the remaining terms enter the right-hand side. Boundary quantities at negative time steps are set to zero in accordance with the zero initial conditions.

After the coupling-matrix contributions associated with all integer delays $k \in [0, \gamma_{\max}]$ in the history term have been evaluated and accumulated, a single backward transformation is performed on the row cluster tree:

$$\boldsymbol{y}_{C_{\mathrm{ch}}}^n += \boldsymbol{E}_{C_{\mathrm{ch}} \leftarrow C_{\mathrm{p}}} \boldsymbol{y}_{C_{\mathrm{p}}}^n \tag{39}$$

Upon reaching a leaf-level row cluster $C_c$, the history contributions at the collocation points are recovered through the row cluster basis: $\boldsymbol{f}_{C_c,\mathcal{H}^2}^n = \boldsymbol{P}_{C_c}\boldsymbol{y}_{C_c}^n$. Combining the results from all leaf clusters gives $\boldsymbol{f}_{\mathcal{H}^2}^n$. For far-field direct-evaluation indices in $\mathcal{D}_k$, the coefficients stored during preprocessing are read and multiplied by the corresponding historical boundary quantities to obtain $\boldsymbol{f}_{\mathrm{f,dir}}^n$. The far-field history term at time step $n$ is then $\boldsymbol{f}_{\mathrm{f}}^n = \boldsymbol{f}_{\mathcal{H}^2}^n + \boldsymbol{f}_{\mathrm{f,dir}}^n$.

The detailed procedure for evaluating the known far-field history contributions is provided in Appendix C. The current-step unknown contributions are handled separately during the iterative solution using the same hierarchical matrix operations.

### 5.4 Computational complexity and storage requirements

This section analyzes the storage requirements and computational complexity per time step for low-frequency, long-wavelength problems. For the following estimates, it is assumed that the delay-partitioned hierarchical block structure remains sparse, with both the number of far-field cluster pairs represented in the $\mathcal{H}^2$ format and the total number of far-field index pairs retained for direct evaluation scaling as $O(N_s)$. This assumption is particularly relevant to the low-frequency, long-wavelength regime considered here, where the propagation distance $c\Delta t$ can be comparable to the computational-domain size. In this regime, only a small number of integer-delay groups are required, allowing many far-field cluster pairs to adopt a common integer delay without extensive subdivision and reducing the regions requiring direct evaluation. The number of geometric near-field elements associated with each collocation point is also assumed to remain bounded.

For storage, the nested structure of the cluster bases requires explicit storage only of the leaf-level cluster bases and parent-child transfer matrices, with costs of $O(n_s N_s)$ and $O({n_s}^2 N_s)$, respectively. For far-field cluster pairs represented in the $\mathcal{H}^2$ format, the single- and double-layer equivalent matrices share coupling matrices. Assuming that the approximation rank after ACA recompression does not exceed $r_s$, the storage requirement for all low-rank coupling matrices is approximately $O(r_s n_s N_s)$. The geometric near field and far-field regions excluded from the $\mathcal{H}^2$ representation are evaluated directly. With a small number of nonzero time-lag matrices and $O(N_s)$ direct-evaluation pairs, this part requires $O(N_s)$ storage. The overall storage requirement is therefore $O[({n_s}^2 + r_s n_s + 1)N_s]$. For fixed $p_s$ and $r_s$, this becomes $O(N_s)$, with a storage constant that increases with $n_s = {p_s}^3$.

For history-term evaluation, the column-cluster auxiliary vectors at each time node are obtained by a single forward transformation once the boundary quantities become available. They are retained within the required history window and reused directly at subsequent time steps. A forward transformation comprises the construction of leaf-level auxiliary vectors and their upward transfer between child and parent clusters, at a cost of $O({n_s}^2 N_s)$. The cluster pairs then use the stored auxiliary vectors to evaluate low-rank intercluster interactions, with a total cost of approximately $O(r_s n_s N_s)$. A single backward transformation is performed on the row cluster tree, also at a cost of $O({n_s}^2 N_s)$.

For far-field index pairs excluded from the $\mathcal{H}^2$ representation, direct evaluation is retained during time marching. Although these interactions are evaluated individually, their total number is assumed to scale as $O(N_s)$, as stated above. Since the corresponding coefficients and delay indices are precomputed, the direct-evaluation cost per time step is $O(\sum_k |\mathcal{D}_k|) = O(N_s)$. Together with the bounded geometric near-field interactions, the total direct-evaluation cost remains $O(N_s)$ per time step. Thus, the history-term evaluation cost per time step is $O[({n_s}^2 + r_s n_s + 1)N_s]$, which likewise reduces to $O(N_s)$ when the relevant parameters are fixed. The linear system is solved using the generalized minimal residual (GMRES) method. With a fixed interpolation order and a bounded number of iterations at each time step, the solution cost per time step is approximately $O(N_s)$.

In summary, under the stated assumptions, the proposed CEM-$\mathcal{H}^2$ method requires $O(N_s)$ storage. With a bounded number of iterations, the combined computational cost of history convolution and iterative system solution over $N_t$ time steps is therefore $O(N_t N_s)$.

# 6 Numerical examples

This section presents numerical examples to assess the computational accuracy, compression efficiency, storage requirements, and computational performance of the proposed CEM-$\mathcal{H}$ and CEM-$\mathcal{H}^2$ methods. Accuracy is evaluated against semi-analytical solutions or conventional TDBEM results. Computational efficiency is assessed in terms of matrix storage, computation time, and the ability to solve large-scale problems. Comparisons with an existing fast method, IFMM [6], are also included where results are available. The methods use the same discretized models and time-integration parameters, and identical iterative convergence criteria are adopted unless otherwise specified. The parallel implementations and common computational parameters are described first, followed by the numerical results.

All methods employ OpenMP-based shared-memory parallelization to improve computational efficiency. For conventional TDBEM, matrix assembly and history-term evaluation are parallelized by assigning matrix rows to different tasks. For CEM-$\mathcal{H}$, far-field matrix compression is parallelized over blocks, and dynamic scheduling is used for history-term evaluation. For larger far-field blocks, the computational tasks can be further partitioned along the target-row direction to improve load balancing. Thread-local contributions are combined after evaluation. For CEM-$\mathcal{H}^2$, the upward and downward passes are parallelized over clusters at the same level, with synchronization between levels. Intercluster coupling operations are grouped by target cluster for parallel execution, while direct evaluation is parallelized over target rows.

The linear system at each time step is solved using GMRES with a relative residual tolerance of $10^{-6}$. The leaf-cluster capacity is set to $n_{\text{leaf}} = 64$ for both CEM-$\mathcal{H}$ and CEM-$\mathcal{H}^2$. All numerical examples are performed on a computer equipped with two Intel Xeon CPUs (Gold 5215, 20 physical cores), and 128 GB of memory. The number of OpenMP threads is specified for each example.

## 6.1 Sphere

The transient scattering example for a rigid sphere considered by Takahashi [6] is used to assess the computational accuracy, memory requirements, and efficiency of CEM-$\mathcal{H}$ and CEM-$\mathcal{H}^2$. The sphere has a diameter of $D = 1$ and is centered at $(0.5,0,0)$. The sound speed in the unbounded exterior medium is $c = 1$, and the sphere surface satisfies the rigid boundary condition $q = \partial p/\partial n = 0$. An incident plane wave propagates in the positive $x$-direction. Both general transient and low-frequency, long-wavelength cases are considered, with a semi-analytical solution serving as the accuracy reference. All methods in this subsection use 16 OpenMP threads.

a. General transient analysis using CEM-$\mathcal{H}$

A compactly supported cosine pulse is used as the incident excitation. The sphere geometry, incident pulse, and boundary mesh are shown in Fig. 5. The incident acoustic pressure is given by

$$p^{\mathrm{in}}(\boldsymbol{x},t)=\begin{cases}A\{1-\cos(\frac{2\pi c\tau}{\lambda})\}, & 0\leq c\tau\leq\lambda\\ 0, & \text{other}\end{cases} \tag{40}$$

where $\tau=t-x/c$, $A$ is the amplitude coefficient, and $\lambda=D/2$ is the spatial support length of the pulse. The peak pulse amplitude is $2A=1$. At the initial time, the incident wavefront is located at $x=0$ and subsequently propagates toward the sphere. The total acoustic pressure $p$ on the sphere surface is computed over the analysis interval $T=4$. For exterior scattering problems, the total pressure $p=p^{\mathrm{in}}+p^{\mathrm{sc}}$ is taken as the unknown boundary quantity, where $p^{\mathrm{sc}}$ is the scattered pressure. The incident pressure $p^{\mathrm{in}}$ is included as a known excitation on the right-hand side of the discrete equation.

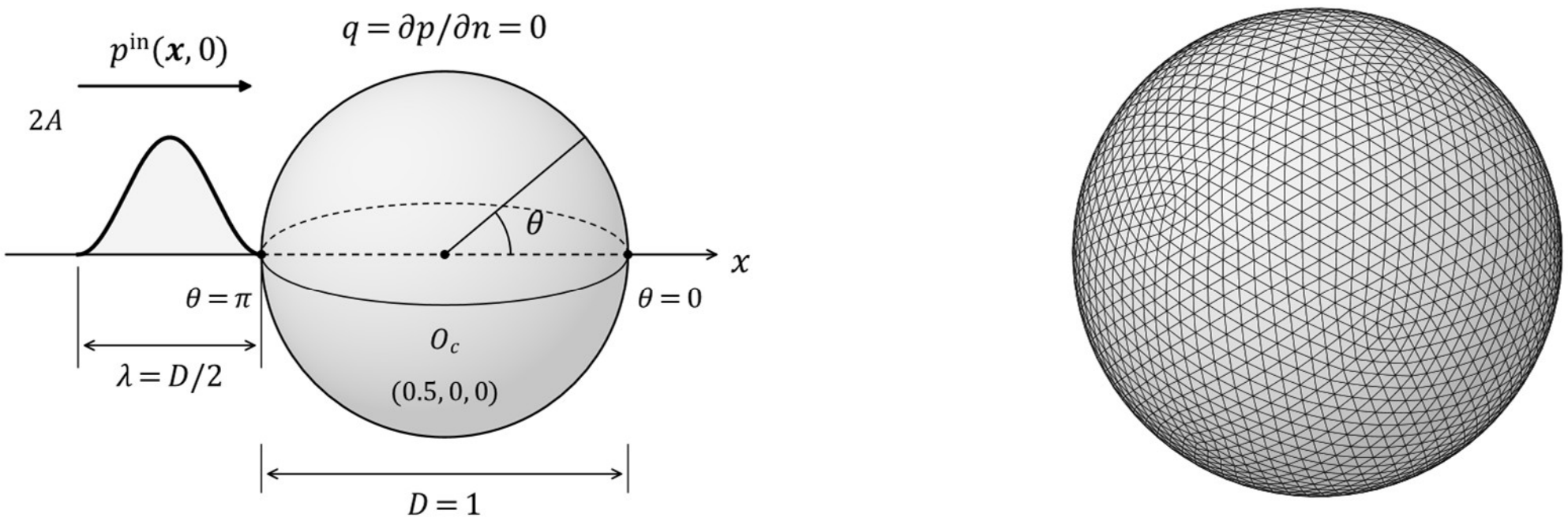

(a) Incident pulse and sphere geometry (b) Boundary mesh ($N_s=5120$)

Fig. 5 Schematic of the transient acoustic scattering example for a sphere

The sphere surface is discretized using constant triangular boundary elements. The time-step size is reduced together with mesh refinement so that the time-step parameter $\beta$ remains of order unity. Numerical results obtained on meshes of different resolutions are first compared with the semi-analytical solution to assess accuracy. The discretization size is then increased further to investigate memory requirements and computational efficiency. The discretization parameters are listed in Tab. 3.

To assess the overall accuracy, the $x$-axis is taken as the polar axis, and $33$ boundary collocation points are selected whose polar angles are closest to $\theta_i=i\pi/32,\ i=0,1,\ldots,32$. The discrete relative $\ell_2$ error with respect to the semi-analytical solution is evaluated at these locations over the entire analysis interval:

$$E_{rel}=\sqrt{\frac{\sum_{0\leq i<3}\ \sum_{0\leq\alpha<N_{\mathrm{t}}}\left[p_h(x(\theta_i),t_\alpha)-p_{ref}(x(\theta_i),t_\alpha)\right]^2}{\sum_{0\leq i<33}\sum_{0\leq\alpha<N_{\mathrm{t}}}\left[p_{ref}(x(\theta_i),t_\alpha)\right]^2}} \tag{41}$$

where $p_h$ denotes the numerical solution and $p_{\mathrm{ref}}$ the semi-analytical solution.

Tab. 3 Discretization parameters for the numerical examples

| $N_s$ | $\Delta t$ | $\beta$ | $N_t$ | $\gamma_{\max}$ |
|---|---|---|---|---|
| 320 | 0.16 | 1.07 | 25 | 7 |
| 1,280 | 0.08 | 1.07 | 50 | 13 |
| 5,120 | 0.04 | 1.06 | 100 | 25 |
| 20,480 | 0.02 | 1.06 | 200 | 50 |
| 81,920 | 0.01 | 1.06 | 400 | 100 |
| 327,680 | 0.005 | 1.06 | 800 | 200 |

Fig. 6 compares the errors of conventional TDBEM and CEM-$\mathcal{H}$ relative to the semi-analytical solution. For CEM-$\mathcal{H}$, PACA relative tolerances of $\varepsilon_{\text{rel}} = 10^{-4}$ and $\varepsilon_{\text{rel}} = 10^{-6}$ are considered, with all other computational parameters held fixed, to examine the influence of the tolerance on overall accuracy. The relative error $E_{\text{rel}}$ decreases consistently for all methods as the number of elements $N_s$ increases. The error curves for the two values of $\varepsilon_{\text{rel}}$ nearly coincide and remain close to those of the conventional method. Thus, $\varepsilon_{\text{rel}} = 10^{-4}$ is sufficient to retain accuracy comparable to conventional TDBEM. Unless otherwise stated, this tolerance is used in all subsequent examples.

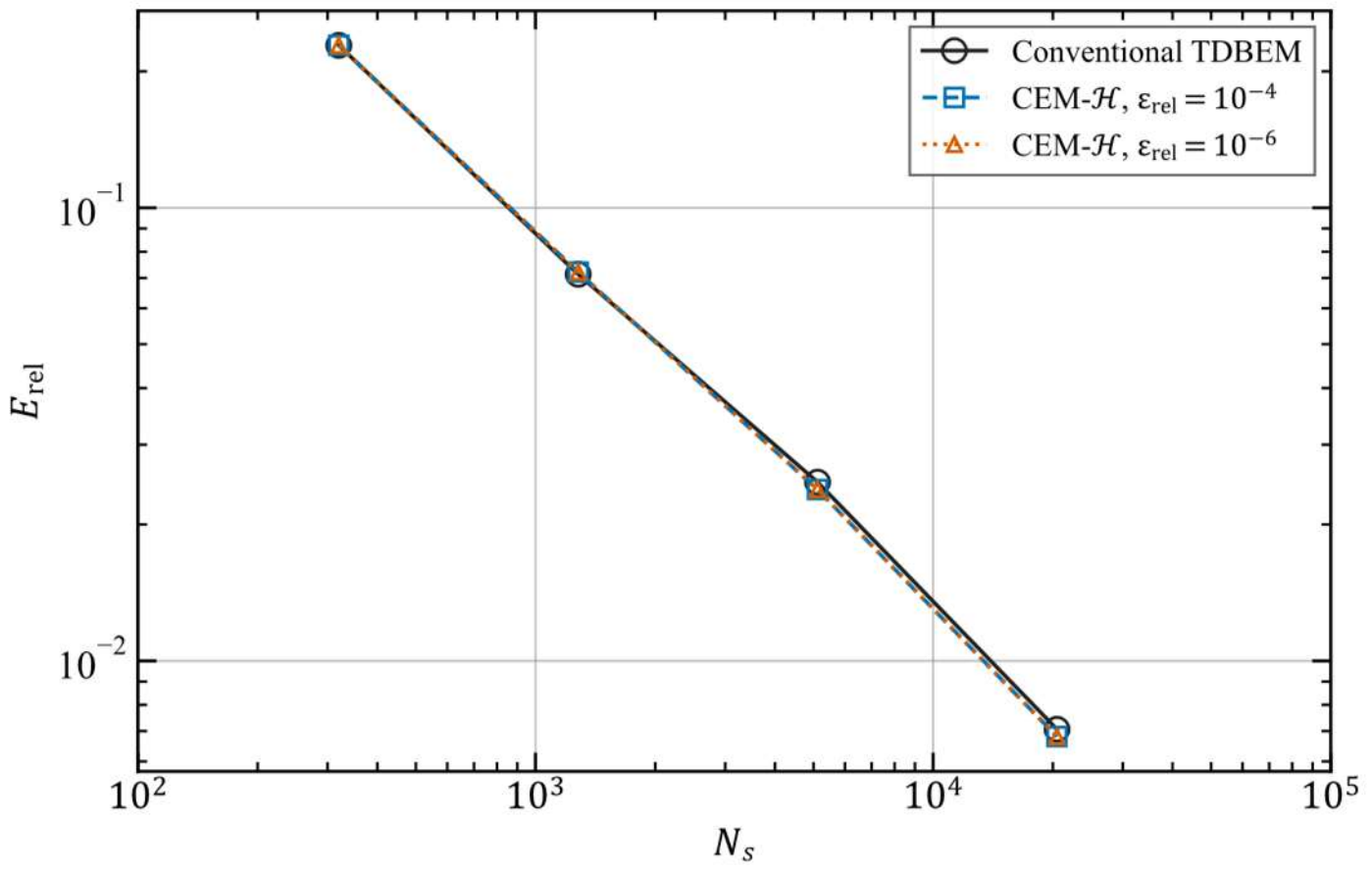

Fig. 6 Relative errors of conventional TDBEM and CEM-$\mathcal{H}$ as a function of the number of elements

Fig. 7 presents the total-pressure time histories at the front $\theta = \pi$ and rear $\theta = 0$ of the sphere, together with the semi-analytical solution, conventional TDBEM results, and IFMM results [6]. As the number of elements increases, the CEM-$\mathcal{H}$ results approach the semi-analytical solution and remain in good agreement with the conventional results obtained using the same discretization. For the same spatial discretization, the IFMM results [6] exhibit some approximation error when relatively few interpolation nodes are used. Increasing the number of interpolation nodes brings the results closer to those of conventional TDBEM and the semi-analytical solution. It should be noted that the time-step sizes used in the present study differ from those in [6], and the influence of this difference is not accounted for in the comparison.

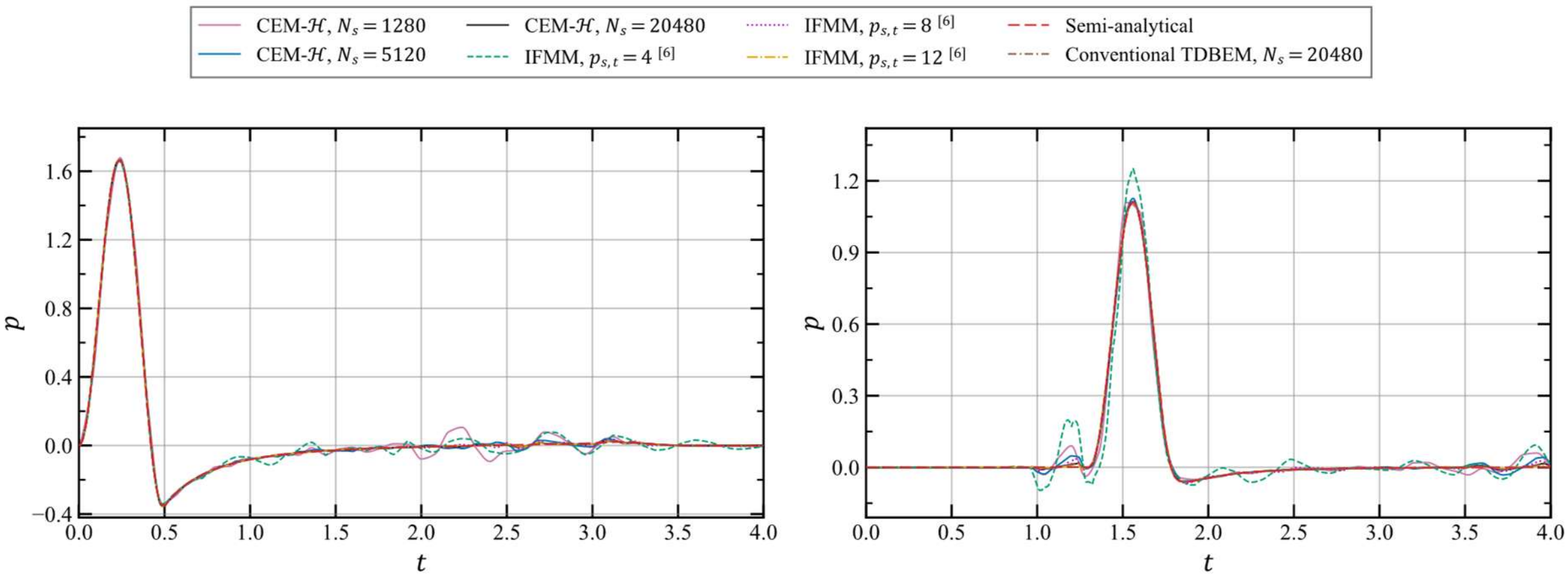

(a) $\theta = \pi$ (b) $\theta = 0$

Fig. 7 Comparison of total-pressure time histories at the front and rear of the sphere

To investigate the influence of the time-step parameter on numerical stability, Fig. 8 compares the total-pressure histories for different values of $\beta$ at $N_s = 1{,}280$. As shown in Fig. 8(a), CEM-$\mathcal{H}$ exhibits rapidly growing late-time oscillations at $\beta = 0.1$ and weaker growing oscillations at $\beta = 0.5$, whereas no noticeable growth is observed for $\beta = 1$, $2$ and $4$ over the simulated interval of $40$ s. This behavior is consistent with the sensitivity to small time-step parameters reported for certain conventional MOT formulations, as discussed in Section 2.2. Fig. 8(b) further compares CEM-$\mathcal{H}$ and conventional TDBEM at $\beta = 0.1$ using $27$ and $192$ quadrature points. Increasing the quadrature resolution suppresses the observed instability in conventional TDBEM but does not eliminate it in CEM-$\mathcal{H}$. These results indicate that CEM-$\mathcal{H}$ is conditionally stable for the investigated discretization, with its stability potentially affected by the additional centroid-based delay approximation. This trend is consistent with the local delay-deviation analysis in Section 3.4: for a fixed spatial mesh, reducing $\beta$ decreases the time step and therefore increases the relative delay deviation $\eta_{ij}^{\tau} = \delta\tau_{ij}/\Delta t$. Since $\delta\tau_{ij}$ is determined by the spatial configuration and sound speed $c$ rather than the quadrature rule, increasing the quadrature resolution does not reduce the centroid-based delay deviation itself.

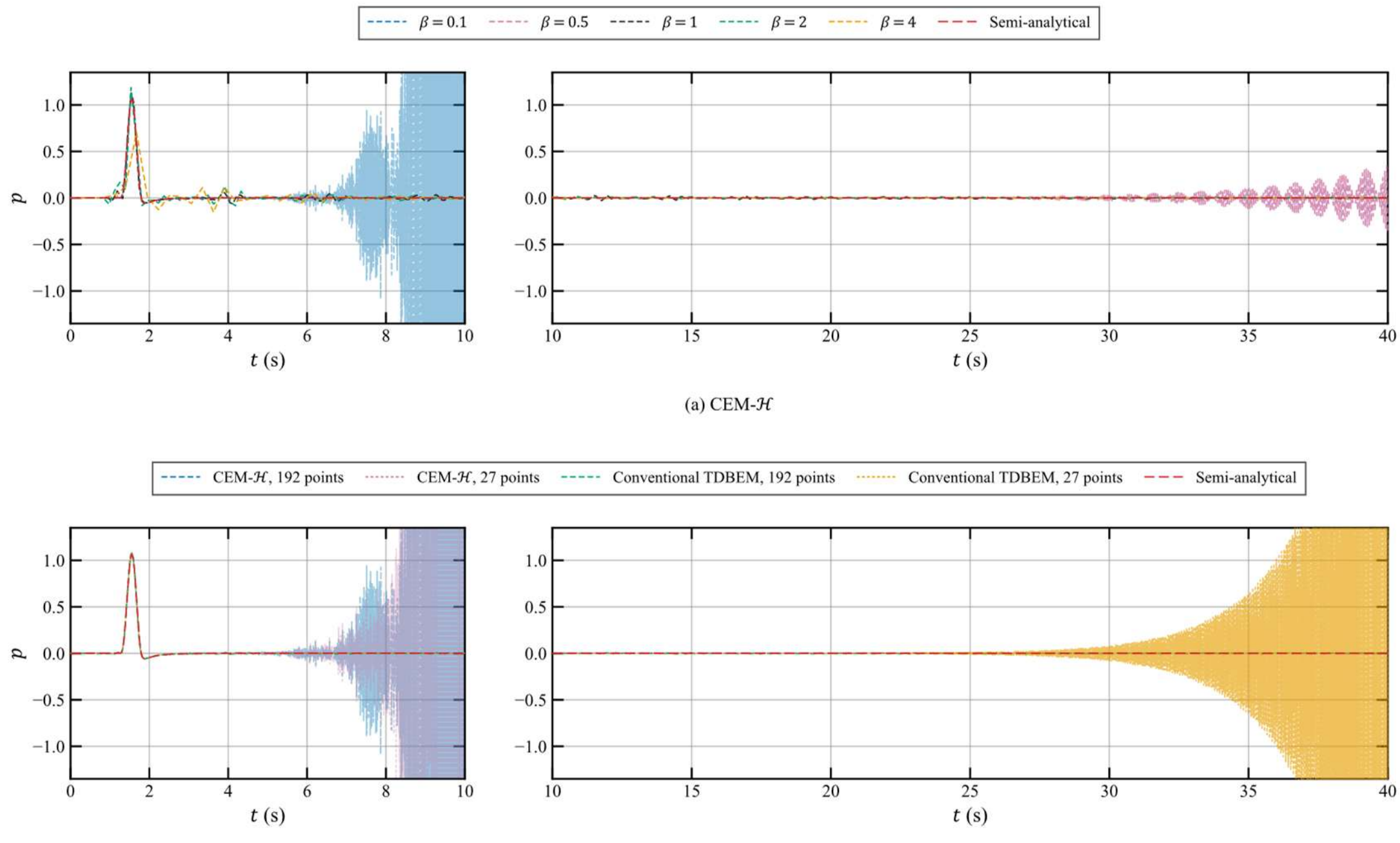

(a) CEM-$\mathcal{H}$

(b) $\beta = 0.1$

Fig. 8 Influence of the time-step parameter $\beta$ on the total-pressure histories at the rear of the sphere ($N_s = 1{,}280$)

Tab. 4 lists the storage requirements and total computation times of the two methods for different problem sizes. At $N_s = 20{,}480$, the estimated memory requirement of conventional TDBEM is approximately $349.0$ GB, whereas CEM-$\mathcal{H}$ requires only approximately $1.6$ GB. The speedup relative to the estimated computation time of the conventional method is approximately $45.4\times$, demonstrating substantial savings in both storage and computation. When the number of elements increases to $N_s = 327{,}680$, the computation can still be completed using approximately $28.6$ GB of memory.

Tab. 4 Memory requirements and total computation times of conventional TDBEM and CEM-$\mathcal{H}$

| $N_s$ | $N_t$ | Conventional TDBEM | | CEM-$\mathcal{H}$ | |
|---|---|---|---|---|---|
| | | Time (s) | Memory (GB) | Time (s) | Memory (GB) |
| 320 | 25 | 0.2 | 0.1 | 0.1 | 0.1 |
| 1,280 | 50 | 5.4 | 0.4 | 1.2 | 0.1 |
| 5,120 | 100 | 219.7 | 11.7 | 12.2 | 0.4 |
| 20,480 | 200 | 12,447.8* | 349.0* | 274.5 | 1.6 |
| 81,920 | 400 | — | — | 9,340.9 | 7.1 |
| 327,680 | 800 | — | — | 350,338.7 | 28.6 |

∗ Estimated values. For conventional TDBEM at $N_s = 20{,}480$, the reported storage is the theoretical requirement for the coefficient matrices, and the total computation time is estimated by extrapolation. In the actual computation, only a portion of the matrices is retained in memory, while the remaining portions must be recomputed.

To examine how storage requirements and computation time vary with problem size, Fig. 9 compares memory usage and average computation time per time step for different element counts. The latter is obtained by dividing the total computation time by the number of time steps. As the number of elements increases, both the memory usage and average time per step of conventional TDBEM exhibit approximately $O(\gamma_{\max} {N_s}^2)$ growth. In contrast, the storage requirement of CEM-$\mathcal{H}$ follows a trend close to $O(N_s \log N_s)$. Although its average time per step grows as $O({N_s}^2)$ for large $N_s$, this increase is slower than that of the conventional method. Consequently, the differences in storage and computation time between the two methods become larger as the problem size increases.

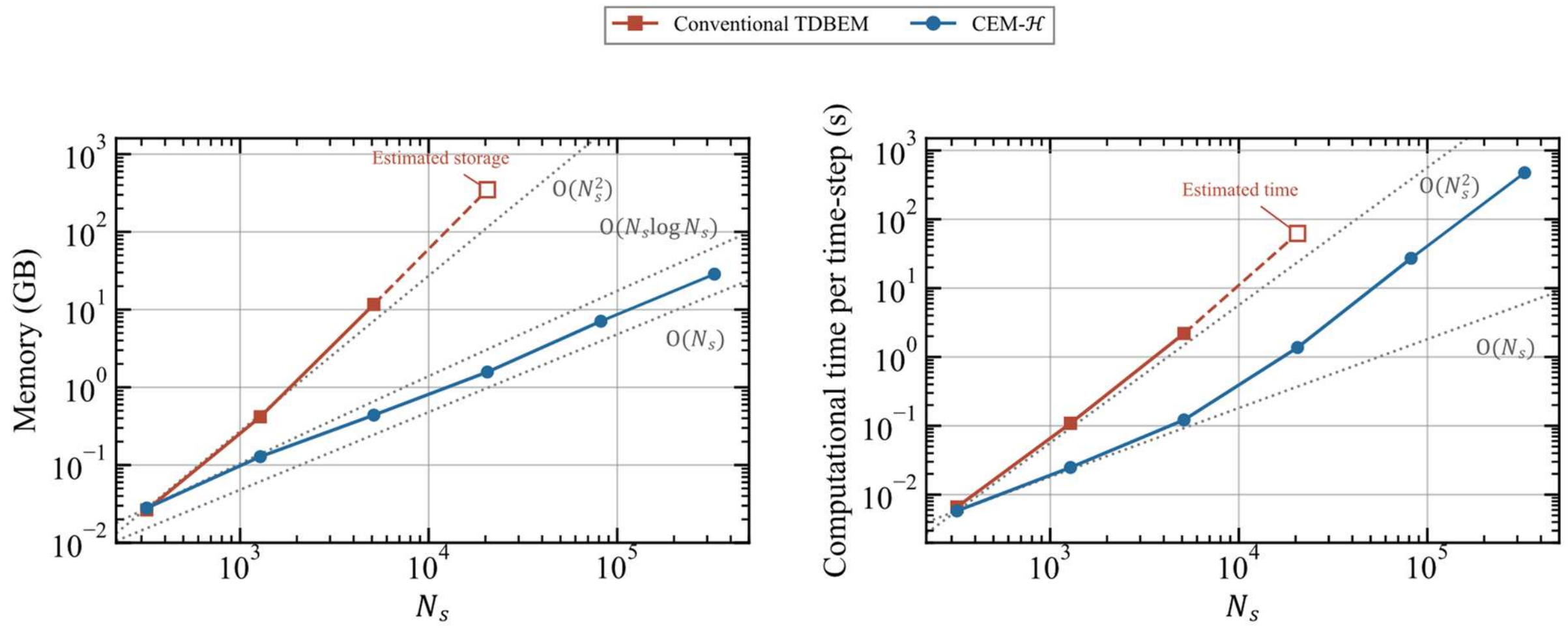

(a) Storage

(b) Computational time

Fig. 9 Storage requirements and average computation time per time step for conventional TDBEM and CEM-$\mathcal{H}$

Finally, the acceleration provided by the local dynamic cluster-pair tree during initial wavefront propagation is examined. Fig. 10 shows its additional storage requirement and the time savings relative to computations without this strategy. As discussed in Section 4.4, the additional storage follows an $O\left({N_s}^{3/2}\right)$ trend but remains small relative to the overall storage requirement of the algorithm. For example, at $N_s = 327{,}680$, the dynamic cluster-pair tree requires only 0.5 GB of additional storage. Over the range $N_s = 5{,}120 - 327{,}680$, the relative time savings during initial wavefront propagation are approximately $15\% - 20\%$.

These results indicate that maintaining a coarse representation of unreached regions, performing local subdivision near the wavefront, and merging fully arrived regions reduce unnecessary accesses in history convolution with modest additional storage and update costs. As shown in Fig. 10(b), the reduction in average time per step scales approximately as $O\left({N_s}^2\right)$. These savings apply only to the initial wavefront propagation stage and do not represent the overall speedup.

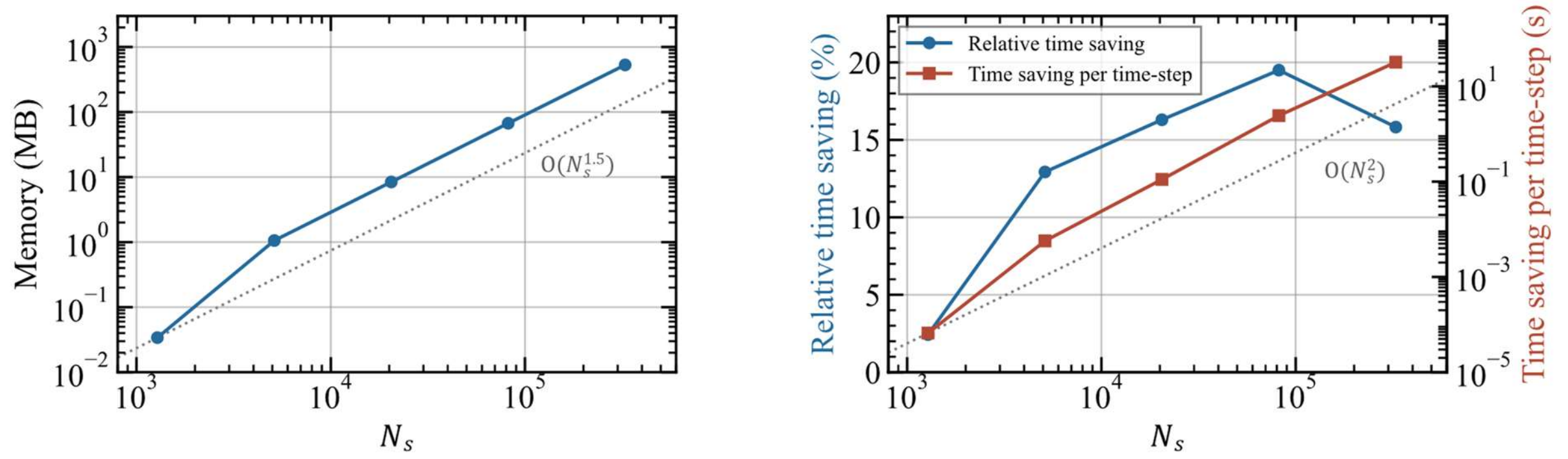

(a) Dynamic tree storage overhead

(b) Time saving during initial wavefront propagation

Fig. 10 Additional storage requirements of the local dynamic cluster-pair tree and time savings during initial wavefront propagation

In summary, for the investigated cases with $\beta$ of order unity, CEM-$\mathcal{H}$ retains accuracy comparable to conventional TDBEM while reducing the storage requirements and computation time of large-scale analyses. The local dynamic cluster-pair tree further improves computational efficiency during initial wavefront propagation.

b. Low-frequency, long-wavelength analysis using CEM-$\mathcal{H}^2$

To assess the accuracy and efficiency of CEM-$\mathcal{H}^2$ in the low-frequency, long-wavelength regime, the incident excitation is changed to a smoothly switched-on sinusoidal plane wave propagating in the positive $x$-direction, as shown in Fig. 11. The incident pressure is $p^{\text{in}}(\boldsymbol{x}, t) = As(\tau)sin(2\pi c\tau/\lambda)$, where $A = D = 1$ and $\lambda = 20D$. The smooth switch-on function is

$$s(\tau) = \begin{cases} 0, & \tau < 0 \\ \frac{1}{2}\left[1 - \cos\left(\frac{\pi\tau}{40}\right)\right], & 0 \le \tau < 40 \\ 1, & \tau \ge 40 \end{cases} \tag{42}$$

The initial acoustic field is zero, and the scattered field satisfies the outgoing radiation condition at infinity. The total analysis time is $T = 160$. In this case, the sphere is much smaller than the incident wavelength, corresponding to a typical low-frequency, long-wavelength regime. The incident wave has a period of 20, the switch-on stage lasts two periods, and the full analysis interval covers eight periods.

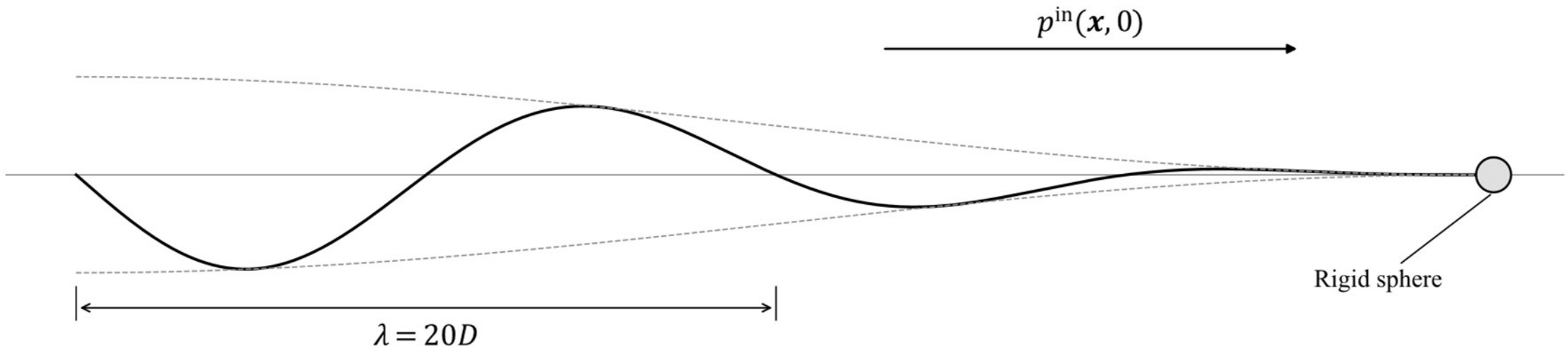

Fig. 11 Schematic of scattering by a rigid sphere under low-frequency, long-wavelength incidence, $\lambda = 20D$

Following the temporal oversampling settings used in the multilevel Cartesian non-uniform grid time domain algorithm (CNGTDA) [15] and the ACE approach of characterizing interaction scales through the domain size and the distance traveled within one time step [19], we introduce the dimensionless propagation-step count $N_D = D/(c\Delta t)$, which represents the number of time steps required for sound to travel one sphere diameter. Values of $N_D = 1,2,4$ are considered to assess the accuracy and efficiency of CEM-$\mathcal{H}^2$ at different time-step sizes. The corresponding parameters are listed in Tab. 5. Since $\lambda = 20D$, these values of $N_D$ give $20, 40$, and $80$ time steps per incident-wave period, respectively, providing sufficient sampling of the sinusoidal carrier [15]. Thus, low-frequency excitation allows relatively large time steps to be used while maintaining adequate temporal resolution. The distance traveled within one time step is then of the same order as the sphere diameter, yielding a computational regime with only a small number of nonzero time-lag matrices.

Tab. 5 Temporal discretization parameters for the low-frequency, long-wavelength sphere example

| $N_D$ | $\Delta t$ | $N_t$ | Time steps per period |
|---|---|---|---|
| 1 | 1 | 160 | 20 |
| 2 | 0.5 | 320 | 40 |
| 4 | 0.25 | 640 | 80 |

Fig. 12 compares the errors of conventional TDBEM and CEM-$\mathcal{H}^2$ relative to the semi-analytical solution, and Tab. 6 lists the corresponding total computation times. For CEM-$\mathcal{H}^2$, the spatial interpolation-node parameter is set to $p_s = 2, 4, 6$. With $p_s = 2$, the computation time is relatively low, but the error increases with the number of elements. This behavior is mainly attributable to the smaller proportion of far-field interactions on coarse meshes, where interpolation has a limited effect. As the mesh is refined, a larger proportion of far-field interactions is evaluated by interpolation, and the insufficient accuracy of low-order interpolation affects the overall solution.

For $p_s = 4$ or $6$, the errors are close to those of the conventional method. Reducing the time-step size significantly decreases the error but also increases the computational cost. The accuracy obtained with $p_s = 4$ is comparable to that with $p_s = 6$, while the computation time is substantially lower. Therefore, within the parameter range investigated in this example, $p_s = 4$ provides a favorable balance between accuracy and efficiency.

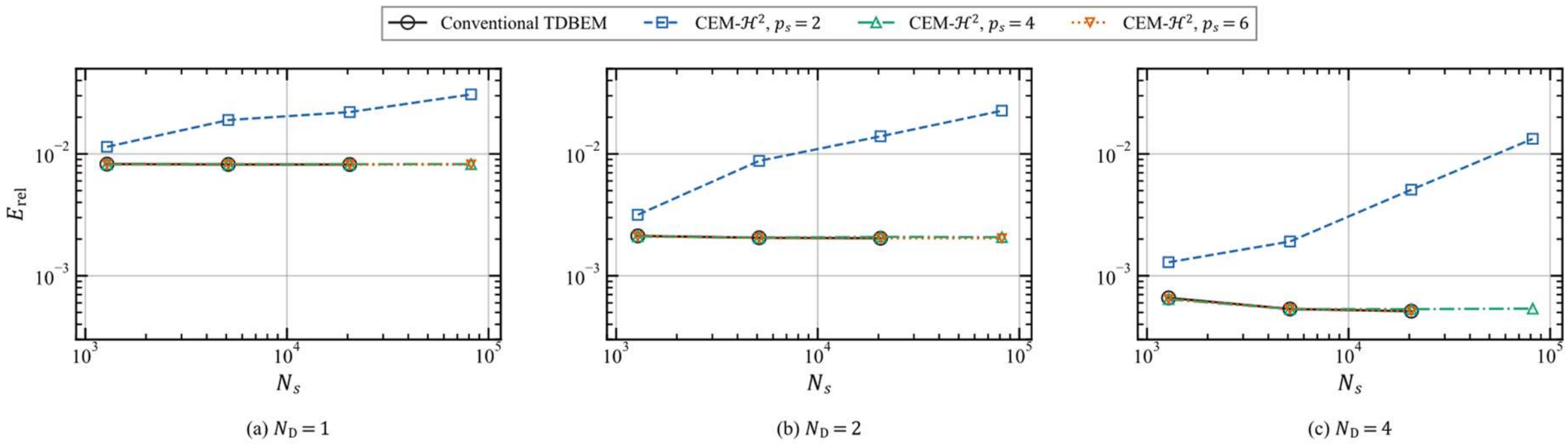

Fig. 12 Relative errors for low-frequency scattering with different time-step parameters and numbers of spatial interpolation nodes

Tab. 6 Total computation times (s) for different propagation-step counts and numbers of spatial interpolation nodes

| $N_s$ | $N_D$ | Conventional TDBEM | CEM-$\mathcal{H}^2$ $p_s = 2$ | $p_s = 4$ | $p_s = 6$ |
|---|---|---|---|---|---|
| | 1 | 3.2 | 2.4 | 6.1 | 13.1 |
| 1,280 | 2 | 5.5 | 2.1 | 3.1 | 5.5 |
| | 4 | 11.2 | 4.3 | 6.3 | 11.3 |
| | 1 | 32.1 | 5.4 | 12.3 | 38.3 |
| 5,120 | 2 | 53.9 | 25.3 | 28.8 | 75.7 |
| | 4 | 121.0 | 37.6 | 44.3 | 115.0 |
| | 1 | 460.2 | 25.1 | 48.0 | 171.4 |
| 20,480 | 2 | 790.4 | 165.3 | 230.3 | 571.3 |
| | 4 | 1,643.5 | 370.6 | 440.6 | 890.3 |
| | 1 | — | 97.2 | 194.2 | 615.8 |
| 81,920 | 2 | — | 1,160.2 | 1,464.9 | 3,202.8 |
| | 4 | — | 2,925.9 | 3,699.0 | — |

Fig. 13 presents the total-pressure time histories at the front $\theta = \pi$ and rear $\theta = 0$ of the sphere for $N_s = 20{,}480$, together with the semi-analytical solution. The propagation-step count is set to $N_D = 1, 4$, and the spatial interpolation-node parameter to $p_s = 2, 6$, to examine the effects of the time-step size and spatial interpolation accuracy. Increasing $N_D$, which corresponds to reducing the time-step size, and appropriately increasing $p_s$ both reduce the deviation from the semi-analytical solution and improve the accuracy of the transient response.

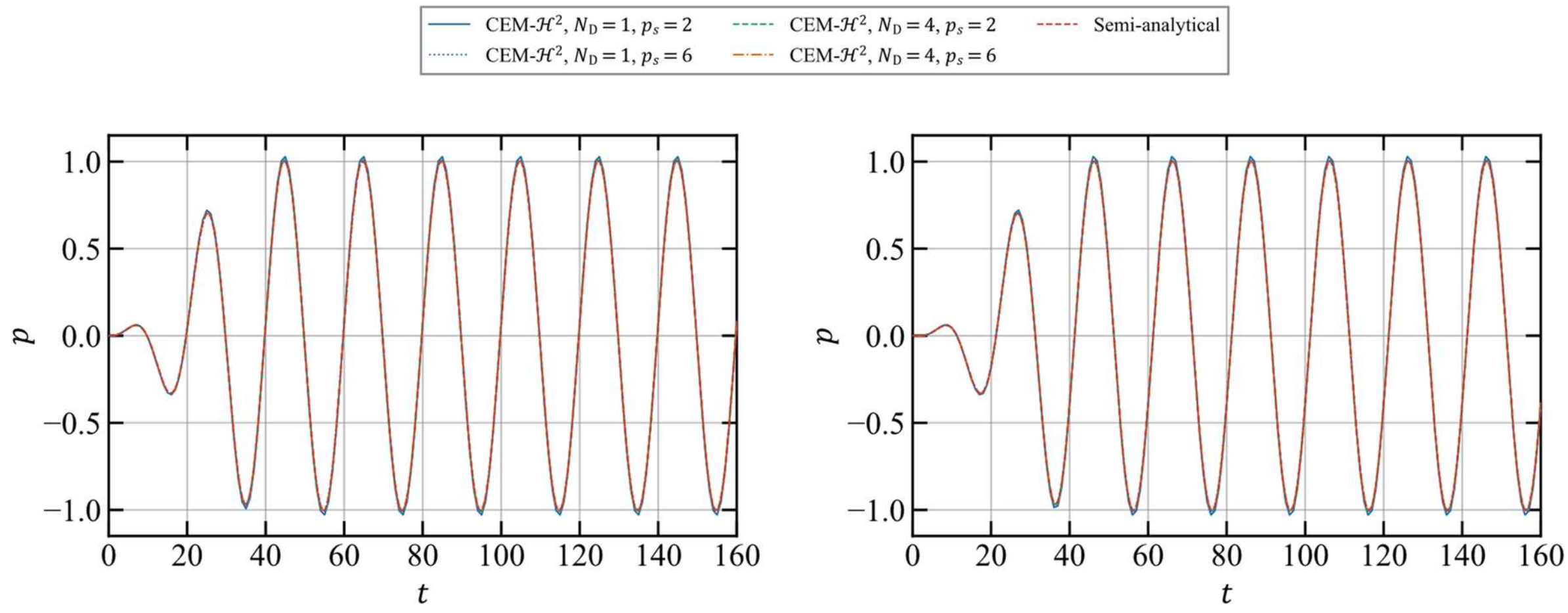

(a) $\theta = \pi$ (b) $\theta = 0$

Fig. 13 Pressure time histories at the front and rear of the sphere for different time-step parameters and numbers of spatial interpolation nodes

To assess large-scale computational performance, $N_D = 1$ and $p_s = 4$ are held fixed while the number of elements is progressively increased. Fig. 14 shows the storage requirements and average computation time per time step. For small problems, the additional overhead of CEM-$\mathcal{H}^2$ limits its advantage in computation time. As the problem size increases, its storage and computational advantages become more pronounced. Overall, both the storage requirement and computation time for large problems exhibit $O(N_s)$ growth. At $N_s = 1{,}310{,}720$, the peak memory usage is approximately $100.5$ GB, and the average computation time per time step is approximately $17.8$ s. These results demonstrate the capability of CEM-$\mathcal{H}^2$ to handle low-frequency transient scattering problems with more than one million boundary elements.

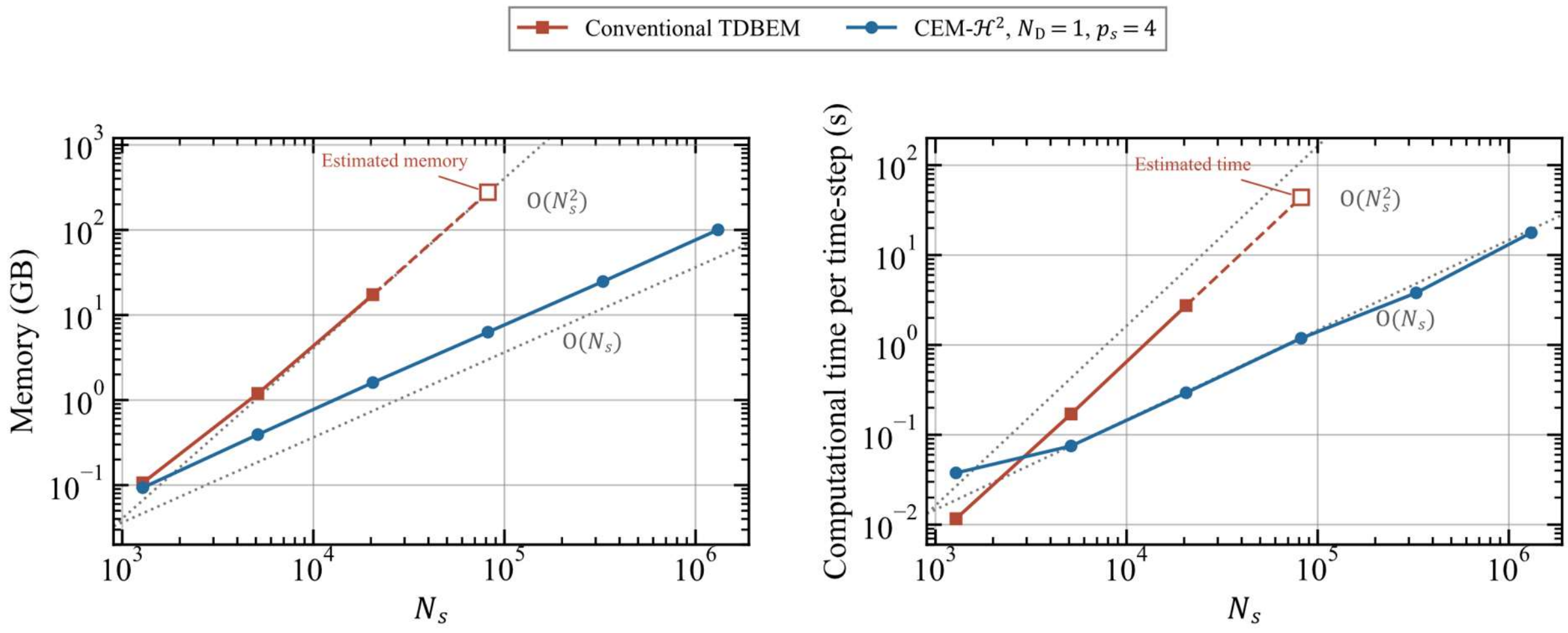

(a) Storage (b) Computational time

Fig. 14 Storage requirements and average computation time per time step for conventional TDBEM and CEM-$\mathcal{H}^2$ in the low-frequency, long-wavelength regime

## 6.2 Underwater vehicle

To further assess the applicability of the proposed methods to complex geometries, an acoustic scattering model of an underwater vehicle is constructed using the publicly available BeTSSi benchmark geometry [41], uniformly scaled by a factor of $1:50$. The scaled model has an overall length of $l = 1.24$ m. The surrounding fluid is treated as an unbounded homogeneous medium with sound speed $c = 1500\ \mathrm{m/s}$ and density $\rho = 10^3\ \mathrm{kg/m^3}$. The vehicle surface satisfies the rigid boundary condition $q = \partial p/\partial n = 0$. The initial acoustic field is zero, and the scattered field satisfies the outgoing radiation condition at infinity.

The incident plane wave propagates in the positive $x$-direction and is modulated by a Gaussian envelope:

$$p^{\mathrm{in}}(\boldsymbol{x}, t) = \begin{cases} p_0 \cos(2\pi f_0 \tau) \exp\left[-\frac{(\tau - t_0)^2}{2\sigma_t^2}\right], & 0 \le \tau \\ 0, & 0 > \tau \end{cases} \tag{43}$$

where $\tau = t - x/c$, $p_0 = 1$ Pa, $f_0$ is the modulation frequency, $\sigma_t$ controls the temporal width of the envelope, and $t_0$ is the time at which the envelope is centered. These parameters are adjusted to construct conventional-frequency and low-frequency, long-wavelength cases, allowing computational behavior at different propagation scales to be examined. In addition to the surface pressure distribution, the scattered-pressure time history is evaluated at an exterior observation point $P = (1.5,0,0)$ on the shadow side of the vehicle. CPU Hyper-Threading is enabled for the examples in this subsection, and all logical cores are used, giving a total of 40 OpenMP threads. The computational model and the incident-wave configuration are shown in Fig. 15.

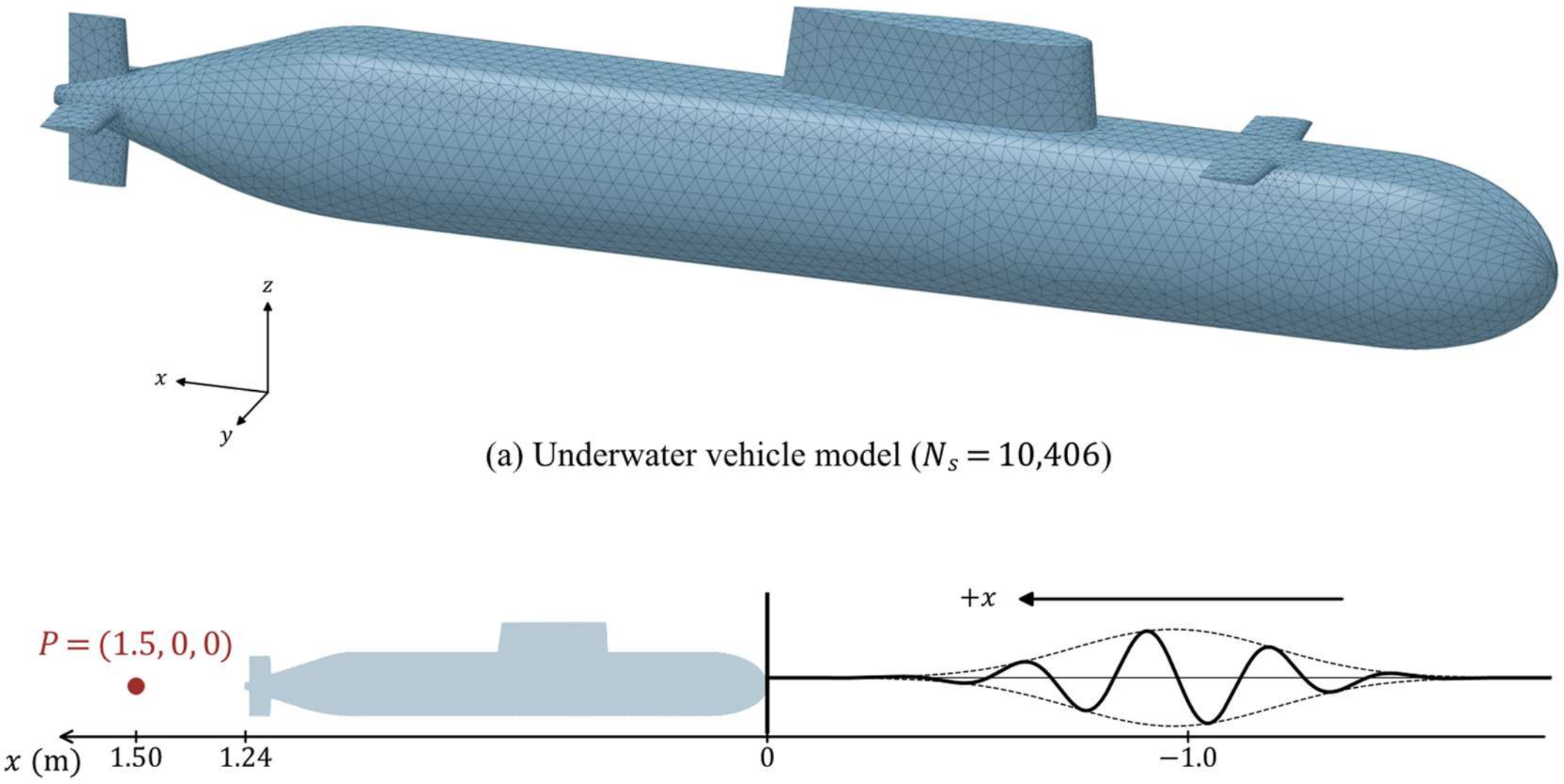

(a) Underwater vehicle model ($N_s = 10{,}406$)

(b) Incident wave at $t = 0$ and observation point ($f_0 = 5$ kHz)

Fig. 15 Underwater vehicle model and schematic of the incident-wave configuration

a. Conventional-frequency case

The center frequency is set to $f_0 = 5$ kHz. Using 1% of the peak spectral amplitude as the bandwidth criterion gives an effective maximum frequency of $f_{\max} = 8$ kHz. The Gaussian-envelope parameters are $\sigma_t = 1.61 \times 10^{-4}$ s and $t_0 = 4\sigma_t = 6.44 \times 10^{-4}$ s. The model surface is discretized into $N_s = 106{,}888$ constant triangular boundary elements. The time-step size is $\Delta t = 2\ \mu\mathrm{s}$, and the total number of time steps is $N_t = 1100$. The computation is performed using CEM-$\mathcal{H}$.

Fig. 16 shows the surface total-pressure distributions at $t = 0.7$, $1.1$, $1.3$ and $1.5$ ms. As the incident pulse propagates along the vehicle, the region of pronounced pressure response moves progressively from the illuminated end toward the stern. At $t = 0.7$ ms, alternating positive and negative pressures are concentrated mainly near the bow and forward section, while the stern response remains weak. By $t = 1.1$ ms, the principal response region has shifted toward the midsection and the sail, and the pressure near the bow has decreased. The response subsequently propagates farther toward the stern and is concentrated mainly around the aft section and tail fins at $t = 1.5$ ms. The alternating positive and negative pressure distributions are consistent with the oscillatory character of the Gaussian-modulated cosine pulse, and the movement of the response region follows the incident-wave propagation direction. After the pulse passes, the pressure over most of the surface gradually decreases.

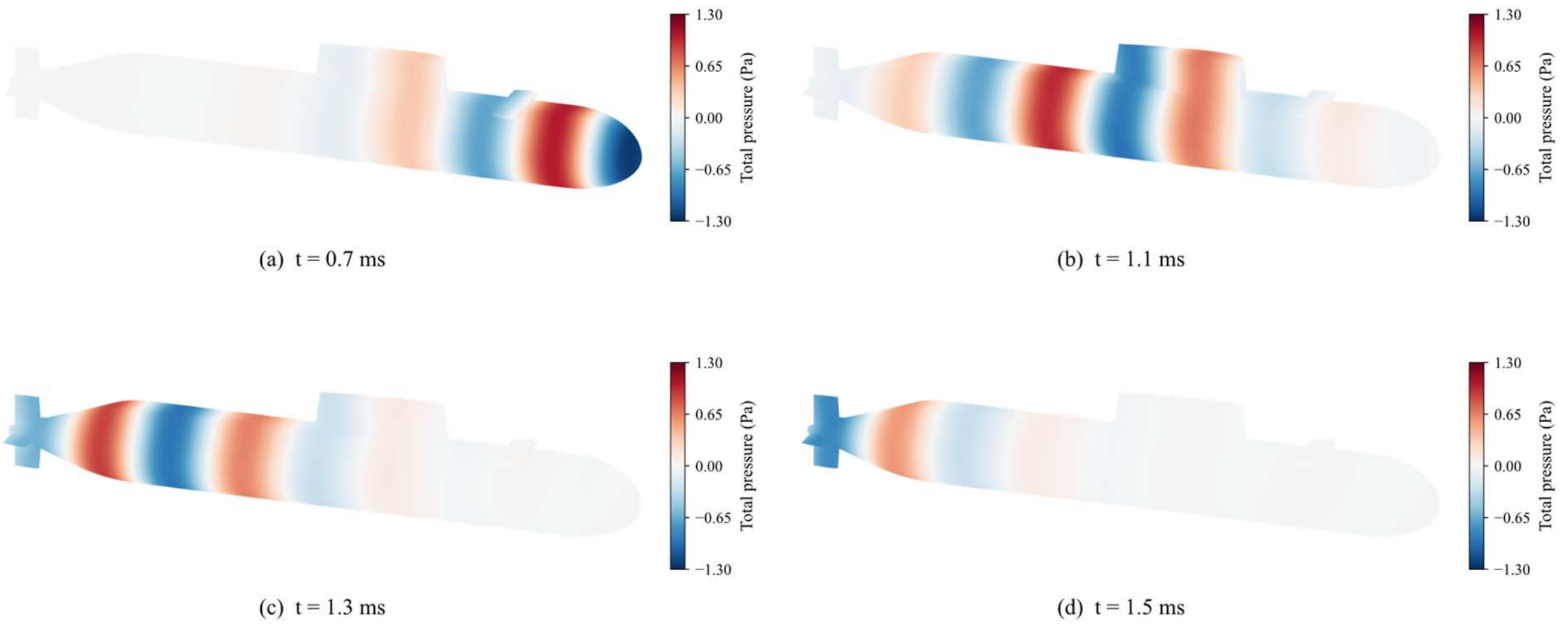

Fig. 16 Transient surface total-pressure distributions on the underwater vehicle under conventional-frequency incidence

Fig. 17 presents the scattered-pressure time history at observation point $P$, together with conventional TDBEM results obtained using a smaller discretization as a reference. The two methods agree well in the waveform during the main response stage, the timing of the peaks, and the subsequent decay. The response at $P$ is initially close to zero and then exhibits alternating positive and negative oscillations. The main response occurs approximately within $1.3 - 2.0$ ms, with the minimum and maximum occurring at $t = 1.56$ ms and $t = 1.64$ ms, respectively. The amplitude gradually decreases after the main response stage. Together with the surface pressure distributions, these results capture the propagation of the pulse along the structure, the development of an exterior response on the shadow side, and its subsequent decay.

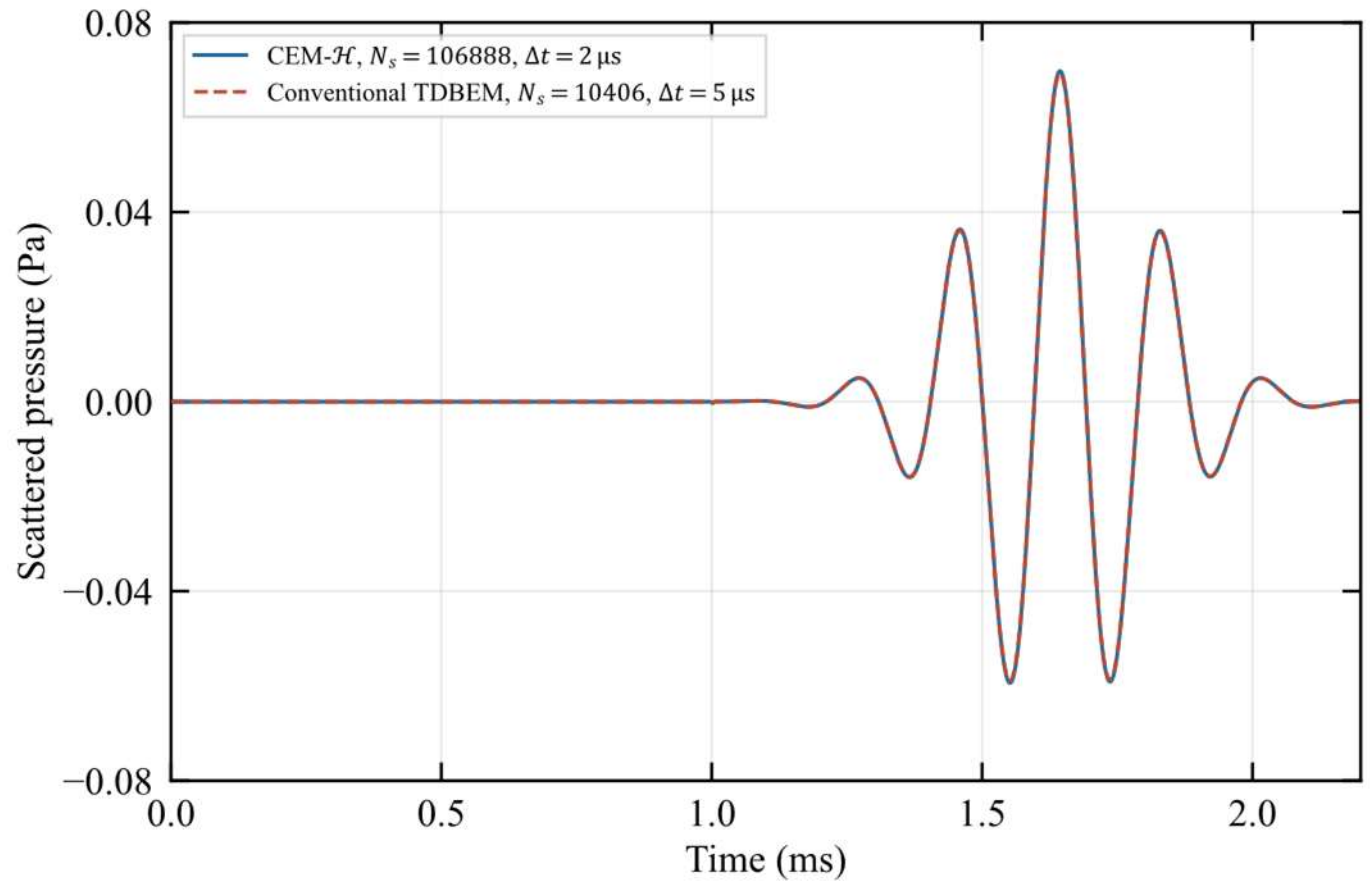

Fig. 17 Scattered-pressure time history at observation point $P$ under conventional-frequency incidence

For this example, CEM-$\mathcal{H}$ requires $15.7$ GB of memory and a total computation time of $76{,}148$ s. It should be noted that this example and the low-frequency, long-wavelength example use different time-step sizes and propagation conditions. Their computation times are jointly influenced by the number of time steps, the number of nonzero time-lag matrices, and the organization of matrix-vector operations.

b. Low-frequency, long-wavelength case

The geometry, medium parameters, boundary conditions, and observation-point location are unchanged, while the incident excitation is replaced by a low-frequency Gaussian pulse. Setting $f_0 = 0$ and taking $f_{\max} = 50$ Hz gives the envelope parameters $\sigma_t = 9.66 \times 10^{-3}$ s and $t_0 = 4\sigma_t = 3.86 \times 10^{-2}$ s. Under the same spectral cutoff criterion, the shortest effective wavelength is $\lambda_{\min} = c/f_{\max} = 30$ m, which is much larger than the model dimensions.

The surface is discretized into $N_s = 1{,}584{,}014$ constant triangular boundary elements. The time-step size is $\Delta t = 1$ ms, and the total number of time steps is $N_t = 100$. The distance traveled within one time step is $c\Delta t = 1.5$ m, exceeding the maximum spatial distance across the model. CEM-$\mathcal{H}^2$ is used with $p_s = 4$.

For this example, CEM-$\mathcal{H}^2$ requires 91.9 GB of memory and a total computation time of 3,362 s. Fig. 18 shows the surface total-pressure distributions during the rising stage, near the peak, and during the decay of the Gaussian pulse. At $t = 20$ and $t = 30$ ms, the total pressure on the vehicle surface increases, with higher pressure near the bow than near the stern. At $t = 40$ ms, the pressure over most of the surface is close to 1 Pa, and the variation along the vehicle is small. By $t = 55$ ms, the overall pressure has decreased substantially, and the pressure near the stern is slightly higher than that near the bow.

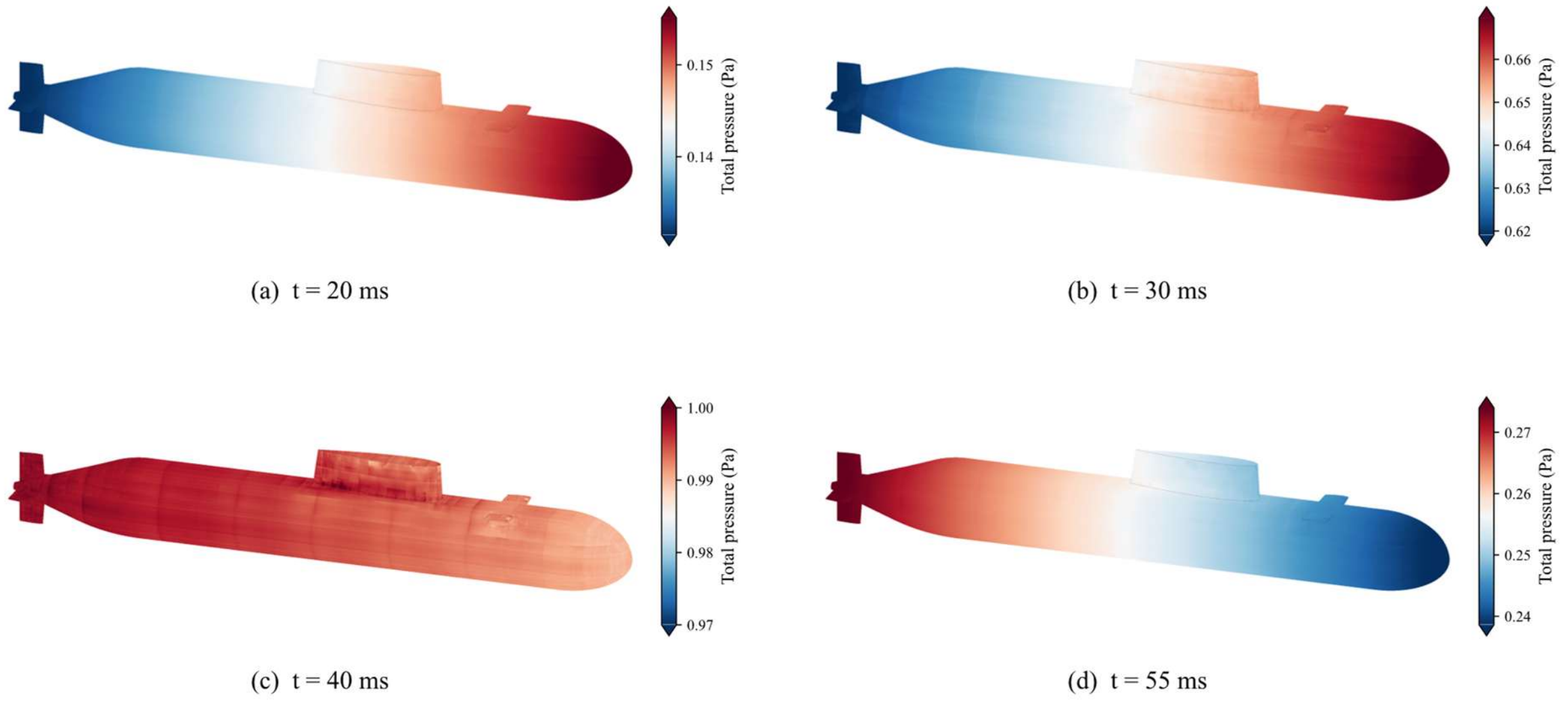

Fig. 18 Transient surface total-pressure distributions on the underwater vehicle under low-frequency Gaussian-pulse incidence

Fig. 19 presents the scattered-pressure time history at observation point $P$ under low-frequency Gaussian-pulse incidence and compares it with conventional TDBEM results obtained using a smaller discretization. The two time histories nearly coincide, with consistent timings of the extrema and overall evolution of the bipolar response. The smaller-scale conventional TDBEM results are included as a reference to demonstrate that the proposed method remains consistent with the conventional method for the complex geometry.

The scattered pressure at $P$ exhibits a bipolar response with a negative lobe followed by a positive lobe. The minimum is approximately $-1.127 \times 10^{-4}$ Pa, occurring at $t = 30$ ms, and the maximum is approximately $1.128 \times 10^{-4}$ Pa, occurring at $t = 49$ ms. For an incident amplitude coefficient of 1 Pa, the scattered pressure at this point is of the order of $10^{-4}$ Pa, substantially smaller than the response in the conventional-frequency case. This indicates a weak scattering response at the selected observation point under the present low-frequency conditions.

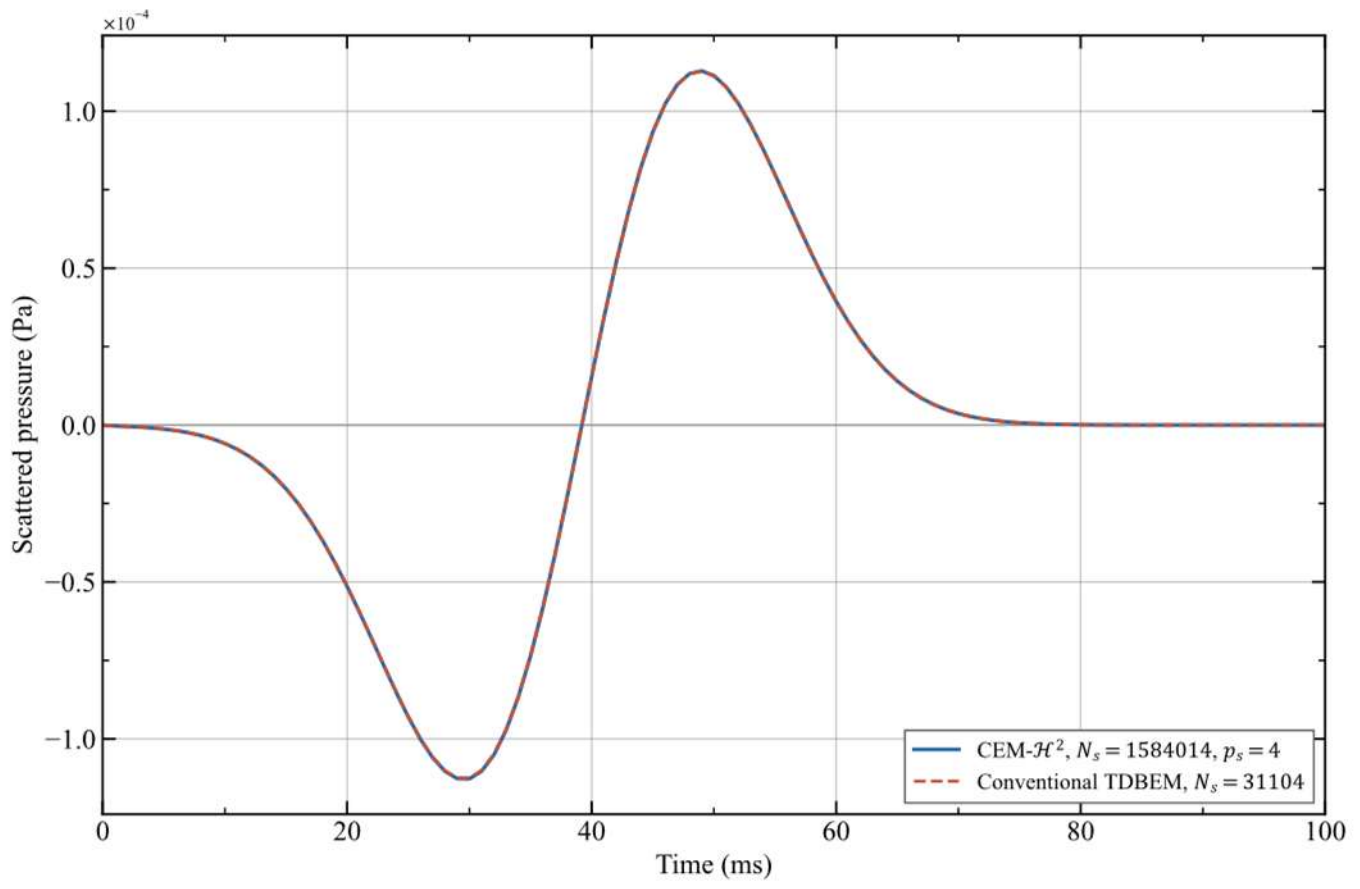

Fig. 19 Scattered-pressure time history at observation point $P$ under low-frequency Gaussian-pulse incidence

# 7 Conclusion

The principal contribution of this work is a time-independent equivalent matrix representation that separates compressible spatial interactions from propagation-delay-dependent historical data access in direct time-marching TDBEM. Based on this representation, two implementations, CEM-$\mathcal{H}$ and CEM-$\mathcal{H}^2$, are developed for different propagation regimes. The main conclusions are summarized as follows.

(1) The far-field retarded interactions are reformulated using three time-independent spatial equivalent matrices, together with the corresponding retarded boundary quantities. This enables a common compressed spatial representation to be reused across propagation delays without additional temporal kernel interpolation. Although the pair-dependent retarded boundary quantities prevent direct application of standard matrix-vector products, the proposed framework provides a basis for incorporating hierarchical spatial compression into time-domain history evaluation. The local error introduced by the centroid-based delay approximation is governed by the propagation-time variation within an element and the temporal variation of the boundary quantities.

(2) For conventional space-time discretizations, CEM-$\mathcal{H}$ combines PACA compression with a local dynamic cluster-pair tree to reduce storage and unnecessary history-data accesses. In the rigid-sphere example with $N_s = 20{,}480$, the storage requirement is reduced from an estimated $349.0$ GB to $1.6$ GB, with an approximately $45.4\times$ speedup relative to the extrapolated conventional computation time. At $N_s = 327{,}680$, the computation requires $28.6$ GB of memory, while the dynamic cluster-pair tree reduces computation time during initial wavefront propagation by approximately $15\% - 20\%$.

(3) For low-frequency, long-wavelength problems, CEM-$\mathcal{H}^2$ combines delay-consistent partitioning with nested cluster bases and shared coupling matrices to enable hierarchical evaluation of retarded interactions. Under the stated sparsity and bounded-parameter assumptions, both storage and computational cost per time step scale as $O(N_s)$, excluding preprocessing. The rigid-sphere example with $N_s = 1{,}310{,}720$ requires approximately $100.5$ GB of memory and $17.8$ s per time step. The underwater-vehicle examples further demonstrate applicability to complex geometries, with CEM-$\mathcal{H}^2$ handling up to $N_s = 1{,}584{,}014$ using $91.9$ GB of memory.

Future work will focus on two aspects. First, the long-time numerical stability and error accumulation of the proposed methods will be systematically investigated, with particular attention to the effects of time-step size, propagation-delay approximations, and spatial compression accuracy. Second, the time-independent equivalent matrix framework will be extended to other transient wave problems, particularly elastodynamics and electromagnetics, where multiple wave speeds, vector-valued kernels, and more complex boundary integral operators introduce additional challenges.

## Acknowledgement

The study was supported by the Beijing Natural Science Foundation (No. 1264055, No. 1252007), the National Key R&D Program of China (No. 2024YFC3017000), and the National Natural Science Foundation of China (No. 72334003).

# Appendix A. Derivation of far-field equivalent matrices for 3D acoustic TDBEM

The derivation of the single-layer equivalent matrix $\boldsymbol{S}$ has been presented in Section 3.3. This appendix provides the detailed derivation of the double-layer equivalent matrices $\boldsymbol{D}_0$ and $\boldsymbol{D}_1$.

Consider the time-domain double-layer integral term on the left-hand side of Eq. (2). Substituting $q^*(\boldsymbol{x},t;\boldsymbol{y},\tau)$ from Eq. (15) into the corresponding temporal convolution and introducing $t_r = t - r/c$ gives

$$\int_0^t q^*(\boldsymbol{x},t;\boldsymbol{y},\tau)p(\boldsymbol{y},\tau)\,d\tau = \frac{(x-y)\cdot\boldsymbol{n}_y}{4\pi}\int_0^t \left[\frac{\delta(t_r-\tau)}{r^3} + \frac{\delta'(t_r-\tau)}{cr^2}\right] p(\boldsymbol{y},\tau)d\tau \tag{A-1}$$

For the term involving the Dirac delta function $\delta(\cdot)$,

$$\int_0^t p\,(\boldsymbol{y},\tau)\delta(t_r-\tau)\,d\tau = p(\boldsymbol{y},t_r) \tag{A-2}$$

For the term involving its derivative $\delta'(\cdot)$, using $\delta'(t_r-\tau) = -\partial\delta(t_r-\tau)/\partial\tau$ and integrating by parts yields

$$\begin{aligned}\int_0^t p\,(\boldsymbol{y},\tau)\delta'(t_r-\tau)\,d\tau &= -\int_0^t p\,(\boldsymbol{y},\tau)\frac{\partial}{\partial\tau}\delta(t_r-\tau)\,d\tau \\ &= -[p(\boldsymbol{y},\tau)\,\delta(t_r-\tau)]_{\tau=0}^{\tau=t} + \int_0^t \dot{p}\,(\boldsymbol{y},\tau)\delta(t_r-\tau)\,d\tau \\ &= -p(\boldsymbol{y},t)\,\delta\left(-\frac{r}{c}\right) + p(\boldsymbol{y},0)\,\delta(t_r) + \int_0^t \dot{p}\,(\boldsymbol{y},\tau)\,\delta(t_r-\tau)\,d\tau \\ &= \dot{p}(\boldsymbol{y},t_r)\end{aligned} \tag{A-3}$$

Since $r > 0$, $\delta(-r/c) = 0$. The zero initial conditions also give $p(\boldsymbol{y},0) = 0$. With the boundary pressure defined to vanish at negative times, corresponding to the period before wavefront arrival, $p(\boldsymbol{y},t_r) = \dot{p}(\boldsymbol{y},t_r) = 0$ when $t - r/c < 0$.

Substituting Eqs. (A-2) and (A-3) into Eq. (A-1) gives

$$\int_0^t q^*(\boldsymbol{x},t;\boldsymbol{y},\tau)p(\boldsymbol{y},\tau)\,d\tau = \frac{(x-y)\cdot\boldsymbol{n}_y}{4\pi}\left[\frac{p(y,t-r/c)}{r^3} + \frac{\dot{p}(y,t-r/c)}{cr^2}\right] \tag{A-4}$$

Thus, the double-layer contribution of the $j$-th boundary element $\Gamma_j$ to the $i$-th collocation equation is

$$D_{ij}p_j = \int_{\Gamma_j}\frac{(x_i-y)\cdot\boldsymbol{n}_j}{4\pi}\left[\frac{p_j(t-r_{ij}(y)/c)}{r_{ij}^3(y)} + \frac{\dot{p}_j(t-r_{ij}(y)/c)}{c\,r_{ij}^2(y)}\right]d\Gamma_y \tag{A-5}$$

For a far-field matrix index pair $(i,j) \in \mathcal{P}_\mathrm{f}$, the propagation delay is determined uniformly using $R_{ij}$. As in Eq. (19), the historical pressure and its time derivative can then be taken outside the spatial integrals:

$$D_{ij}p_j \approx \left[\int_{\Gamma_j}\frac{(x_i-y)\cdot\boldsymbol{n}_j}{4\pi r_{ij}^3(y)}d\Gamma_y\right]p_j\left(t-\frac{R_{ij}}{c}\right) + \frac{1}{c}\left[\int_{\Gamma_j}\frac{(x_i-y)\cdot\boldsymbol{n}_j}{4\pi r_{ij}^2(y)}d\Gamma_y\right]\dot{p}_j\left(t-\frac{R_{ij}}{c}\right) \tag{A-6}$$

This yields the double-layer equivalent matrices $\boldsymbol{D}_0$ and $\boldsymbol{D}_1$, whose expressions are given in Eq. (21). With these definitions, the far-field double-layer contribution can be written compactly as

$$D_{ij}p_j \approx D_{0,ij}p_j\left(t-\frac{R_{ij}}{c}\right) + \frac{D_{1,ij}}{c}\dot{p}_j\left(t-\frac{R_{ij}}{c}\right) \tag{A-7}$$

where $D_{0,ij}$ corresponds to the retarded-pressure term, while $D_{1,ij}$ corresponds to the retarded pressure time-derivative term. The factor $1/c$ is not included in the definition of $D_{1,ij}$, so both equivalent matrices depend only on the spatial geometry and integration kernels.

# Appendix B. Error analysis of the far-field equivalent matrices

B.1 Error in propagation-delay evaluation

For any far-field matrix index pair $(i,j) \in \mathcal{P}_\mathrm{f}$, define the propagation time from the collocation point $\boldsymbol{x}_i$ to an arbitrary integration position $\boldsymbol{y} \in \Gamma_j$, together with the centroid-based propagation time, as follows:

$$\tau_{ij}(\boldsymbol{y}) = \frac{r_{ij}(\boldsymbol{y})}{c}$$
$$\tau_{ij}^{c} = \frac{R_{ij}}{c} \tag{B-1}$$

For the $g$-th quadrature point $\boldsymbol{y}_{jg}$, let $r_{ijg} = r_{ij}(\boldsymbol{y}_{jg})$. Its propagation-time deviation from the centroid-based value is

$$\delta\tau_{ijg} = \left|\tau_{ij}(\boldsymbol{y}_{jg}) - \tau_{ij}^{c}\right| = \frac{|r_{ijg} - R_{ij}|}{c} \tag{B-2}$$

The element-wise theoretical maximum propagation-time deviation is defined as

$$\delta\tau_{ij} = \max_{\boldsymbol{y}\in\Gamma_j}\left|\tau_{ij}(\boldsymbol{y}) - \tau_{ij}^{c}\right| = \frac{1}{c}\max_{\boldsymbol{y}\in\Gamma_j}\left|r_{ij}(\boldsymbol{y}) - R_{ij}\right| \tag{B-3}$$

Consequently, $\delta\tau_{ijg} \le \delta\tau_{ij},\ g = 1, \dots, N_g$.

For a specific quadrature rule, one may also define the realized quadrature-point maximum $\delta\tau_{ij}^{Q} = \max_{1\le g\le N_g} \delta\tau_{ijg}$, which satisfies $\delta\tau_{ij}^{Q} \le \delta\tau_{ij}$. The derivation below uses the element-wise theoretical bound $\delta\tau_{ij}$, so that the resulting estimate is independent of the particular quadrature-point distribution. The corresponding dimensionless measure is

$$\eta_{ij}^{\tau} = \frac{\delta\tau_{ij}}{\Delta t} \tag{B-4}$$

This quantity measures the maximum element-wise propagation-time deviation relative to one time step. Suppose that the underlying historical boundary quantity $\phi_j(t)$ is absolutely continuous over the relevant local retarded-time interval and has an essentially bounded first time derivative. For the $g$-th quadrature point, the fundamental theorem of calculus gives

$$\phi_j\left(t - \frac{r_{ijg}}{c}\right) - \phi_j\left(t - \frac{R_{ij}}{c}\right) = \int_{t-R_{ij}/c}^{t-r_{ijg}/c} \dot{\phi}_j(\xi)\mathrm{d}\xi \tag{B-5}$$

Taking absolute values and applying the integral triangle inequality yields

$$\left|\int_{t-R_{ij}/c}^{t-r_{ijg}/c} \dot{\phi}_j(\xi)\mathrm{d}\xi\right| \le \int_{\min(t-R_{ij}/c,\,t-r_{ijg}/c)}^{\max(t-R_{ij}/c,\,t-r_{ijg}/c)} \left|\dot{\phi}_j(\xi)\right|\mathrm{d}\xi \tag{B-6}$$

Let $I_{ij}(t)$ denote the smallest interval containing $t - R_{ij}/c$ and all retarded times $t - r_{ij}(\boldsymbol{y})/c$ with $\boldsymbol{y} \in \Gamma_j$. By the definition of the $L^{\infty}$-norm, $\left|\dot{\phi}_j(\xi)\right| \le \|\dot{\phi}_j\|_{L^{\infty}(I_{ij}(t))}$ almost everywhere on this interval. Combining Eqs. (B-2) and (B-3) with Eq. (B-6) gives

$$\left|\phi_j\left(t - \frac{r_{ijg}}{c}\right) - \phi_j\left(t - \frac{R_{ij}}{c}\right)\right| \le \delta\tau_{ijg}\|\dot{\phi}_j\|_{L^{\infty}(I_{ij}(t))} \le \delta\tau_{ij}\|\dot{\phi}_j\|_{L^{\infty}(I_{ij}(t))} \tag{B-7}$$

Eq. (B-7) provides an upper bound on the local error in the historical boundary quantity introduced by the centroid-based delay approximation. Because $\delta\tau_{ij}$ is defined over the complete source element, this bound applies to every quadrature point used in the spatial integration.

By the reverse triangle inequality, $\left|r_{ij}(\boldsymbol{y}) - R_{ij}\right| \le \|\boldsymbol{y} - \boldsymbol{y}_j^{\mathrm{cen}}\|$, so that $\delta\tau_{ij} \le \rho_j/c$, where $\rho_j = \max_{\boldsymbol{y}\in\Gamma_j} \| \boldsymbol{y} - \boldsymbol{y}_j^{\mathrm{cen}} \|$. Thus, the theoretical propagation-time deviation is controlled directly by the spatial size of the source element.

B.2 Error in far-field element contributions

Let $s_{ijg}$, $d_{0,ijg}$, and $d_{1,ijg}$ denote the quadrature-point coefficients corresponding to the equivalent matrices $\boldsymbol{S}$, $\boldsymbol{D}_0$, and $\boldsymbol{D}_1$, respectively. When propagation delays are evaluated at quadrature points, the discrete integral contribution associated with $(i, j) \in \mathcal{P}_{\mathrm{f}}$ is

$$F_{ij}^{\text{trad}}(t)=\sum_{g=1}^{N_g}\left[s_{ijg}q_j\left(t-\tau_{ijg}\right)-d_{0,ijg}p_j\left(t-\tau_{ijg}\right)-\frac{d_{1,ijg}}{c}\dot{p}_j\left(t-\tau_{ijg}\right)\right] \tag{B-8}$$

With the centroid-based delay approximation, the corresponding contribution becomes

$$F_{ij}^{\text{eq}}(t)=\sum_{g=1}^{N_g}\left[s_{ijg}q_j\left(t-\tau_{ij}^c\right)-d_{0,ijg}p_j\left(t-\tau_{ij}^c\right)-\frac{d_{1,ijg}}{c}\dot{p}_j\left(t-\tau_{ij}^c\right)\right] \tag{B-9}$$

Define $E_{ij}(t)=F_{ij}^{\text{trad}}(t)-F_{ij}^{\text{eq}}(t)$, and assume that the underlying boundary histories satisfy $q_j\in C^1$ and $p_j\in C^2$ over the local retarded-time interval $I_{ij}(t)$. These smoothness assumptions refer to the underlying boundary histories used in the local consistency analysis, rather than to their globally piecewise-linear discrete representations. Consequently, the estimate below should be interpreted as a local smoothness-based bound; it does not by itself provide a bound for derivative jumps associated with the temporal interpolation knots.

Applying Eq. (B-7) separately to $q_j$, $p_j$, and $\dot{p}_j$, and then using the triangle inequality, gives

$$\left|E_{ij}(t)\right|\leq c_{ij}(t)\,\delta\tau_{ij}=c_{ij}(t)\,\eta_{ij}^{\tau}\Delta t \tag{B-10}$$

where

$$\begin{aligned}c_{ij}(t)=&\sum_{g=1}^{N_g}\left|s_{ijg}\right|\,\|\dot{q}_j\|_{L^\infty\left(I_{ij}(t)\right)}+\sum_{g=1}^{N_g}\left|d_{0,ijg}\right|\,\|\dot{p}_j\|_{L^\infty\left(I_{ij}(t)\right)}\\&+\frac{1}{c}\sum_{g=1}^{N_g}\left|d_{1,ijg}\right|\,\|\ddot{p}_j\|_{L^\infty\left(I_{ij}(t)\right)}\end{aligned} \tag{B-11}$$

Eq. (B-10) shows that, provided the boundary quantities are sufficiently smooth over the local retarded-time interval, the local error bound for the quadrature-evaluated far-field element contribution introduced by the centroid-based delay approximation is proportional to the maximum element-wise propagation-time deviation $\delta\tau_{ij}$.

The following conclusions can therefore be drawn.

(1) When the spatial integral coefficients and the required time derivatives of the boundary quantities are bounded, the upper bound on the local element-contribution error introduced by the centroid-based delay approximation is proportional to $\delta\tau_{ij}$. Since $\delta\tau_{ij}\leq\rho_j/c$, the bound decreases as the source-element size decreases.

(2) The error also depends on the temporal variation of the boundary quantities. Faster temporal variation increases the coefficient $c_{ij}(t)$; consequently, high-frequency signals and steep wavefronts are more sensitive to a given propagation-delay deviation.

(3) The dimensionless measure $\eta_{ij}^{\tau}$ characterizes the maximum element-wise propagation-time deviation relative to the time-step size. When $\eta_{ij}^{\tau}\ll 1$ and the boundary quantities vary slowly over the corresponding history interval, the local error introduced by the centroid-based delay approximation is expected to remain small.

It should be emphasized that this appendix analyzes only the local error introduced by the centroid-based delay approximation. Spatial quadrature error and piecewise-linear temporal interpolation error are not included, and no global error accumulation over the fully discrete time-marching procedure is estimated here. The purpose of the analysis is to provide a consistency bound that identifies the principal quantities governing the additional error introduced by replacing quadrature-point propagation delays with the centroid-based delay in the far-field equivalent representation.

# Appendix C. Supporting derivations and implementation details of CEM-$\mathcal{H}^2$

C.1 Interpolation nodes and basis functions

On the reference interval $[-1,1]$, $p_s$ Chebyshev nodes are selected, and the corresponding Lagrange basis functions are defined as

$$\hat{\xi}_a=\cos\frac{(2a-1)\pi}{2p_s} \tag{C-1}$$

$$\ell_a(\hat{x}) = \prod_{\substack{d=1 \\ d\neq a}}^{p_s} \frac{\hat{x}-\hat{\xi}_d}{\hat{\xi}_a-\hat{\xi}_d} \tag{C-2}$$

These functions satisfy $\ell_a(\hat{\xi}_d) = \delta_{ad}$, where $\delta_{ad}$ is the Kronecker delta.

For a spatial cluster $C$, a nondegenerate interpolation bounding box $\tilde{B}_C = \prod_{\nu=1}^{3}\left[x_{\nu,C}^{-}, x_{\nu,C}^{+}\right]$ is used to enclose the collocation points and source-element integration domains within the cluster. Here, $\nu = 1,2,3$ denotes the three coordinate directions, $\prod$ represents the Cartesian product of the three intervals, and $x_{\nu,C}^{-}$ and $x_{\nu,C}^{+}$ are the lower and upper coordinate bounds, respectively. The mapping from physical coordinates to reference coordinates is

$$\hat{x}_\nu = \frac{2(x_\nu - x_{\nu,C}^{-})}{x_{\nu,C}^{+} - x_{\nu,C}^{-}} - 1 \tag{C-3}$$

The $n_s$ interpolation nodes and basis functions are obtained by taking tensor products in three dimensions:

$$[\boldsymbol{\xi}_{\boldsymbol{a}}^{C}]_\nu = x_{\nu,C}^{-} + \frac{\hat{\xi}_{a_\nu}+1}{2}\left(x_{\nu,C}^{+} - x_{\nu,C}^{-}\right) \tag{C-4}$$

$$L_{\boldsymbol{a}}^{C}(\boldsymbol{x}) = \prod_{\nu=1}^{3} \ell_{a_\nu}(\hat{x}_\nu) \tag{C-5}$$

where $\boldsymbol{a} = (a_1, a_2, a_3)$ specifies the interpolation-node indices in the three coordinate directions.

C.2 Two-sided interpolation of time-independent equivalent matrices

Following the definitions of the time-independent equivalent matrices in Section 3.3, this subsection describes the construction of their $\mathcal{H}^2$-matrix representations. The interpolation below is applied only to the spatial kernels.

Consider a far-field block $b$. For notational simplicity, the superscript $b$ on the clusters is omitted below. Within the interpolation regions of the row and column clusters, the smooth scalar spatial kernel $K(\boldsymbol{x},\boldsymbol{y})$ is approximated by two-sided interpolation:

$$K(\boldsymbol{x},\boldsymbol{y}) \approx \sum_{a=1}^{n_s}\sum_{d=1}^{n_s} L_a^{C_c}(\boldsymbol{x}) K\left(\boldsymbol{\xi}_a^{C_c}, \boldsymbol{\xi}_d^{C_s}\right) L_d^{C_s}(\boldsymbol{y}) \tag{C-6}$$

Evaluating at the collocation point $\boldsymbol{x}_i$ and integrating over source element $\Gamma_j$ gives

$$\left[\boldsymbol{A}_{K,b}\right]_{ij} = \int_{\Gamma_j} K(\boldsymbol{x}_i, \boldsymbol{y})\mathrm{d}\Gamma_y \approx \sum_{a=1}^{n_s}\sum_{d=1}^{n_s} L_a^{C_c}(\boldsymbol{x}_i) K\left(\boldsymbol{\xi}_a^{C_c}, \boldsymbol{\xi}_d^{C_s}\right) \int_{\Gamma_j} L_d^{C_s}(\boldsymbol{y})\mathrm{d}\Gamma_y \tag{C-7}$$

Using the row and column cluster bases defined in Section 5.1 and setting $\left[\boldsymbol{G}_{K,b}\right]_{ad} = K\left(\boldsymbol{\xi}_a^{C_c}, \boldsymbol{\xi}_d^{C_s}\right)$ yields $\boldsymbol{A}_{K,b} \approx \boldsymbol{P}_{C_c}\boldsymbol{G}_{K,b}\boldsymbol{Q}_{C_s}^{T}$. For $K(\boldsymbol{x},\boldsymbol{y}) = 1/(4\pi r)$, the corresponding coupling matrix is $\boldsymbol{G}_b$, giving the single-layer equivalent matrix representation $\boldsymbol{S}_b \approx \boldsymbol{P}_{C_c}\boldsymbol{G}_b\boldsymbol{Q}_{C_s}^{\mathrm{T}}$. The column cluster basis is evaluated by integration over the source elements: $[\boldsymbol{Q}_C]_{jd} \approx \sum_g w_{jg} L_d^{C}(\boldsymbol{y}_{jg})$. For the double-layer matrix $\boldsymbol{D}_{0,b}$, the spatial kernel is the source-side normal derivative of the scalar kernel above. Differentiating Eq. (C-6) in this direction, with differentiation acting only on the column-cluster basis functions, and using the normal-derivative column cluster basis $\boldsymbol{Q}_C^{D}$ defined in Section 5.1 gives $\boldsymbol{D}_{0,b} \approx \boldsymbol{P}_{C_c}\boldsymbol{G}_b\left(\boldsymbol{Q}_{C_s}^{D}\right)^{\mathrm{T}}$. The normal-derivative column cluster basis is obtained by integration over the source elements: $[\boldsymbol{Q}_C^{D}]_{jd} \approx \sum_g w_{jg}\boldsymbol{n}_j \cdot \nabla_y L_d^{C}(\boldsymbol{y}_{jg})$. Thus, the single- and double-layer equivalent matrices share the same coupling matrix $\boldsymbol{G}_b$, differing only in their column cluster bases.

C.3 Derivation of the cluster-basis nesting relations

When the same $p_s$ is used at every level, a parent-cluster basis function has degree at most $p_s - 1$ in each coordinate direction. It can therefore be represented exactly using the tensor-product Lagrange basis of a child cluster:

$$L_a^{C_{\mathrm{p}}}(\boldsymbol{x}) = \sum_{d=1}^{n_s} L_a^{C_{\mathrm{p}}}\left(\boldsymbol{\xi}_d^{C_{\mathrm{ch}}}\right) L_d^{C_{\mathrm{ch}}}(\boldsymbol{x}) \tag{C-8}$$

Using the transfer matrix defined in Eq. (33), evaluation at the collocation points of the child cluster and integration over its source elements give, respectively,

$$\left[\boldsymbol{P}_{C_{\rm p}}\right]_{ia} = \sum_{d=1}^{n_s}\left[\boldsymbol{P}_{C_{\rm ch}}\right]_{id}\left[\boldsymbol{E}_{C_{\rm ch}\leftarrow C_{\rm p}}\right]_{da}$$
$$\left[\boldsymbol{Q}_{C_{\rm p}}\right]_{ja} = \sum_{d=1}^{n_s}\left[\boldsymbol{Q}_{C_{\rm ch}}\right]_{jd}\left[\boldsymbol{E}_{C_{\rm ch}\leftarrow C_{\rm p}}\right]_{da} \tag{C-9}$$

where $i,j \in K_{C_{\rm ch}}$. These expressions establish the cluster-basis nesting relations used in Section 5.1. The column-cluster forward transformation uses the transpose of the transfer matrix to aggregate child-cluster auxiliary vectors into the parent cluster, whereas the row-cluster backward transformation uses the transfer matrix to propagate parent-cluster contributions to the child clusters. The normal-derivative column cluster basis $\boldsymbol{Q}_C^D$ satisfies the same relations and therefore shares the parent-child transfer matrices with $\boldsymbol{Q}_C$.

C.4 Derivation of delay-consistent retarded interactions

For a cluster pair $b \in \mathcal{B}_k$, the fractional-delay parameter at quadrature point $g$ is given by $\alpha_{ijg} = r_{ijg}/(c\Delta t) - k$. Using the piecewise-linear temporal interpolation in Eq. (34), the single-layer contribution can be expressed in terms of the equivalent matrices $\boldsymbol{S}_b$ and $\boldsymbol{S}_{r,b}$ as follows:

$$\boldsymbol{f}_{b,S}^n = \boldsymbol{S}_b\left[(k+1)\boldsymbol{q}_{C_S^b}^{n-k} - k\boldsymbol{q}_{C_S^b}^{n-k-1}\right] + \frac{1}{c\Delta t}\boldsymbol{S}_{r,b}\left[\boldsymbol{q}_{C_S^b}^{n-k-1} - \boldsymbol{q}_{C_S^b}^{n-k}\right] \tag{C-10}$$

where $\boldsymbol{S}_{r,b}$ is the distance-weighted single-layer matrix defined in Section 5.3.

For the double-layer contribution, the distance-weighted equivalent matrix $\boldsymbol{D}_{0,r,b}$ is introduced, with entries defined by

$$D_{0,r,ij} = \int_{\Gamma_j} r_{ij}(\boldsymbol{y})\frac{(x_i-y)\cdot\boldsymbol{n}_j}{4\pi r_{ij}^3(y)}d\Gamma_{\boldsymbol{y}} = D_{1,ij} \tag{C-11}$$

Therefore,

$$\boldsymbol{D}_{0,r,b} = \boldsymbol{D}_{1,b} \tag{C-12}$$

Incorporating the fractional-delay weights into the double-layer contribution and accounting for the pressure time-derivative term yields

$$\begin{aligned}\boldsymbol{f}_{b,D}^n \approx \boldsymbol{D}_{0,b}\left[(k+1)\boldsymbol{p}_{C_S^b}^{n-k} - k\boldsymbol{p}_{C_S^b}^{n-k-1}\right] + \frac{1}{c\Delta t}\boldsymbol{D}_{0,r,b}\left[\boldsymbol{p}_{C_S^b}^{n-k-1} - \boldsymbol{p}_{C_S^b}^{n-k}\right] \\ + \frac{1}{c\Delta t}\boldsymbol{D}_{1,b}\left[\boldsymbol{p}_{C_S^b}^{n-k} - \boldsymbol{p}_{C_S^b}^{n-k-1}\right]\end{aligned} \tag{C-13}$$

The second term on the right-hand side arises from the distance-dependent temporal interpolation weights, whereas the third term originates from the pressure time derivative. Since $\boldsymbol{D}_{0,r,b} = \boldsymbol{D}_{1,b}$, these two terms cancel exactly when the same quadrature-point distances and common integer-delay interval are used.

Consequently, the double-layer contribution reduces to $\boldsymbol{f}_{b,D}^n \approx \boldsymbol{D}_{0,b}\left[(k+1)\boldsymbol{p}_{C_S^b}^{n-k} - k\boldsymbol{p}_{C_S^b}^{n-k-1}\right]$, as shown in Eq. (36).

For leaf cluster pairs accepted through the centroid-based delay-consistency criterion, the same algebraic cancellation holds within the assigned time interval. The interval-assignment approximation and the associated linear extrapolation are discussed in Section 5.2.

C.5 History-convolution evaluation algorithm

The history-term evaluation procedure using the delay-partitioned $\mathcal{H}^2$-matrix representation is summarized in Algorithm C. The algorithm combines the reuse of transformed historical boundary quantities, intercluster coupling operations, direct evaluation of the remaining far-field interactions, and a final backward transformation on the row cluster tree.

**Algorithm C History-convolution evaluation using $\mathcal{H}^2$-matrices**

**Input**: Delay index tables $\mathcal{G}_k = (\mathcal{B}_k, \mathcal{D}_k)$; current time step $n$, assuming $n > \gamma_{\max}$.

Initialize $\boldsymbol{f}_{\mathrm{f}}^{n} = 0$, $\boldsymbol{f}_{\mathcal{H}^2}^{n} = 0$, $\boldsymbol{f}_{\mathrm{f,dir}}^{n} = 0$.

for $k = 0, \min\{n, \gamma_{\max}\}$ do

  if $\mathcal{B}_k \neq \emptyset$ then

    Retrieve or compute the transformed auxiliary vectors at $n-k$ and $n-k-1$, and cache newly computed results.

    for each $b \in \mathcal{B}_k$ do

      Accumulate the corresponding known contributions into $\boldsymbol{y}_{C_c^b}^{n}$ using Eq. (38).

    end for

  end if

  Evaluate the history contributions associated with $\mathcal{D}_k$ at time steps $n-k$ and $n-k-1$, and accumulate them into $\boldsymbol{f}_{\mathrm{f,dir}}^{n}$.

end for

Perform the backward transformation using Eq. (39).

for each leaf-level row cluster $C_c$ do

  $\boldsymbol{f}_{C_c,\mathcal{H}^2}^{n} = \boldsymbol{P}_{C_c} \boldsymbol{y}_{C_c}^{n}$

end for

$\boldsymbol{f}_{\mathrm{f}}^{n} = \boldsymbol{f}_{\mathcal{H}^2}^{n} + \boldsymbol{f}_{\mathrm{f,dir}}^{n}$

**Output**: Far-field history term $\boldsymbol{f}_{\mathrm{f}}^{n}$ at time step $n$.

# References

[1] A.A. Ergin, B. Shanker, E. Michielssen, Fast transient analysis of acoustic wave scattering from rigid bodies using a two-level plane wave time domain algorithm, Journal of the Acoustical Society of America, 106 (1999) 2405–2416, https://doi.org/10.1121/1.428077.

[2] C. Thirard, J.M. Parot, On a way to save memory when solving time domain boundary integral equations for acoustic and vibroacoustic applications, Journal of Computational Physics, 348 (2017) 744–753, https://doi.org/10.1016/j.jcp.2017.07.040.

[3] H.W. Jang, J.G. Ih, On the instability of time-domain acoustic boundary element method due to the static mode in interior problems, Journal of Sound and Vibration, 332 (2013) 6463–6471, https://doi.org/10.1016/j.jsv.2013.07.018.

[4] Z.Y. Gao, Z.L. Li, Y.J. Liu, A time-domain boundary element method using a kernel-function library for 3D acoustic problems, Engineering Analysis with Boundary Elements, 161 (2024) 103–112, https://doi.org/10.1016/j.enganabound.2024.01.001.

[5] A. Aimi, G. Di Credico, H. Gimperlein, C. Guardasoni, Adaptive time-domain boundary element methods for the wave equation with Neumann boundary conditions, Computers & Mathematics with Applications, 198 (2025) 196–213, https://doi.org/10.1016/j.camwa.2025.08.020.

[6] T. Takahashi, An interpolation-based fast-multipole accelerated boundary integral equation method for the three-dimensional wave equation, Journal of Computational Physics, 258 (2014) 809–832, https://doi.org/10.1016/j.jcp.2013.11.008.

[7] N.A. Gumerov, R. Duraiswami, A broadband fast multipole accelerated boundary element method for the three dimensional Helmholtz equation, Journal of the Acoustical Society of America, 125 (2009) 191–205, https://doi.org/10.1121/1.3021297.

[8] N.A. Gumerov, R. Duraiswami, Fast multipole accelerated boundary element methods for room acoustics, Journal of the Acoustical Society of America, 150 (2021) 1707–1720, https://doi.org/10.1121/10.0006102.

[9] J. Dölz, H. Harbrecht, S. Kurz, S. Schöps, F. Wolf, A fast isogeometric BEM for the three dimensional Laplace- and Helmholtz problems, Computer Methods in Applied Mechanics and Engineering, 330 (2018) 83–101, https://doi.org/10.1016/j.cma.2017.10.020.

[10] F. Amlani, S. Chaillat, A. Loseille, An efficient preconditioner for adaptive Fast Multipole accelerated Boundary Element Methods to model time-harmonic 3D wave propagation, Computer Methods in Applied Mechanics and Engineering, 352 (2019) 189–210, https://doi.org/10.1016/j.cma.2019.04.026.

[11] C. Jelich, W.C. Zhao, H.B. Chen, S. Marburg, Fast multipole boundary element method for the acoustic analysis of finite periodic structures, Computer Methods in Applied Mechanics and Engineering, 391 (2022) 114528, https://doi.org/10.1016/j.cma.2021.114528.
[12] A. Brancati, M.H. Aliabadi, I. Benedetti, Hierarchical Adaptive Cross Approximation GMRES Technique for Solution of Acoustic Problems Using the Boundary Element Method, CMES-Computer Modeling in Engineering & Sciences, 43 (2009) 149–172, https://doi.org/10.3970/cmes.2009.043.149.
[13] R.Y. Li, Y.J. Liu, W.J. Ye, A fast direct boundary element method for 3D acoustic problems based on hierarchical matrices, Engineering Analysis with Boundary Elements, 147 (2023) 171–180, https://doi.org/10.1016/j.enganabound.2022.11.035.
[14] X. Xie, Y.J. Liu, An adaptive model order reduction method for boundary element-based multi-frequency acoustic wave problems, Computer Methods in Applied Mechanics and Engineering, 373 (2021) 113532, https://doi.org/10.1016/j.cma.2020.113532.
[15] J. Meng, A. Boag, V. Lomakin, E. Michielssen, A multilevel Cartesian non-uniform grid time domain algorithm, Journal of Computational Physics, 229 (2010) 8430–8444, https://doi.org/10.1016/j.jcp.2010.07.026.
[16] A.A. Ergin, B. Shanker, E. Michielssen, Fast evaluation of three-dimensional transient wave fields using diagonal translation operators, Journal of Computational Physics, 146 (1998) 157–180, https://doi.org/10.1006/jcph.1998.5908.
[17] A. Boag, V. Lomakin, E. Michielssen, Nonuniform grid time domain (NGTD) algorithm for fast evaluation of transient wave fields, IEEE Transactions on Antennas and Propagation, 54 (2006) 1943–1951, https://doi.org/10.1109/Tap.2006.877186.
[18] A.A. Ergin, B. Shanker, E. Michielssen, Fast analysis of transient acoustic wave scattering from rigid bodies using the multilevel plane wave time domain algorithm, Journal of the Acoustical Society of America, 107 (2000) 1168–1178, https://doi.org/10.1121/1.428406.
[19] M. Vikram, B. Shanker, Fast evaluation of time domain fields in sub-wavelength source/observer distributions using accelerated Cartesian expansions (ACE), Journal of Computational Physics, 227 (2007) 1007–1023, https://doi.org/10.1016/j.jcp.2007.08.017.
[20] M. Carley, G. Ghorbaniasl, Fast computation of time-dependent acoustic fields, Journal of the Acoustical Society of America, 140 (2016) 3963–3970, https://doi.org/10.1121/1.4968018.
[21] L.X. Ying, G. Biros, D. Zorin, A kernel-independent adaptive fast multipole algorithm in two and three dimensions, Journal of Computational Physics, 196 (2004) 591–626, https://doi.org/10.1016/j.jcp.2003.11.021.
[22] W. Fong, E. Darve, The black-box fast multipole method, Journal of Computational Physics, 228 (2009) 8712–8725, https://doi.org/10.1016/j.jcp.2009.08.031.
[23] T. Takahashi, M. Tanigawa, N. Miyazawa, An enhancement of the fast time-domain boundary element method for the three-dimensional wave equation, Computer Physics Communications, 271 (2022) 108229, https://doi.org/10.1016/j.cpc.2021.108229.
[24] W. Hackbusch, A sparse matrix arithmetic based on $\mathscr{H}$-matrices: Part I: Introduction to $\mathscr{H}$-matrices, Computing, 62 (1999) 89–108, https://doi.org/10.1007/s006070050015.
[25] M. Bebendorf, Approximation of boundary element matrices, Numerische Mathematik, 86 (2000) 565–589, https://doi.org/10.1007/Pl00005410.
[26] M. Bebendorf, S. Rjasanow, Adaptive low-rank approximation of collocation matrices, Computing, 70 (2003) 1–24, https://doi.org/10.1007/s00607-002-1469-6.
[27] L. Grasedyck, W. Hackbusch, S. Le Borne, Adaptive geometrically balanced clustering of $\mathscr{H}$-matrices, Computing, 73 (2004) 1–23, https://doi.org/10.1007/s00607-004-0068-0.
[28] W. Hackbusch, S. Börm, Data-sparse approximation by adaptive $\mathscr{H}^2$-matrices, Computing, 69 (2002) 1–35, https://doi.org/10.1007/s00607-002-1450-4.
[29] H. Harbrecht, M. Peters, Comparison of fast boundary element methods on parametric surfaces, Computer Methods in Applied Mechanics and Engineering, 261 (2013) 39–55, https://doi.org/10.1016/j.cma.2013.03.022.
[30] S. Börm, J. Henningsen, $\mathscr{H}^2$-matrices for translation-invariant kernel functions, Engineering Analysis with Boundary Elements, 175 (2025) 106190, https://doi.org/10.1016/j.enganabound.2025.106190.
[31] J. Tibaut, M. Schanz, J. Ravnik, Fast Boundary-Domain Integral Method with the $\mathscr{H}^2$-matrix formulation for large scale numerical investigations, Engineering Analysis with Boundary Elements, 138 (2022) 1–12, https://doi.org/10.1016/j.enganabound.2022.01.019.
[32] S. Börm, L. Grasedyck, W. Hackbusch, Introduction to hierarchical matrices with applications, Engineering Analysis with Boundary Elements, 27 (2003) 405–422, https://doi.org/10.1016/S0955-7997(02)00152-2.

[33] D. Seibel, Boundary element methods for the wave equation based on hierarchical matrices and adaptive cross approximation, Numerische Mathematik, 150 (2022) 629–670, https://doi.org/10.1007/s00211-021-01259-8.
[34] A. Aimi, L. Desiderio, G. Di Credico, Partially pivoted ACA based acceleration of the energetic BEM for time-domain acoustic and elastic waves exterior problems, Computers & Mathematics with Applications, 119 (2022) 351–370, https://doi.org/10.1016/j.camwa.2022.05.024.
[35] R. Chen, S. Bin Sayed, N. Alharthi, D. Keyes, H. Bagci, An explicit marching-on-in-time scheme for solving the time domain Kirchhoff integral equation, Journal of the Acoustical Society of America, 146 (2019) 2068–2079, https://doi.org/10.1121/1.5125259.
[36] G.R. Cowper, Gaussian quadrature formulas for triangles, International Journal for Numerical Methods in Engineering, 7 (1973) 405–408, https://doi.org/10.1002/nme.1620070316.
[37] W.Y. Zhou, Y.Q. Chen, A mixed cell compressed sparse row for time domain boundary element method in elastodynamics, Advances in Engineering Software, 192 (2024) 103633, https://doi.org/10.1016/j.advengsoft.2024.103633.
[38] T. Saitoh, F. Chikazawa, S. Hirose, Convolution quadrature time-domain boundary element method for 2-D fluid-saturated porous media, Applied Mathematical Modelling, 38 (2014) 3724–3740, https://doi.org/10.1016/j.apm.2014.02.009.
[39] M. Schanz, V.L. Keshava, H. De Gersem, Comparison of $\mathscr{H}$-matrix- and FMM-based 3D-ACA for a time-domain boundary element method, Computational Mechanics, 78 (2026) 563–578, https://doi.org/10.1007/s00466-026-02768-9.
[40] W.Y. Zhou, Z.P. Wu, Y. Liu, Y. Chen, Time domain boundary element method for semi-infinite domain problems using CSR storage method, Engineering Analysis with Boundary Elements, 147 (2023) 267–275, https://doi.org/10.1016/j.enganabound.2022.11.027.
[41] J.V. Venås, T. Kvamsdal, Isogeometric boundary element method for acoustic scattering by a submarine, Computer Methods in Applied Mechanics and Engineering, 359 (2020) 112670, https://doi.org/10.1016/j.cma.2019.112670.